\documentclass[12pt]{article}
\usepackage{amssymb}
\usepackage{array}
\usepackage{mathrsfs}
\usepackage{amsfonts}
\usepackage{graphicx}
\usepackage{indentfirst}
\usepackage{colortbl,dcolumn}
\usepackage{color}
\usepackage{amsmath}
\allowdisplaybreaks[4]
\usepackage{psfrag}
\usepackage{booktabs}
\usepackage{cite}
\usepackage[hidelinks]{hyperref}
\usepackage{amsthm}
\allowdisplaybreaks
\numberwithin{equation}{section}
\usepackage[top=1in, bottom=1in, left=0.8in, right=0.8in]{geometry}
\renewenvironment{proof}%
  {\begin{trivlist}%
   \item[\hskip\labelsep]% 移除 "Proof." 标题
   \ignorespaces}%
  {\qed\end{trivlist}}% 在结束时添加 QED 符号

\usepackage{caption}
\usepackage{appendix} 
\usepackage{graphicx}
\usepackage{subcaption}
\usepackage{subcaption}
\usepackage{float}
\newcommand{\dd}{\text{d}}
\newtheorem{thm}{Theorem}[section]
\newtheorem{defn}[thm]{Definition}
\newtheorem{lem}[thm]{Lemma}
\newtheorem{prop}[thm]{Proposition}

\newtheorem{assumption}[thm]{Assumption}

\begin{document}

\title{Total variation error bounds for Schr\"odinger-F\"ollmer samplers
%Non-asymptotic error bounds for Schr\"odinger-F\"ollmer samplers in the total variation distance
%Efficient multimodal sampling via Schr\"odinger-F\"ollmer diffusion with temperatures
%High dimensional sampling based on the Schr\"odinger-F\"ollmer diffusion with temperatures: new efficient algorithm
%New efficient Schr\"odinger-F\"ollmer Sampler
\footnotemark[2] {}}
%%%
%%%%%%%%%%%%%%
	\author{
          \ 
	Xiaojie Wang,
    Xiaoyan Zhang
		\\
		\footnotesize School of Mathematics and Statistics,   HNP-LAMA, Central South University, Changsha 410083, Hunan, China
	}

\maketitle
\footnotetext[2]{This work was supported by Natural Science Foundation of China (12471394, 12371417, 12071488) and Hunan Basic Science Research Center for Mathematical Analysis (2024JC2002). 
       %(12071488, 12371417, 11971488)
       %(12471394, 12071488, 12371417). \\
                 E-mail addresses:
 x.j.wang7@csu.edu.cn, 
 x.y.zhang@csu.edu.cn.

 }

\begin{abstract}
Sampling from high-dimensional, multimodal target distributions plays a fundamental role in scientific computing, statistics and machine learning. 
Recently, the Schr\"odinger-F\"ollmer sampler (SFS) has been
%emerged 
proposed by Huang et al. (IEEE Trans. Inform. Theory, 2025) as an efficient sampling algorithm that transports the degenerate distribution at $t=0$ to the target distribution at $t=1$, without relying on long-time ergodicity. 
In this paper, we establish non-asymptotic error bounds for SFS in the total variation distance by deriving regularization estimates for the associated Kolmogorov equation via the Bismut-Elworthy-Li formula.
When the drift can be computed exactly, a convergence rate of order $\mathcal{O}( h|\ln h|)$  was obtained under certain smoothness conditions on the drift with a uniform time step-size $h>0$.
When the drift is approximated by Monte Carlo methods, we obtain a convergence rate of order
$\mathcal{O}( h|\ln h|+ 1/M)$, where $M$ is the number of samples used in the Monte Carlo approximation. 
This yields a sharper dependence on $M$ than the existing bound $\mathcal{O}(1/\sqrt{M})$ in the $L^2$-Wasserstein distance in the aforementioned paper, by turning the total variation error to the weak error analysis in the setting of non-smooth test functions.
Numerical experiments are presented to illustrate the convergence rate with respect to the time step-size $h$ and the Monte Carlo sample size $M$ and the effectiveness of the proposed algorithm.\\

%\newpage 
%Two new Schr\"odinger-F\"ollmer samplers (SFS) with temperatures are introduced and an enhanced convergence rate of order $\mathcal{O}(h)$ is obtained in both the $L^2$-Wasserstein and total variation distances under certain smoothness conditions on the drift. In particular, when the drift is approximated by a Monte Carlo method, weak error estimates against bounded test functions yield a Monte Carlo approximation error of order $\mathcal{O}(1/{M})$ in the total variation distance, improving upon the existing order $\mathcal{O}(1/\sqrt{M})$ in the $L^2$-Wasserstein distance.\\
%We introduce two new Schr\"odinger-F\"ollmer samplers (SFS) with temperatures and investigate their convergence rate in both the $L^2$-Wasserstein and total variation distances. Further, an enhanced convergence rate of order $\mathcal{O}(h)$ is obtained in both distances under certain smoothness conditions on the drift. In particular, when the drift is approximated by Monte Carlo method, weak error estimates against bounded test functions yield a Monte Carlo approximation error of order $\mathcal{O}(1/{M})$ in the total variation distance, improving upon the existing order $\mathcal{O}(1/\sqrt{M})$ in the $L^2$-Wasserstein distance.\\


\end{abstract}

\section{Introduction}
\noindent
Sampling from complex, high-dimensional multimodal distributions
of the form 
\begin{align}
\label{eq:target-distribution}
    \mu(\mathrm{\,d} x) 
    \propto 
     e^
    {-V(x)} 
    \mathrm{\,d} x,
    \quad
    x \in \mathbb{R}^d,
    \
    d \gg 1,
\end{align}
where 
$V:\mathbb{R}^d \rightarrow \mathbb{R}$ is a potential function, is increasingly relevant in a number of application
domains, ranging from Bayesian inference, statistical physics, to machine learning, computational biology among many others. 
This topic has recently attracted considerable attention from the research community, with many fundamental sampling algorithms developed, including Markov Chain Monte Carlo (MCMC) methods, Annealed Importance Sampling, among others.
%A large number of research literature on sampling methods have emerged in recent years.
%
%Especially, there is good deal of intensive research work on the Markov Chain Monte Carlo (MCMC) methods, including unadjusted Langevin Monte Carlo (LMC) algorithm\cite{DS2017,DAME2017}, the Stochastic Gradient Langevin Dynamics (SGLD)\cite{barkhagen2021stochastic,ZDS2023,li2022sharp}, the celebrated Metropolis Hastings (MH) algorithm\cite{haario2001adaptive,bou2013nonasymptotic}, the underdamped Langevin dynamics \cite{guillin2024error,schuh2024convergence} and the references therein.
%
A typical approach of sampling is overdamped Langevin Monte Carlo (LMC):
\begin{equation}
\label{eq:ULA}
    Y_{n+1}
    =
    Y_n
    -\nabla V(Y_n) h
    +
    \sqrt{2}
    \Delta W_n, 
    \quad Y_0=x_0,
\end{equation}
where $h>0$ is the uniform time step-size and $\Delta W_n := W_{t_{n+1}}- W_{t_{n}}$.
In practice, the LMC \eqref{eq:ULA} is an Euler discretization of the following Langevin stochastic differential equation (SDE):
\begin{equation}
\label{eq:unadjusted-Langevin-SDE}
    \mathrm{\, d} X_t
    =
    -\nabla V(X_t) 
    \mathrm{\, d} t
    +
    \sqrt{2} 
    \mathrm{\, d} W_t, 
    \quad X_0=x_0, 
    \quad t>0,
\end{equation}
where $-\nabla V(\cdot)$ is the drift term and $\left\{W_t\right\}_{t \geq 0}$ is a standard $d$-dimensional Brownian motion process. 
Under certain (weak) convexity type conditions, the Langevin diffusion \eqref{eq:unadjusted-Langevin-SDE} admits the target distribution \eqref{eq:target-distribution} as its invariant distribution.
Therefore, by taking a small step-size for the discretization and a large number of iterations, LMC \eqref{eq:ULA} will
yield an approximate sample from $\mu$.
%In accordance with the specified conditions, the Langevin diffusion process $\left\{X_t\right\}_{t \geq 0}$ in \eqref{eq:unadjusted-Langevin-SDE} admits an invariant distribution $\mu(x) = e^{ - V (x)} / C,\quad x\in \mathbb{R}^d$, where $C$ is the normalising constant \cite{bakry2008rate}. 
%To asymptotically sample from $\mu$, one classic example is the unadjusted Langevin Monte Carlo algorithm, which also corresponds to the well-known Euler-Maruyama scheme of the Langevin SDE \eqref{eq:unadjusted-Langevin-SDE}, denoted by $\left\{Y_n\right\}_{n \geq 0}$ as\begin{equation}
%\label{eq:ULA}
%    Y_{n+1}
%    =
%    Y_n
%    -\nabla U(Y_n) h
%    +
%    \sqrt{2 h} \xi_{n+1}, 
%    \quad Y_0=x_0,
%\end{equation}
%where $h \in(0,1)$ represents the uniform timestep and $\xi_k=\left(\xi_{1, k}, \ldots, \xi_{d, k}\right)^T, k \in \mathbb{N}$, are independent and identically distributed (i.i.d.) standard $d$-dimensional Gaussian vectors.

%
Non-asymptotic convergence analyses for the LMC algorithm \eqref{eq:ULA} have been well established in the literature under the assumptions that the function $ \nabla V(x)$ is globally Lipschitz continuous and that the potential function $V(x)$ is strongly convex, see \cite{CB2018,DS2017,DMM2019,DM2017,SZ2019}.
%Nevertheless, the strongly convex assumption is extremely tough for the gradient of function $V(x)$ to satisfy.
%
The fact that most real-world target distributions are not log-concave poses a significant challenge and motivates the development of more robust samplers for such complex problems.
As a result, a number of related studies has analyzed the convergence of LMC under weaker assumptions, including the contractivity at infinity condition \cite{PWW25,NNZ2025,liu2025non,MMS2020} and the log-Sobolev inequality condition \cite{MWF2022,Li2025Asharp,YW2025,CELAZ2025,LS2025} and so on. 
Such conditions can be viewed as weak convexity conditions, used to ensure the long-time dynamics such as ergodicity of the Langevin SDEs \eqref{eq:unadjusted-Langevin-SDE}.
Despite the availability of LMC algorithm, sampling from high-dimensional or multimodal target distributions remains challenging, due to the need for long Markov chains to ensure adequate mixing and convergence.
In fact, LMC \eqref{eq:ULA} is highly sensitive to its initial state as the sampled points are prone to becoming trapped in a single mode and struggle to escape (see \cite[subsection 5.2]{Wang2026Multimodal} for details).
%\\
%{\it Instead of the Langevin sampling over infinite time, can one do efficient sampling over unit interval $[0, 1]$ without ergodicity, to avoid the (weak) convexity condition?}

%{\it It is natural to ask whether one can further construct an efficient sampling algorithm over unit interval $[0, 1]$ that eliminates the need for long-time simulation. }\\

To address this issue, some researchers have made some efforts (see \cite{AMVE2023,zhang2025stochastic,Albergo2025Stochastic,MRjiao,DYJ2023,RHBA2023,Wang2026Multimodal,Jiao2026Convergence}).
Very recently, the authors of \cite{MRjiao} proposed a novel Schr\"odinger-F\"ollmer sampler (SFS) based on the following Schr\"odinger F\"ollmer diffusion
\begin{equation}
\label{eq:sfp}
    \mathrm{\,d} X_t
    =
    \nabla \log \big( Q_{1-t} g(X_t) \big) 
    \mathrm{\,d} t
    +
    \mathrm{\,d} W_t, 
    \quad X_0=0, 
    \quad t \in ( 0, 1 ],
\end{equation}
where $g$ denotes the Radon-Nikodym derivative of $\mu$ with respect to the $d$-dimensional standard 
    Gaussian distribution
$\gamma^d$,
%the $d$-mensional Gaussian distribution $N(0,\beta\,\mathbf{I}_d)$, 
%or the ratio of the density of $\mu$ over the density of $\gamma^d$, i.e., 
given by \eqref{eq:Radon-Nikodym-derivative-g} 
and $\{Q_t\}_{t \in [0, 1]}$ is the heat semigroup defined by \eqref{eq:semi-group} with $\beta=1$,
which transports the degenerate distribution $\delta_0$ at $t=0$ to the target distribution 
$
\mu(\mathrm{\,d} x) 
    \propto 
     e^
     {-V(x)} 
    \mathrm{\,d} x
$
at $t=1$.
Unlike the Langevin samplers discussed above, SFS is an efficient method that evolves over the fixed time interval $[0,1]$. Moreover, it does not rely on any weak convexity or dissipativity condition to guarantee ergodicity. Numerical experiments suggest that SFS can outperform vanilla Langevin samplers, especially for multimodal target distributions (see \cite[Subsection 5.2]{Wang2026Multimodal} for details).
However, the drift of \eqref{eq:sfp} does not always have closed forms, and may need to be approximated by Monte Carlo methods.
In this case, the sampling error of SFS comes not only from the time discretization, but also the Monte Carlo approximation of the drift term.
And a convergence rate of order $\mathcal{O}(\sqrt{h}+\tfrac{1}{\sqrt{M}})$ in the $L^2$-Wasserstein distance was obtained there for the Euler discretization with uniform step-size $h>0$ under the Lipschitz condition of the drift term and strong convexity condition of the potential function, where $M$ is the number of samples used in the Monte Carlo approximation.
Furthermore, the authors of \cite{Jiao2026Convergence} obtained the convergence rate of order $\mathcal{O}(\sqrt{h}+\tfrac{1}{\sqrt{M}})$ in the $L^2$-Wasserstein distance without strong convexity condition.
In \cite{Wang2026Multimodal}, the authors introduced a variant of Schr\"odinger-F\"ollmer diffusion with temperatures as follows:
\begin{align}
\label{eq:intro-sf-beta}
    \mathrm{\,d} X_t
    =
    \beta\,\nabla \log Q^\beta_{1-t} g_{\beta}(X_t)
    \mathrm{\,d} t
    +
    \sqrt{\beta}
    \mathrm{\,d} W_t, 
    \quad X_0=0, 
    \quad t \in ( 0, 1 ],
\end{align}
where $g_{\beta}$ is
given by \eqref{eq:Radon-Nikodym-derivative-g} 
and $\{Q^\beta_t\}_{t \in [0, 1]}$ is the heat semigroup defined by \eqref{eq:semi-group}.
By carefully handling the singularity of the time derivative of the drift, the authors proposed a novel approach of error analysis to obtain an enhanced convergence rate of order $\mathcal{O}(h+ \tfrac{1}{\sqrt{M}})$ in the $L^2$-Wasserstein distance for the time discretization when the drift term $f_\beta$ lacks a closed form expression.
To speed up the sampling, the authors of \cite{Lin2025SRK} proposed an efficient stochastic Runge-Kutta Schr\"odinger-F\"ollmer sampler (SRKSFS) and proved that it achieves a convergence rate of order $\mathcal{O} ( h^{3/2} |\ln h|+\tfrac{1}{\sqrt{M}}) $ in the $L^2$-Wasserstein distance, considerably improving the order $\mathcal{O}(h)$ of the existing Euler type sampler. 
%However, all analyses for SFS algorithms work in the $L^2$-Wasserstein distance. 
%
In \cite{endo2024weak}, the authors yields convergence for SFS with exact drift in total variation under suitable assumptions, but does not provide an explicit non-asymptotic convergence rate. 
To the best of our knowledge, convergence rates for SFS in the total variation distance, have not yet been established.
%This may be challenging, as it seems different from our proof technique.
Bearing the above error bounds in mind, the following interesting question thus arises:\\

{\it Can an almost first-order error bound still be established \textbf{in the total variation distance} under suitable conditions?
}\\

In this work, we are concerned with the non-asymptotic convergence rate for the Schr\"odinger-F\"ollmer sampler to the target distribution \eqref{eq:target-distribution} in the total variation distance under suitable conditions (see Assumption \ref{ass:higher-g-lipschitz}). 
The key idea of this paper is to convert non-asymptotic bounds in the total variation distance between the law of the numerical approximation and the target distribution into weak error estimates for bounded, possibly non-smooth, test functions $\phi\in B_b(\mathbb{R}^d)$.
It is worth noting that our analysis is substantially more delicate than those in related works, where sufficiently smooth test functions are considered \cite{Pang2024Linear,wang2024weak,brehier2023approximation,fang2020adaptive}.
To overcome difficulties caused by the low regularity of the test functions, we rely on the Bismut-Elworthy-Li (BEL) formula to explore more delicate regularity estimates of the associated Kolmogorov equation (see Subsection \ref{subsec:Kolmogorov-equation} for details).
Combining these estimates with a careful treatment of the singular time derivative of the drift, we obtained a convergence rate of order $\mathcal{O}(d h|\ln h |)$ for the time discretization in the total variation distance in the setting the drift has a closed form. 
When the drift can not be calculated exactly, we used inexact drift based on Monte Carlo approximations and derived an error bound of order $\mathcal{O}(d h |\ln h|) + \mathcal{O}(\tfrac{d}{{M}})$ in the total variation distance, where $M$ is the number of samples used in the Monte Carlo approximation. 
To control the Monte Carlo approximation error, we use the unbiasedness of the numerator and denominator estimators.
After applying a Taylor expansion to the drift, the first-order terms vanish under conditional expectation, and the remaining second-order terms yield a convergence rate of order $\mathcal{O}(1/M)$ (see Section \ref{section:proof-theorem} for details), significantly improving the existing Monte Carlo convergence rate $\mathcal{O}(1/\sqrt{M})$ obtained in the $L^2$-Wasserstein distance in \cite{MRjiao,Jiao2026Convergence,Wang2026Multimodal}.
To conclude, our contribution can be summarized as follows:
\begin{itemize}

 \item When the drift can be exactly calculated, the following error bound in the total variation distance is established:
    \begin{align*}
        \big\|
        \mathcal{L}aw
        \left(Y_{1}\right)- \mu
        \big\|_{\mathrm{TV}}
        \leq 
        C d h |\ln h|, 
    \end{align*}
    where the constant $C$ is independent of $d$ and $h$, depending only on $L_g,\rho$ and $\beta$.
 
 \item  In the case that the drift can not be calculated exactly, the following error bound in the total variation distance is obtained:
 \[
    \big\|
    \mathcal{L}aw
    (\tilde{Y}_1)-\mu
    \big\|_{\mathrm{TV}}
    \leq 
    C 
    d
    h|\ln h|
    +
    \tfrac{C d}{M},
\]
  where $M$ is the number of samples used in the Monte Carlo estimator of the drift term, and the constant $C$ is independent of $d,h$ and $M$, depending only on $L_g,\rho$ and $\beta$. Moreover, the error caused by the Monte Carlo approximation is of order $\mathcal{O}(1/M)$ in total variation distance, yielding a sharper dependence on $M$ than the $\mathcal{O}(1/\sqrt{M})$ bounds in the $L^2$-Wasserstein distance established in \cite{MRjiao,Jiao2026Convergence,Wang2026Multimodal}.

  \item A variety of numerical experiments are reported to validate the effectiveness of SFS for sampling from target distributions. In particular, increasing the temperature can lead to a noticeable improvement when sampling from multimodal distributions. This phenomenon is intuitively explained by the fact that Markov chains at higher temperatures are able to cross energy barriers more easily (see Subsection \ref{subsec:Mog} for details).
 
\end{itemize}

The rest of this article is organized as follows. The next section formulates the simple background of Schr\"odinger-F\"ollmer sampler.
In Section \ref{sec:sfs}, we present assumptions and main results of this paper, whose proofs are given in Section \ref{section:proof-theorem}. 
In Section \ref{section:Preliminary-results}, we present some preliminary results. 
Some numerical tests are shown to illustrate our theoretical findings in Section \ref{section:Numerical-experiment}. 
Concluding remarks are given in the last section.
Finally, the appendix contains the detailed proof of several auxiliary lemmas and propositions.

\section{Schr\"odinger-F\"ollmer diffusion with temperatures}
\subsection{Notation}
Throughout this paper, we use $\mathbb{N}$ to denote the set of all positive integers and let $d \in \mathbb{N}$.
For short, we denote $[m]:= \{1,...,m\}$ and $[m]_0:= \{0, 1,...,m\}$.
Let $\|\cdot\|$ and $\langle\cdot, \cdot\rangle$ denote the Euclidean norm and the inner product of vectors in $\mathbb{R}^d$, respectively. 
Adopting the same notation as the vector norm, we denote $\|A\|:=\|A\|_{\mathrm{op}}=\sup_{\|v\|=1}\|Av\|$ as the operator norm of a matrix $A \in \mathbb{R}^{d \times d}$.
We use $\mathbf{0}_d,\mathbf{1}_d \in \mathbb{R}^d$ to denote the zero vector and the all-ones vector, respectively.
Let $\mathbf{I}_d \in \mathbb{R}^{d\times d}$ denote the identity matrix.
%For any $r\in \mathbb{R}^+$, let $\lfloor r\rfloor$ denote the greatest integer less than or equal to $r$.
%
We denote the $d$-dimensional standard Gaussian distribution by $\gamma^d$.
For two numbers $a$ and $b$, we denote $a \vee b:=\max\{a, b\}$ and $a \wedge b:=\min\{a, b\}$. 
We write $\tilde{\mathcal{O}}(\cdot)$ to mean that it hides polylogarithmic factors in $d$ and $1/\epsilon$.

Let $\{W_t\}_{t \in [0, 1] }$ be a standard $d$-dimensional Brownian motion process, defined on a filtered probability space $
\big(
\Omega_W, 
\mathcal{F}^W,
%\left\{\mathcal{F}^W_t\right\}_{t \in[0, 1]}, 
\mathbb{P}_W, 
\{\mathcal{F}^W_t\}_{t \in[0,1]}
\big)
$ 
satisfying the usual conditions.
Let $ \pi_h $ be a temporal grid of the form 
$$
\pi_h:=\{t_n=n \cdot h, 
\quad n\in [N]_0, \text { with } h=1 / N\}.
$$
Let $\{\xi_{t_{n,j}}\}_{ n \in[N]_0, \,j \in \mathbb{N}}$ be an i.i.d. family of standard Gaussian distributed random variables, independent of $\left\{W_t\right\}_{t \in [0, 1]}$,  defined on an additional filtered probability space 
$
\big(
\Omega_\xi, 
\mathcal{F}^\xi,
\mathbb{P}_\xi,
\{\mathcal{F}^\xi_{t}\}_{t \in[0,1]}\big)$,
where $\mathcal{F}_t^\xi$ is the $\sigma$ algebra generated by $\bigl\{\xi_{t_{k,j}}:t_k\le t,\ j\in\mathbb N\bigr\}$. 
By $\mathbb{E}_W$ and $\mathbb{E}_\xi$, we denote expectations in these two probability spaces 
$
\big(
\Omega_W, 
\mathcal{F}^W,
%\left\{\mathcal{F}^W_t\right\}_{t \in[0, 1]}, 
\mathbb{P}_W
\big)
$
and
$
\left(
\Omega_\xi, 
\mathcal{F}^\xi,
%\left(\mathcal{F}_n^\xi\right)_{n \in \mathbb{N}}, 
\mathbb{P}_\xi\right),
$
respectively.
By introducing the following product probability space
\begin{align*}       (\Omega,\mathcal{F},\mathbb{P},\{\mathcal{F}_t\}_{t\in[0,1]})
:=
\big(
\Omega_W\otimes\Omega_{\xi}, \mathcal{F}^W
\otimes
\mathcal{F}^{\xi},
\mathbb{P}_W
\otimes\mathbb{P}_{\xi},
\mathcal{F}^W_t
\otimes
\mathcal{F}^\xi_t
\big),
\end{align*}
%Moreover, for each temporal grid $\pi_h$ a discrete-time filtration $\left(\mathcal{F}_n\right)_{n \in[N]_0}$ is given by $$\mathcal{F}_n:=\mathcal{F}_{t_n}^W \otimes \mathcal{F}_{t_n}^\xi, \quad \text { for } n \in[N].$$
%
%Let $\xi_{1,n},\cdots,\xi_{M,n}$ be i.i.d. $\gamma^d,$ where $M \geq 1$ is sufficiently large.
%
we use $\mathbb{E} $ to mean the expectation in the product probability space and $L^r\left(\Omega, \mathbb{R}^d\right), r \geq 1$, to denote the family of $\mathbb{R}^d$-valued random variables $\eta$ satisfying $\mathbb{E} \left[\|\eta\|^r\right] := \mathbb{E}_W \big( \mathbb{E}_\xi \left[\|\eta\|^r\right] \big) <\infty$. 

Let $C_b(\mathbb{R}^d, \mathbb{R})$ (resp. $B_b(\mathbb{R}^d)$) denote the Banach space of all continuous and bounded mappings (resp. Borel bounded mappings) $ \varphi:\mathbb{R}^d \rightarrow \mathbb{R}$ endowed
with the norm $|\varphi|_0:=\sup _{x \in\mathbb{R}^d} |\varphi(x)|$.
%In addition, denote by $\mathrm{Lip}(\mathbb{R}^d)$ the Banach spaces of all Lipschitz continuous mappings $\phi$ endowed with the norm %$|\phi|_{\mathrm{\rm Lip}}
%:=
%\sup_{x\ne y}|\phi(x)-\phi(y)|
%\|x-y\|^{-1}$.
%
For $k \geq 0$, denote by $C^k\left(\mathbb{R}^d, \mathbb{R}\right)$ the set of functions from $\mathbb{R}^d$ to $\mathbb{R}$ which has continuous 0 -th, $\ldots, k$-th order derivatives, further denote by $C_b^k\left(\mathbb{R}^d, \mathbb{R}\right)$ the set of functions from $\mathbb{R}^d$ to $\mathbb{R}$ which has bounded continuous $0$-th, $\ldots, k$-th order derivatives.
For $y \in C^3\left(\mathbb{R}^d, \mathbb{R}\right)$ and $v_1, v_2, v_3, x \in \mathbb{R}^d$, we denote
\begin{align*}
        \nabla_{v_1} y(x) 
        & =
        \lim_{\varepsilon \rightarrow 0} \frac{y(x+\varepsilon v_1)
        -y(x)}
        {\varepsilon}, \\
        \nabla_{v_2} 
        \nabla_{v_1} y(x) 
        & =
        \lim_{\varepsilon \rightarrow 0} \frac{\nabla_{v_1} y\left(x+\varepsilon v_2\right)
        -
        \nabla_{v_1} y(x)}{\varepsilon},\\
        \nabla_{v_3} 
        \nabla_{v_2} 
        \nabla_{v_1} y(x)
        & =
        \lim_{\varepsilon \rightarrow 0} \frac{\nabla_{v_2} \nabla_{v_1} y\left(x+\varepsilon v_3\right)
        -\nabla_{v_2} 
        \nabla_{v_1} y(x)}{\varepsilon},
\end{align*}
as the directional derivatives of $y$. 
If the function $y$ is differentiable at the point 
$x \in \mathbb{R}^d$, then the directional derivative exists along any nonzero vector $v\in \mathbb{R}^d$. 
In this case, we have
\begin{align*}
    \nabla_{v}y(x)=
    \left \langle
    \nabla y(x), v
    \right\rangle.
\end{align*}
One knows 
$\nabla y(x) \in \mathbb{R}^d, \nabla^2 y(x) \in \mathbb{R}^{d \times d}, 
\nabla^3 y(x) \in$ $\mathbb{R}^{d \times d \times d}$. 
Moreover, we define the operator norm of $\nabla^k y(x), k=1,2,3$ by
\begin{align}
\label{SFS-W2-eq:operator-norm-def}
    \big\|
    \nabla^k y(x)
    \big\|
    :=
    \big\|
    \nabla^k y(x)
    \big\|_{\mathrm{op}}
    =
    \sup _{\|v_i\|=1, i=1, \ldots, k}
    \big\|
    \nabla_{v_k} 
    \ldots 
    \nabla_{v_1} y(x)
    \big\|.
\end{align} 
For the vector-valued function $\mathbf{u}: \mathbb{R}^d \rightarrow \mathbb{R}^{\ell}, \mathbf{u}=\big(u_{(1)}, \ldots, u_{(\ell)}\big)^T$, its first order partial derivative is considered as the Jacobian matrix as
\begin{align*}
    D \mathbf{u}=\left(\begin{array}{ccc}
\frac{\partial u_{(1)}}{\partial x_1} & \cdots & \frac{\partial u_{(1)}}{\partial x_d} \\
\vdots & \ddots & \vdots \\
\frac{\partial u_{(\ell)}}{\partial x_1} & \cdots & \frac{\partial u_{(\ell)}}{\partial x_d}
\end{array}\right)_{\ell \times d} .
\end{align*}
For any $v_1 \in \mathbb{R}^d$, one knows $D(\mathbf{u}) v_1 \in \mathbb{R}^{\ell}$ and one can define $D^2 \mathbf{u}\left(v_1, v_2\right)$ as
\begin{align*}
    D^2 \mathbf{u}\left(v_1, v_2\right):=D\left(D(\mathbf{u}) v_1\right) v_2, \quad \forall v_1, v_2 \in \mathbb{R}^d .
\end{align*}
In the same manner, one can define
\begin{align*}
    D^3 \mathbf{u}\left(v_1, v_2, v_3\right):=D\left(D\left(D(\mathbf{u}) v_1\right) v_2\right) v_3, \quad \forall v_1, v_2, v_3 \in \mathbb{R}^d,
\end{align*}
and for any integer $k \geq 3$ the $k$-th order derivatives of the function $\mathbf{u}$ can be defined recursively.
We remark that the partial derivatives of the scalar valued function can be covered by the special case $\ell=1$.

%%%%%%%%%%%%%%%%%%
%
\iffalse
Given the Banach spaces $\mathcal{X}$ and $\mathcal{Y}$, we denote by $L(\mathcal{X}, \mathcal{Y})$ the Banach space of bounded linear operators from $\mathcal{X}$ into $\mathcal{Y}$. 
%
Then the partial derivatives of the function $\mathbf{u}$ can also be regarded as the operators
\begin{equation}
    \begin{gathered}
         D \mathbf{u}(\cdot)(\cdot): \mathbb{R}^d \rightarrow L\left(\mathbb{R}^d, \mathbb{R}^{\ell}\right), \\
         D^2 \mathbf{u}(\cdot)(\cdot, \cdot): \mathbb{R}^d \rightarrow L\left(\mathbb{R}^d, L\left(\mathbb{R}^d, \mathbb{R}^{\ell}\right)\right) \cong L\left(\mathbb{R}^d \otimes \mathbb{R}^d, \mathbb{R}^{\ell}\right)
    \end{gathered}
\end{equation}
and
\begin{equation}
    D^3 \mathbf{u}(\cdot)(\cdot, \cdot \cdot \cdot): \mathbb{R}^d \rightarrow L\left(\mathbb{R}^d, L\left(\mathbb{R}^d, L\left(\mathbb{R}^d, \mathbb{R}^{\ell}\right)\right)\right) \cong L\left(\left(\mathbb{R}^d\right)^{\otimes 3}, \mathbb{R}^{\ell}\right) .
\end{equation}
%
For any $k \in \mathbb{N}$, let $C_b^k\left(\mathbb{R}^d\right)$ be the subspace of $C_b\left(\mathbb{R}^d\right)$ consisting of all functions with bounded partial derivatives $D^i \varphi(x), 1 \leq i \leq k$. In what follows, we use the letter $C$ to denote generic constants, independent of both the step-size $h \in(0,1)$ and the dimension $d$. Also, let $\mathbf{1}_H$ be the indicator function of a set $H$.
\fi
%%%%%%%%%%%%%%%%

%
Let $\mathcal{B}\left(\mathbb{R}^d\right)$ be the Borel $\sigma$-field of $\mathbb{R}^d$ and $\mathcal{P}\left(\mathbb{R}^d\right)$ be the space of all probability distributions on $\left(\mathbb{R}^d, \mathcal{B}\left(\mathbb{R}^d\right)\right)$.
Let $\nu_1$ and $\nu_2$ be two probability measures defined on $\left(\mathbb{R}^d, \mathcal{B}\left(\mathbb{R}^d\right)\right)$. 
%and let $\mathcal{D}\left(\nu_1, \nu_2\right)$ represent the collection of couplings $\nu$ on $\left(\mathbb{R}^{2 d}, \mathcal{B}\left(\mathbb{R}^{2 d}\right)\right)$ whose first and second marginal distributions are $\nu_1$ and $\nu_2$, respectively.
%
Further, the total variation distance between two Borel probability distributions $\nu_1$ and $\nu_2$ is defined by
\begin{equation}
\label{eq:tv-definition}
\big\|
\nu_1-\nu_2
\big\|_{\mathrm{TV}}
:=
\sup_{\substack{\varphi\in B_b(\mathbb{R}^d), \varphi\neq 0}}
\frac{
\left|
\int_{\mathbb{R}^d}\varphi(x)\nu_1(\mathrm{d}x)
-
\int_{\mathbb{R}^d}\varphi(x)\nu_2(\mathrm{d}x)
\right|
}
{
|\varphi|_0
}.
\end{equation}
Recalling that for any 
$\varphi \in B_b\left(\mathbb{R}^d\right)$, 
there exists a sequence 
$\{\varphi_k \}_{k \in \mathbb{N}}$ of bounded and continuous functions which converges boundedly and pointwise to $\varphi$. Namely, 
$\sup_{k \in \mathbb{N}} \sup_{x \in \mathbb{R}^d} |\varphi_k(x)| < \infty$ and 
$\lim _{k \rightarrow \infty} \varphi_k(x) \rightarrow \varphi(x)$, for all 
$x \in \mathbb{R}^d$ (see \cite{Brehier2024Total} or \cite[Proposition 4.2]{ethier2009markov} for details). 
Therefore, the total variation distance between two Borel probability distributions $\nu_1$ and $\nu_2$ has the representation as follows
\begin{align}
    \big\|
\nu_1-\nu_2
\big\|_{\mathrm{TV}}
& =
    \sup_{\varphi \in C_b(\mathbb{R}^d),\varphi \neq 0}
    \frac{
    \left|
    \int_{\mathbb{R}^d} \varphi(x)
    \nu_1(\mathrm{\, d} x)
    -
    \int_{\mathbb{R}^d} \varphi(x) 
    \nu_2(\mathrm{\, d} x)
    \right|
    }
    {|\varphi|_0}
    \notag\\
    & =
    \sup_{\varphi \in C_b(\mathbb{R}^d),|\varphi|_0 \leq 1}
    \left|
    \int_{\mathbb{R}^d} \varphi(x)
    \nu_1(\mathrm{\, d} x)
    -
    \int_{\mathbb{R}^d} \varphi(x) 
    \nu_2(\mathrm{\, d} x)
    \right|.
\end{align}
%%and the $\mathcal{W}_1$-distance between two Borel probability distributions $\mu_1$ and $\mu_2$ has the representation as follows,
%%\begin{align}
%    \mathcal{W}_1 \left(\mu_1,\mu_2\right)
%    & : =
%    \sup_{\varphi \in C(\mathbb{R}^d),|\varphi|_{Lip} \leq 1}
%    \left|
%    \int_{\mathbb{R}^d} \varphi(x) 
%    \mu_1(\mathrm{\, d} x)
%%%    -
%%    \int_{\mathbb{R}^d} \varphi(x) 
%    \mu_2(\mathrm{\, d} x)
%%    \right|.
%\end{align}
Consequently,  if $X$ is a $\mathbb{R}^d$-valued random variable, the probability distribution of $X$ is denoted by $\mathcal{L}aw(X)$.
Moreover, for the probability distribution of $X$ and $Y$, one also has
\begin{align}
\label{eq:definition-tv-rv}
        \big\|
\mathcal{L}aw(X)-\mathcal{L}aw(Y)
\big\|_{\mathrm{TV}}
        & :=
        \sup _{\varphi \in C_b(\mathbb{R}^d),|\varphi|_0 \leq 1}
        \big|
        \mathbb{E}
        [\varphi(X)]
        -
        \mathbb{E}
        [\varphi(Y)]
        \big|.
%        \nonumber\\
%        \mathcal{W}_1 \big(
%        \mathcal{L}aw(X),
%        \mathcal{L}aw(Y)
%        \big) 
%        & :=
%       \sup _{\varphi \in \mathrm{Lip}(\mathbb{R}^d),|\varphi|_{\mathrm{Lip}}
%        \leq 1}
%        \big|
%        \mathbb{E}
%        [\varphi(X)]
%        -
%        \mathbb{E}
%        [\varphi(Y)]
%        \big|.
\end{align}

\subsection{The Schr\"odinger-F\"ollmer process (SFP) with temperatures}
In this subsection, we revisit the Schr\"odinger-F\"ollmer diffusion process with temperatures in the literature.
We first make the following assumption.

\begin{assumption}
\label{ass:abs-continuous}
    Let $\mu \in \mathcal{P}(\mathbb{R}^d)$ denote the target distribution of interest, which is absolutely continuous with respect to the $d$-dimensional 
    %standard 
    Gaussian distribution denoted by $\mathcal{N}(0,\beta\,\mathbf{I}_d)$. 
\end{assumption}
Recall the Schr{\"o}dinger-F{\"o}llmer process with temperatures introduced in the latest  work \cite{Wang2026Multimodal}:
%$0<\beta<\infty$:
\begin{align}
    \label{eq:sde}
     \mathrm{\,d} X_t
     =
     f_{\beta}(X_t,t)
     \mathrm{\,d} t
     +
     \sqrt{\beta}
     \mathrm{\,d} W_t,
     \quad  X_0=0, 
     \quad t \in(0,1],
\end{align}
where $\beta \in (0, \infty)$ is a temperature parameter
and the drift $f_{\beta}$ is given by
\begin{align}
\label{eq:f-expression}
    f_{\beta}(x, t) 
    :=
    \beta~
    \nabla \log \big( Q^{\beta}_{1-t} g_{\beta}(x) \big)
    ,
    \quad  x \in \mathbb{R}^d, \: t \in[0,1],
%    \: \xi \sim \gamma^d.
\end{align}
with $g_{\beta}$ being the Radon-Nikodym derivative of $\mu$ with respect to $\mathcal{N}(0,\beta\,\mathbf{I}_d)$:
%or the ratio of the density of $\mu$ over the density of $N(0,\beta\,\mathbf{I}_d)$:
\begin{align}
\label{eq:Radon-Nikodym-derivative-g}
    g_{\beta}(x)
    :=
    \frac{\mathrm{d}\mu}
    {\mathrm{d} \mathcal{N} (0,\beta\,\mathbf{I}_d)}(x)
    =
    \tfrac{(2\pi\beta)^{d/2}}{C}
    \exp 
    \big(
    -V(x)
    +
    \tfrac{\|x\|^2}{2\beta}
    \big), 
    \quad x \in \mathbb{R}^d,
\end{align}
and $\{Q_t^{\beta}\}_{t \in [0, 1]}$ being the heat semigroup defined by
\begin{align}
\label{eq:semi-group}
    Q_t^{\beta} g_{\beta}(x)
    :=
    \mathbb{E}_{\xi}
    \big[
     g_{\beta}(x+\sqrt{t\beta} \,\xi)
    \big],
    \quad  t \in[0,1],
    \quad  \xi \sim \gamma^d.
\end{align}
%%%%%%%%%%%%%%%%%%%%%%%%%%%%%%%%%%%%%%%%
\iffalse
%where the normalizing constant $C > 0$ might be unknown.
%Noting that $g_{\beta}$ is a differentiable function satisfying  
$\mathbb{E}_{\xi}
\big[
\nabla g_{\beta}(x+\sqrt{(1-t)\beta}\,\xi)
\big]< +\infty$,
one can apply Stein's lemma \cite[Lemma 3.6.5]{garthwaite2002statistical} to obtain, for $t \in [0, 1)$,
\begin{align}
\label{eq:stein-lemma}
    \mathbb{E}_{\xi}
    \big[
    \nabla g_{\beta}(x+\sqrt{(1-t)\beta} \,\xi)
    \big]
    =
    \tfrac{1}{\sqrt{(1-t)\beta}}
    \,
    \mathbb{E}_{\xi}
    \big[
    \xi \,
    g_{\beta}\big(x+\sqrt{(1-t)\beta} \,\xi \big)
    \big],
\end{align}
which enables us to avoid the calculation of $\nabla g_{\beta}$ and 
%
recast the drift term as
\begin{align*}
    f_{\beta}(x, t)
    & =
    \frac{\beta
    \mathbb{\,E}_{\xi}
    [\nabla g_{\beta}(x+\sqrt{(1-t)\beta} \,\xi)]}{
    \mathbb{E}_{\xi}
    [g_{\beta}(x +\sqrt{(1-t)\beta} \,\xi)]}\\
    & =
    \frac{
    \beta\mathbb{\,E}_{\xi}
    \big[
    \xi \,g_{\beta}(x+\sqrt{(1-t)\beta} \,\xi)
    \big]}
    {\mathbb{E}_{\xi}
    \big[
    g_{\beta}(x+\sqrt{(1-t)\beta} \,\xi)
    \big]
    \cdot \sqrt{(1-t)\beta}},
    \quad \xi \sim \gamma^d,
\end{align*}
for $t \in [0, 1)$. 
\fi
%%%%%%%%%%%%%%%%%%%%%%%%%%%%%%%%%%%%%%%%%
In the special case $\beta
=1$, the Schr\"odinger-F\"ollmer process \eqref{eq:sde} with temperatures reduces into the usual Schr\"odinger-F\"ollmer process in \cite{MRjiao}.
%
%Moreover, we introduce a new notation $X_t^x$ for $t \in[0, 1]$ denoting the solution of SDE \eqref{SFS-W2-eq:f-sde} satisfying the initial condition $X_0^x=X_0=x$. 
%In this work, we set $X_0^x=X_0=x=0$, particularly.
%
To ensure that the $\operatorname{SDE}$ \eqref{eq:sde} admits a unique strong solution, we make the following assumption.
\begin{assumption}
    \label{ass:g-nabla-g-lip}
    Suppose that the Radon-Nikodym derivative $g_{\beta}$ and its gradient $\nabla g_{\beta}$ are $L_g$-Lipschitz continuous,
    and $g_{\beta}$ is uniformly bounded below by a positive constant $\rho$:
    \begin{align}
    \label{eq:g-lower-bounded}
        g_{\beta} \geq \rho >0.
    \end{align}
\end{assumption}

Under the above assumptions, we have the following proposition, which follows from \cite{lehec2013representation,MRjiao}.
 
\begin{prop}
    \label{prop:unique-strong-solution}
     Let Assumptions \ref{ass:abs-continuous}, \ref{ass:g-nabla-g-lip} hold. Then for any $t \in[0,1]$, the drift coefficient $f_{\beta}(\cdot,t):\mathbb{R}^d \rightarrow \mathbb{R}^d$ of SDE \eqref{eq:sde} is Lipschitz continuous and uniformly bounded. That is, there exist constants $L_f,\hat{L}_f > 0$, independent of $d$, such that, for any $x,y\in \mathbb{R}^d$ and $ t \in [0,1]$,
    \begin{align}
    \label{eq:f-lipschitz-condition}
        \big\|f_{\beta}(x,t)-f_{\beta}(y,t)\big\|
        \leq
        L_f\|x-y\|,
    \end{align}
    and
    \begin{align}
    \label{eq:f-linear-growth}
        \big\|f_{\beta}(x,t)\big\|
        & \leq
        \hat{L}_f,
    \end{align}
    where $L_f := \big(1+\tfrac{L_g}{\rho}
    \big)
    \tfrac{\beta L_g}{\rho},
    \,
    \hat{L}_f := \tfrac{\beta L_g}{\rho}$,
    and $L_g, \rho$ come from Assumption \ref{ass:g-nabla-g-lip}.
    Then the Schr\"odinger-F\"ollmer diffusion \eqref{eq:sde} has a unique strong solution $\{X_t\}_{t \in[0,1]}$ satisfying $X_1 \sim \mu$.
\end{prop}

Proposition \ref{prop:unique-strong-solution} shows that,  
%one can start from $X_0=0$ and update the values of $\{X_t: 0\leq $ $t \leq 1\}$ according to the Sch\"olinger-F\"ollmer SDE \eqref{SFS-W2-eq:f-sde} in continuous time, then the value $X_1$ has the desired distributional property, that is, $X_1 \sim \mu$.  
%
the Schr\"odinger-F\"ollmer process with  temperatures, as defined in \eqref{eq:sde}, transports the degenerate distribution $\delta_0$ at $t=0$ to the target distribution 
$
\mu(\mathrm{\,d} x) 
    \propto 
     e^
     {-V(x)} 
    \mathrm{\,d} x
$
at $t=1$.
The proof can be obtained by adapting the argument presented in \cite[Appendix B]{MRjiao}.

To guarantee the Lipschitz continuity of the drift $f_\beta$, we impose the uniformly lower bound assumption \eqref{eq:g-lower-bounded} on $g_\beta$. Although this assumption may seem restrictive, it is standard in the literature; see, for example, \cite{MRjiao,Jiao2026Convergence,lehec2013representation,RHBA2023}. To alleviate this restriction, one may regularize the target distribution $\mu$ by mixing it with the reference Gaussian distribution $\mathcal N(0,\beta\, \mathbf I_d)$, following the idea in \cite{MRjiao}. More precisely, for $\eta\in(0,1)$, define
$\mu_\eta := (1-\eta)\mu+\eta\mathcal N(0,\beta\mathbf I_d)$.
The corresponding density ratio with respect to $\mathcal N(0,\beta \, \mathbf I_d)$ is then
$
g_{\beta,\eta}
:=
\frac{d\mu_\eta}{d\mathcal N(0,\beta\mathbf I_d)}
=
(1-\eta)g_\beta+\eta.
$
It follows immediately that $g_{\beta,\eta}\ge \eta>0,$ and hence the lower bound condition \eqref{eq:g-lower-bounded} is satisfied. Moreover, if $g_\beta$ and $\nabla g_\beta$ are Lipschitz continuous, then so are $g_{\beta,\eta}$ and $\nabla g_{\beta,\eta}$.

%One can prove it by following the same lines of \cite[Appendix B]{MRjiao}.
%In addition, the Schr\"odinger-F\"ollmer process with temperatures have been used to solve a certain optimal problem (see \cite{DYJ2023}). 
%From a theoretical point of view, the introduce of the temperature parameter $\beta$ does not make any difference. Nevertheless, as revealed by numerical results, raising the temperature is crucial and makes a significant difference in sampling from high-dimensional multimodal distributions (see subsection \ref{subsec:Mog} for details).
%
%{\color{red}The intuition behind this is that the Markov chains at higher temperatures can cross energy barriers more easily, and hence they act as a bridge between isolated modes \cite{QW2024}.} 
%This happens because 
%This interesting finding can be explained that  
%%%%%%%%%%%%%%%%%%%%%%%

\section{Schr\"odinger-F\"ollmer samplers with temperatures}
\label{sec:sfs}
\noindent
To sample from the target distribution,  it suffices to solve SDE \eqref{eq:sde} to attain $X_1 \sim \mu$.
A practical way is to  discretize the continuous-time process and obtain the approximate solutions.

\subsection{Schr\"odinger-F\"ollmer sampler with exact drift}
If the target distribution $\mu$ has an easy form, such as the following  Gaussian mixture distributions \cite{LLMW2021},
\begin{equation}
        \mu
        =
        \sum_{i=1}^\kappa \theta_i 
        \,
        \mathcal{N}(\alpha_i, \Sigma_i),
        \quad \sum_{i=1}^\kappa \theta_i=1 
        \text { and } 
        0 \leq \theta_i \leq 1, 
        %\quad i=1, \ldots, \kappa,
\end{equation}
where $\kappa \in \mathbb{N}$ is the number of mixture components, $\mathcal{N}(\alpha_i, \Sigma_i)$ is the $i$-th Gaussian component with mean $\alpha_i \in \mathbb{R}^d$ and covariance matrix $\Sigma_i \in \mathbb{R}^{d \times d}$.
Then, after some tedious calculations, the drift $f_{\beta}(x,t)$ can be
calculated explicitly (see \cite{Wang2026Multimodal,MRjiao} for details).
%Next, we gives the following two variants of the Schr\"odinger-F\"ollmer sampler.
%
In this case, we use the Euler-Maruyama discretization of \eqref{eq:sde} with a fixed step-size. 
Then the Euler-Maruyama discretization of \eqref{eq:sde} is given by
\begin{align}
    \label{eq:Euler-scheme-no-mc}
    Y_{t_{n+1}}
    =
    Y_{t_n}
    +
    h f_{\beta}(Y_{t_n}, t_n)
    +
    \sqrt{\beta}
    \Delta W_{n},
    \quad n\in [N-1]_0,
    \quad Y_0=0,
\end{align}
where 
$\Delta W_{n} 
:= 
W_{t_{n+1}}
-
W_{t_n}, 
\,n \in [N-1]_0, \,N \in \mathbb{N}$.
In order to establish the main result, we now proceed to make an assumption as follows. 

\begin{assumption}
\label{ass:higher-g-lipschitz}
   For $\beta \in (0, \infty)$, we assume $g_{\beta}$ 
%    defined by \eqref{SFS-W2-eq:Radon-Nikodym-derivative-g} 
    $\in {C}^4(\mathbb{R}^d,\mathbb{R})$ and there exist constants $L_g,\rho > 0$ such that
    $g_{\beta},\nabla g_{\beta},\nabla^2 g_{\beta}$ and $\nabla^3 g_{\beta}$ are $L_g$-Lipschitz continuous and \eqref{eq:g-lower-bounded} holds.
\end{assumption}

The next proposition gathers some useful regularity estimates for the drift $f_\beta $ that follow directly from Assumption \ref{ass:higher-g-lipschitz}.
The proof follows the same argument as that of \cite[Proposition 3.2]{Wang2026Multimodal} and is therefore omitted here.

\begin{prop}
\label{prop:f-in-C^3}
    Let Assumptions \ref{ass:abs-continuous} and \ref{ass:higher-g-lipschitz} hold. Then, for any $t \in [0,1]$, the drift coefficient $f_{\beta}(\cdot,t):\mathbb{R}^d
    \rightarrow \mathbb{R}^d$ of SDE \eqref{eq:sde} is three times continuously differentiable with bounded partial derivatives.
    There exists a constant $L_f > 0$, which is independent of $d,t$, for any $ x,v_1,v_2,v_3 \in \mathbb{R}^d$, such that
    \begin{align*}
    %\label{eq:f-c-3}
             \|Df_{\beta}(x,t)v_1\| 
             \leq & 
             L_f\|v_1\|,\\
             \|D^2f_{\beta}(x,t)(v_1,v_2)\| 
             \leq &
             L_f
             \|v_1\|
             \cdot\|v_2\|,
    \end{align*}
    and 
    \begin{align*}
        \|D^3f_{\beta}(x,t)(v_1,v_2,v_3)\| 
        \leq 
        L_f
        \|v_1\|
        \cdot\|v_2\|
        \cdot\|v_3\|.
    \end{align*}
    Moreover, for any $x \in \mathbb{R}^d$, the drift coefficient $f_{\beta}(x,\cdot): [0,1) \rightarrow \mathbb{R}^d$ of SDE \eqref{eq:sde} is continuously differentiable.
    There exists a constant $\widetilde{L}_f > 0$, which is independent of $d,t$, for any $t \in [0,1)$, such that
    \begin{align}
    \label{eq:f-t-1-2}
        \big\|
        \partial_t f_{\beta}(x, t)
        \big\|
        \leq 
        \widetilde{L}_f
        \,
        d^\frac{1}{2}
        \frac{1}
        {\sqrt{1-t}}.
    \end{align}
\end{prop}

%%%%%%%%%%%%%%%%%%%%%%%%%%%%%%%%%%%%%%%%
\iffalse
\begin{assumption}
\label{ass:f-in-C^3}
    For any $t \in [0,1]$, the drift coefficient $f_{\beta}(\cdot,t):\mathbb{R}^d
    \rightarrow \mathbb{R}^d$ of SDE \eqref{eq:sde} is three times continuously differentiable with bounded partial derivatives.
    There exists a constant $L_f > 0$, which is independent of $d,t$, for any $ x,v_1,v_2,v_3 \in \mathbb{R}^d$, such that
    \begin{align*}
    %\label{eq:f-c-3}
             \|Df_{\beta}(x,t)v_1\| 
             \leq & 
             L_f\|v_1\|,\\
             \|D^2f_{\beta}(x,t)(v_1,v_2)\| 
             \leq &
             L_f
             \|v_1\|
             \cdot\|v_2\|,
    \end{align*}
    and 
    \begin{align*}
        \|D^3f_{\beta}(x,t)(v_1,v_2,v_3)\| 
        \leq 
        L_f
        \|v_1\|
        \cdot\|v_2\|
        \cdot\|v_3\|.
    \end{align*}
    Moreover, for any $x \in \mathbb{R}^d$, the drift coefficient $f_{\beta}(x,\cdot): [0,1) \rightarrow \mathbb{R}^d$ of SDE \eqref{eq:sde} is continuously differentiable.
    There exists a constant $\widetilde{L}_f > 0$, which is independent of $d,t$, for any $t \in [0,1)$, such that
    \begin{align}
    \label{eq:f-t-1-2}
        \big\|
        \partial_t f_{\beta}(x, t)
        \big\|
        \leq 
        \widetilde{L}_f
        \,
        d^\frac{1}{2}
        \frac{1}
        {\sqrt{1-t}}.
    \end{align}
\end{assumption}
\fi
%%%%%%%%%%%%%%%%%%%%%%%%%%%%%%%%%%%%%%%

Note that Proposition \ref{prop:f-in-C^3} immediately implies that, for any $x,y\in \mathbb{R}^d$ and $t \in [0,1]$,
%\begin{align}
%        \|D^2f_{\beta}(x,t)(v_1,v_2)
%        -
%        D^2f_{\beta}(y,t)(v_1,v_2)\|
%        \leq &
%        L_f 
%        \|x-y\|
%        \cdot\|v_1\|
%        \cdot\|v_2\|,
%        \nonumber\\ 
        %\quad \forall x,y,v_1,v_2 \in \mathbb{R}^d,\\
%        \|Df_{\beta}(x,t)v_1
%        -
%        Df_{\beta}(y,t)v_1\|
%        \leq &
%        L_f 
%        \|x-y\|
%        \cdot\|v_1\|, 
%        %\quad \forall x,y,v_1 \in \mathbb{R}^d,\\
%\end{align}
%and 
$$
    \|f_{\beta}(x,t)-f_{\beta}(y,t)\|
    \leq
    L_f\|x-y\|,
    %\quad \forall x, y\in \mathbb{R}^d,
$$
and
\begin{align}
\label{eq:one-sided-Lipschitz}
    \big\langle x-y,
    f_{\beta}(x,t)-f_{\beta}(y,t)
    \big\rangle 
    &\leq 
    L_f\|x-y\|^2.
    %\quad \forall x\in \mathbb{R}^d,
\end{align}
%and the linear growth condition,
%\begin{align}
%\label{eq:b-linear-growth}
%    \|f_{\beta}(x,t)\| 
%    &\leq 
%    \hat{L}_f  
%    +
%    L_f\|x\|.
%    %\quad \forall x\in \mathbb{R}^d,
%\end{align}
Moreover, the drift term $f_\beta$ is Lipschitz continuous with respect to the spatial variable and $\frac{1}{2}$-H\"older continuous with respect to the time variable. Namely, there exists a constant
$\bar{L}_f>0$ such that, for any $x,y\in\mathbb{R}^d$ and $t,s\in[0,1]$,
\begin{align}
    \label{eq:b-lipschitz-in-x-houder-in-t}
     \|f_{\beta}(x,t)-f_{\beta}(y,s)\|
     \leq 
     \bar{L}_f
     \big(
     \|x-y\|
     +
     d^\frac{1}{2}|t-s|^\frac{1}{2}
     \big).
     %\quad  \forall x,y\in \mathbb{R}^d, 
     %\quad \forall t,s\in [0,1].
\end{align}
%It is noted that $L_f \text{ and }\hat{L}_f$ are generic finite constants with value that could change upon each appearance, but will not depend upon the dimension $d$ and the step-size $h$.
%In other words, Assumption \ref{ass:f-in-C^3} can imply Assumption \ref{ass:f-lipschitz}.
%
This can be proved under Assumptions \ref{ass:abs-continuous} and \ref{ass:g-nabla-g-lip} (see Remark III.1 and Appendix D of \cite{MRjiao}).
Thanks to the above assumptions, we are able to derive a non-asymptotic error bound in the total variation distance for the SFS algorithm with exact drift \eqref{eq:Euler-scheme-no-mc}.

\begin{thm}
[Main result: error bounds with exact drift in the total variation distance]
\label{thm:exact-f}
    Let $N \in\mathbb{N}$ and let Assumptions \ref{ass:abs-continuous} and \ref{ass:higher-g-lipschitz} hold.
    Let $\{{Y}_{t_n}\}_{n \in[N]_0}$ be defined as \eqref{eq:Euler-scheme-no-mc} using the uniform step-size $h=\tfrac{1}{N}$. Then there exists a constant $C$ independent of $d$ and $h$, such that,
\begin{align*}
%    & (a)
    \
    \big\|
    \mathcal{L}aw
    \big(
    {Y}_1
    \big)
    - 
    \mu
    \big\|_{\mathrm{TV}}
    \leq
    C d
    h|\mathrm{ln}h|.
%    & (b)\
%    \mathcal{W}_1
%    \big(
%    \mathcal{L}aw
%    \big(
%    \tilde{Y}_1\big), \mu\big)
%    \leq
%    C d^2 h
%    +
%    \frac{Cd}{M}.
\end{align*}
\end{thm}
See Section \ref{section:proof-theorem} for proof of this theorem. 
As a direct consequence, we obtain the following result on the number of iterations, whose proof is straightforward and therefore omitted.

\begin{prop}
    Let Assumptions \ref{ass:abs-continuous} and \ref{ass:higher-g-lipschitz} hold. To achieve a given precision level $\epsilon>0$ in the total variation distance, a required number of iterations of the SFS algorithm \eqref{eq:Euler-scheme-no-mc} is of order $\tilde{\mathcal{O}}\big(\tfrac{d}{\epsilon}\big)$.
\end{prop}

\subsection{ Schr\"odinger-F\"ollmer sampler with inexact drift}
In general, however, deriving an explicit expression for the drift term $f_{\beta}$ is intractable when the target distribution $\mu$ is complex.
Thus, we approximate it using Monte Carlo methods.
%For $n \in \mathbb{N}$, let $\{\xi_{j}\}_{j\in[m_n]}$ and $\{\hat{\xi}_{j}\}_{j\in[\hat{m}_n]}$ be two families of independent standard Gaussian distributed random variables on the probability space $ \left( \Omega_\xi, \mathcal{F}^\xi,
%\left(\mathcal{F}_n^\xi\right)_{n \in \mathbb{N}}, \mathbb{P}_\xi\right)$ with 
%
In the following, 
we can approximate the drift term $f_{\beta}(x,t_n)$ by $\tilde{f}_{\beta}:\Omega_{\xi} \times \mathbb{R}^d \times \pi_h \rightarrow \mathbb{R}^d$, for any $n \in [N-1]_0$,
\begin{align}
\label{eq:drift-est}
\tilde{f}_{\beta}
    (x,t_n)
    & =
    \frac{\frac{\beta }{M} 
    \sum_{j=1}^M
    \left[
    \nabla 
    g_{\beta}
    (x+\sqrt{(1-t_n)\beta} 
    \,
    {\xi}_{t_{n,j}})
    \right]}
    {\frac{1}{M} 
    \sum_{j=1}^M
    \left[
    g_{\beta}
    (x+\sqrt{(1-t_n)\beta} \,{\xi}_{t_{n,j}})
    \right]}
    \quad \xi_{t_{n,j}} \sim \gamma^d,
\end{align}
%or 
%\begin{equation}
%\tilde{f}_{\beta}
%    (x,t_n)
%     =
%    \frac{\frac{\beta }{M} 
%    \sum_{j=1}^M
%    \left[
%    \xi_j
%    \,
%    g_{\beta}(x+\sqrt{(1-t_n)\beta} \,{\xi}_{j})
%    \right]}
%    {\frac{1}{m} 
%    \sum_{j=1}^m
%    \left[
%    g_{\beta}
%    (x+\sqrt{(1-t_n)\beta} \,{\xi}_{j})
%    \right]
%    \sqrt{(1-t)\beta}}, 
%    \quad \xi_{j} \sim \gamma^d,
%\end{equation}
%where the second equality stands due to Stein's lemma.
where, for $M \in \mathbb{N}$, $\{\xi_{t_{n,j}}\}_{j\in [M]}$ 
is an independent sequence of independent standard Gaussian distributed random variables.
It should
be noted that the Gaussian random variables are updated at each simulation time.
Then, the Euler-Maruyama discretization of \eqref{eq:sde} can be written as,
\begin{equation}
\label{eq:e-m-mc}
   \tilde{Y}_{t_{n+1}}
   =
   \tilde{Y}_{t_n}
   +
   h \tilde{f}_{\beta}
   (\tilde{Y}_{t_n}, t_n)
   +
   \sqrt{\beta} \Delta W_n, 
   \quad n\in[N-1]_0,
   \quad \tilde{Y}_0=0,
\end{equation}
where 
$\Delta W_{n} 
:= 
W_{t_{n+1}}
-
W_{t_n}, 
\,n \in [N-1]_0, \,N \in \mathbb{N}$.
The continuous-time version of Euler algorithm can be rewritten as follows,
\begin{equation}
\label{eq:Euler-Maruyama-contimuous-version}
     \tilde{Y}_s
     =
     \tilde{Y}_{t_n}
     + 
     (s-t_n) 
     \tilde{f}_{\beta}
     (\tilde{Y}_{t_n}, t_n)
     + 
     \sqrt{\beta}(W_s-W_{t_n}), 
     \quad n\in[N-1]_0,
     \quad \tilde{Y}_0=0,
\end{equation}
where $s\in [t_n,t_{n+1}] $. In the following proposition, we show that the inexact drift $\tilde{f}_{\beta}$ is also continuously differentiable and uniformly bounded.

\begin{prop}
\label{prop:tilde-Df-bound}
Let Assumptions \ref{ass:abs-continuous} and \ref{ass:higher-g-lipschitz} hold. Then 
for any $n\in [N-1]_0$, the estimator $\tilde{f}_{\beta}
    (\cdot,t_n)
    %\Omega_{\xi}\times \mathbb{R}^d \times \pi_h 
    \rightarrow \mathbb{R}^d$ is continuously differentiable with bounded partial derivatives and uniformly bounded.
    There exists a constant $L_f, \hat{L}_f > 0$, which is independent of $d,t$, for any $ x,v_1\in \mathbb{R}^d$, such that
    \begin{align}
    \label{eq:tilde-Df-bound}
        \big\|
        D \tilde{f}_\beta(x,t_n)v_1
        \big\|
        \leq
        L_f\|v_1\|,
        %\quad
        %\forall \, x,y \in \mathbb{R}^d,
        %\quad \text{for} ~ n \in \{0,1,\cdots,N\}.
    \end{align} 
    and 
    \begin{align}
    \label{eq:tilde-f-bound}
        \big\|\tilde{f}_\beta(x,t_n)\big\|
        & \leq
        \hat{L}_f,
    \end{align}
    where $L_f := \big(1+\tfrac{L_g}{\rho}
    \big)
    \tfrac{\beta L_g}{\rho},
    \,
    \hat{L}_f := \tfrac{\beta L_g}{\rho}$,
    and $L_g, \rho$ come from Assumption \ref{ass:g-nabla-g-lip}.
\end{prop}
The proof of Proposition \ref{prop:tilde-Df-bound} is straightforward due to Assumptions \ref{ass:abs-continuous}, \ref{ass:higher-g-lipschitz} and an argument similar to that used in \cite[Proposition 3.2]{Wang2026Multimodal}. We therefore omit the details.
Based on the above assumptions, we are able to derive a non-asymptotic error bound in the total variation distance for the SFS algorithm with inexact drift \eqref{eq:e-m-mc}, whose proof is deferred to Section \ref{section:proof-theorem}.

\begin{thm}
[Main result: error bounds with inexact drift in the total variation distance]
\label{thm-main-rerult-w0-w1}
    Let $N \in\mathbb{N}$ and let Assumptions \ref{ass:abs-continuous} and \ref{ass:higher-g-lipschitz} hold.
    Let $\{\tilde{Y}_{t_n}\}_{n \in[N]_0}$ be defined as \eqref{eq:e-m-mc} using the uniform step-size $h=\tfrac{1}{N} \in (0,1/4]$. Then there exists a constant $C$ independent of $d,h$ and $M$, such that,
\begin{align*}
    \big\|
    \mathcal{L}aw
    \big(
    \tilde{Y}_1
    \big)
    - 
    \mu
    \big\|_{\mathrm{TV}}
    \leq
    C d h|\mathrm{ln}h|
    +
    \frac{Cd}{M},
%    & (b)\
%    \mathcal{W}_1
%    \big(
%    \mathcal{L}aw
%    \big(
%    \tilde{Y}_1\big), \mu\big)
%    \leq
%    C d^2 h
%    +
%    \frac{Cd}{M}.
\end{align*}
where $M$ is the number of samples used in the Monte Carlo estimator $\tilde{f}_\beta$.
\end{thm}

As a direct consequence, we obtain the following result on the number of iterations and the number of samples of $\xi$ used in the Monte Carlo estimator, whose proof can be found in Appendix \ref{Appendix:mixing-time}.

\begin{prop}
\label{prop-mixing-time}
    Let Assumptions \ref{ass:abs-continuous} and \ref{ass:higher-g-lipschitz} hold. To achieve a given precision level $\epsilon>0$ in the total variation distance, a required number of iterations of the SFS algorithm \eqref{eq:e-m-mc} is of order $\tilde{\mathcal{O}}\big(\tfrac{d}{\epsilon}\big)$, the number of Monte Carlo samples used at each time step is of order
    $\mathcal{O}\big(\tfrac{d}{\epsilon}\big)$, and the total number of Monte Carlo samples used is of order
    $\tilde{O}\big(\tfrac{d^2}{\epsilon^2}\big)$.
\end{prop}

Under the assumptions that the drift is Lipschitz continuous and that the potential $V(x)$ is strongly convex, the authors of \cite{MRjiao} derived a convergence rate of order $\mathcal{O}(\sqrt{dh})$ for the Euler discretization of the Schr\"odinger-F\"ollmer diffusion \eqref{eq:sde} with $\beta=1$ in the $L^2$-Wasserstein distance. 
Furthermore, the authors of \cite{Jiao2026Convergence} derived a convergence rate of order $\mathcal{O}(\sqrt{dh})$ for the Euler discretization without strongly convex of the potential $V(x)$.
Very recently, the authors of \cite{Wang2026Multimodal} proved an improved error bound $\mathcal{O}(dh)$ in the $L^2$-Wasserstein distance by developing a novel approach of error analysis for the time discretization.
%However, in order to attain the first order convergence, the price to pay is usually imposing an additional smoothness assumption on the drift $f_\beta$.
%
To the best of our knowledge, all analyses for Schr\"odinger-F\"ollmer sampler work in the $L^2$-Wasserstein distance. 
In this work, we establish an almost first-order convergence rate in the total variation distance.
This may be challenging, as it seems different from previous proof technique.

Regarding the error caused by the Monte Carlo approximation, we obtain a bound of order $\mathcal{O}(1/M)$ in the total variation distance, improving the existing order $\mathcal{O}(1/\sqrt{M})$ in the $L^2$-Wasserstein distance. 
We would like to emphasize that this improvement is specific to the weak error analysis and relies on the Taylor expansion argument (see Section \ref{section:proof-theorem} for details).

\section{Preliminary results}
\label{section:Preliminary-results}
\noindent
In this section, we aim to derive uniform moment estimates for the exact solution
of \eqref{eq:sde} and the SFS algorithm with exact drift \eqref{eq:e-m-mc}, along with uniform regularity of the associated
Kolmogorov equation.

\subsection{A priori estimates of the SFP and SFS}
%We start with the uniform moments lemma for $\left\{X_t\right\}_{0 \leq t \leq 1}$, defined by \eqref{eq:sde}.
Similar to \cite[Lemma A.3, Lemma A.4]{MRjiao}, the process $\{X_t\}_{ t \in [0, 1]}$ defined by \eqref{eq:sde} has uniform moments and $\frac{1}{2}$-H\"older continuity.

\begin{lem}
\label{lem:x_t-moments-bounded}
   Let Assumptions \ref{ass:abs-continuous} and \ref{ass:g-nabla-g-lip} hold. Then the solution of SDE  \eqref{eq:sde}, denoted by $\{X_t\}_{t\in [0,1]}$ satisfies, for any $ t \in [0,1]$ and $0 \leq t_1 \leq t_2 \leq 1$,
   \begin{align}
        \label{SFS-W2-eq:x_t-moments-bounded}
        \mathbb{E}
        {_W}
        \left[
        \left\|
        X_t
        \right\|^2
        \right] 
        \leq 
        { M_1} d,
    \end{align}
    and 
    \begin{align}
    \label{SFS-W2-eq:x-1/2-houlder-continuous}
         \mathbb{E}{ _W}
         \left[
         \left\|
         X_{t_2}-X_{t_1}
         \right\|^2
         \right] 
         \leq 
         { M_2}
         \,
         d(t_2-t_1),
     \end{align}
where 
$
M_1 := 2(2\hat{L}^2_f + \beta) \exp(4 L^2_f)
$
and 
$
M_2 := 8 L^2_f \exp(4 L^2_f) (2\hat{L}^2_f + \beta) + 4 \hat{L}^2_f + 2\beta,
$
with constants $L_f$ and $\hat{L}_f$ given by \eqref{eq:f-lipschitz-condition} and \eqref{eq:f-linear-growth}, respectively.
\end{lem}

Next, we present some regularity estimates of the process $\big\{\tilde{Y}_t\big\}_{t \in [0,1]} $ defined in \eqref{eq:Euler-Maruyama-contimuous-version}.

\begin{lem}
\label{lem:numerical-solution-monents-bounded}
   Let Assumptions \ref{ass:abs-continuous} and \ref{ass:g-nabla-g-lip} hold. Let $\{\tilde{Y}_s\}_{s \in[t_n,t_{n+1}]}$ be defined as \eqref{eq:Euler-Maruyama-contimuous-version}, $n\in[N-1]_0$, where $N \in \mathbb{N}, \text{ and } h=\tfrac{1}{N}.$
   Then, the following estimates hold true,
\begin{equation}
    \label{eq:Y_s-bounded}
        \mathbb{E}
        \Big[
          \big\|
            \tilde{Y}_s
          \big\|^{2}
        \Big] 
        \leq
        C d,
\end{equation}
and
\begin{equation}
\label{eq:Ys-Yn}
    \begin{aligned}
       \mathbb{E}
        \Big[
          \big\|
         \tilde{Y}_s-\tilde{Y}_{t_n}
         \big\|^{2}
        \Big] 
        & \leq 
        C d h ,\\
       \mathbb{E}
       \Big[
          \big\|
       f_{\beta}(\tilde{Y}_s,s)
       -
       f_{\beta}(\tilde{Y}_{t_n},t_n)
       \big\|^{2}
        \Big] 
       & \leq 
       C d h,
  \end{aligned}
\end{equation}
where $C$ is independent of $d$ and $h$,  depending only on $L_f, \hat{L}_f$ and $\beta$.
\end{lem}

The proof of Lemma \ref{lem:numerical-solution-monents-bounded} is straightforward due to Assumption \ref{ass:g-nabla-g-lip}, \cite[Lemma A.7]{MRjiao} and the H\"older inequality. Similar estimates can be
found in \cite[Lemma 5.6]{Pang2024Linear} without the dependence with respect to the dimension $d$. As a result, we omit the proof here.

\subsection{Kolmogorov equation and regularization estimates}
\label{subsec:Kolmogorov-equation}
For any $\phi\in C_b(\mathbb R^d)$, let $v(s,x):=\mathbb E[\phi(X_1^{s,x})]$.
By the classical Kolmogorov backward equation for
SDE \eqref{eq:sde}, $v$ satisfies
\begin{equation}
\label{eq:kol-v}
\left\{\begin{array}{l}
\partial_s v(s,x)
+
Dv(s,x)f_\beta(x,s)
+
\tfrac{\beta}{2}
\sum_{j=1}^d
D^2v(s,x)(e_j,e_j)
=0,  
\quad  s \in [0,1),
\quad  x \in \mathbb{R}^d, 
\\
v(1,x)=\phi(x),
\quad  x \in \mathbb{R}^d, 
\end{array}\right.
\end{equation}
where $\{e_j\}_{j\in [d]}$ denotes the standard orthonormal basis of $\mathbb R^d$  (see \cite[p.~83]{mao2007stochastic} for details).
To carry out the weak error analysis, we define %function $u(\cdot, \cdot):[0, 1] \times$ $\mathbb{R}^d \rightarrow \mathbb{R}$ defined by
\begin{align}
\label{eq:u-definition}
    u(t, x)
    :=
    \mathbb{E}
    \big[
    \phi
    (X_1^{1-t,x})
    \big]
    =
    v(1-t,x),
    \qquad (t,x)\in[0,1]\times\mathbb{R}^d.
\end{align}
Using \eqref{eq:kol-v} together with the chain rule, we see that $u(t,x)$ is the unique classical solution of the following time-reversed Kolmogorov equation associated with SDE \eqref{eq:sde}:
\begin{align}
\label{eq:Kolmogorov-equation}
    \left\{\begin{array}{l}
    \partial_t u(t, x)
    =
    D u(t, x) f_{\beta}(x,1-t)
    +
    \frac{\beta}{2}
    \sum_{j=1}^d 
    D^2 u(t, x)
    (e_j, e_j), 
    \quad t\in (0,1], 
    \quad x \in \mathbb{R}^d, 
    \\
    u(0, x)=\phi(x)
    \quad  x \in \mathbb{R}^d. 
\end{array}\right.
\end{align}
This is one of the essential tools in the proof of the weak approximation result.
In addition, let us revisit the mean-square differentiability of random functions, quoted from \cite[Definition 1.1 in Chapter 2]{milstein2004stochastic} or \cite[Definition 2.2]{wang2024weak} as follows.

\begin{defn}
[Mean-square differentiable]
Let $\Psi:\Omega\times\mathbb{R}^d\to\mathbb{R}$ and
$\psi_i:\Omega\times\mathbb{R}^d\to\mathbb{R}$ be random functions satisfying
\begin{equation}
\label{eq:mean-square-differentiable}
\lim_{\tau_1\to 0}
\mathbb{E}
\left[
\left|
\tfrac{1}{\tau_1}
\big[
\Psi(x+\tau_1 e_i)-\Psi(x)
\big]
-\psi_i(x)
\right|^2
\right]
=0,
\qquad
\forall i\in [d],
\end{equation}
where $e_i$ is the unit vector in $\mathbb{R}^d$ with the $i$-th element being $1$.
Then $\Psi$ is called mean-square differentiable, with
$\psi=(\psi_1,\ldots,\psi_d) $being the derivative of $\Psi$ at $x$ in the mean-square sense.
Moreover, we denote
$
\mathcal{D}_{(i)}\Psi=\psi_i,
$ and 
$
\mathcal{D}\Psi(x)=\psi.
$
\end{defn}

We begin with the following estimates on the mean-square derivative $\mathcal{D}X^{s,\cdot}_t$.

\begin{lem}
\label{lem:X-mean-square-differential}
Let Assumptions \ref{ass:abs-continuous} and \ref{ass:higher-g-lipschitz} hold. Then, for any $x,v_1,v_2,v_3\in \mathbb{R}^d$ and $0 \leq s\leq t \leq 1$, 
\begin{align}
\label{eq:DX-estimation}
        \|\mathcal{D} X^{s,x}_t v_1\| 
        & \leq 
        e^{L_f (t-s)}
        \|v_1\|,
        \nonumber\\
        \|\mathcal{D}^2 X^{s,x}_t(v_1, v_2)\|_{L^2(\Omega,\mathbb{R}^d)}
        &\leq
        C e^{(3 L_f +1 / 2) 
        (t-s)}
        \|v_1\|
        \cdot\|v_2\|,
\end{align}
and 
\begin{align}
\label{eq:D^3X-estimation}
    \|\mathcal{D}^3 X^{s,x}_t(v_1, v_2,v_3)\|_{L^2(\Omega,\mathbb{R}^d)}
    \leq 
    C e^{(5 L_f +1) (t-s)} 
    \|v_1\|
    \cdot\|v_2\|
    \cdot\|v_3\|,
\end{align}
where $L_f > 0$ comes from \eqref{eq:f-lipschitz-condition}.
\end{lem}
The proof of the above lemma is postponed to Appendix \ref{appendix:D^3-x-proof}.
Notice that in general if $\phi \in C_b^k(\mathbb{R}^d)$ and the solution $X^{s,x}_t$ is $k$-times mean-square differentiable, then $u(t,\cdot)$ is $k$-times differentiable by using the chain rule (see \cite{Pang2024Linear,wang2024weak} for details). 
Moreover an explicit formula for the derivatives of $u(t,\cdot)$ can be given, in terms of the derivatives of $\phi$ and of $X^{s,x}_t$. 
Actually, we want to give an expression for the derivative of $u(t,\cdot)$ in terms of $\phi$ only, and not of its derivatives.
To this purpose, we recall the Bismut-Elworthy-Li formula (see \cite[Proposition 9.22]{da2014introduction} or \cite{elworthy1994formulae,ZXC2013,Zhang2016Stochastic} for details). 
%Such a formula, which applies when $X^x_t$ is mean-square differentiable and the diffusion term is non-degenerate, specific definitions are as follows.
\begin{lem}
[Bismut-Elworthy-Li formula]
\label{lem:Bismut-Elworthy-Li-formula}
Let $\{X_t^{s,x}\}_{0 \leq s\leq t \leq 1}$ be the solution of SDE \eqref{eq:sde}
starting from $x$ at time $s$. If $\{X_{t}^{s,x}\}_{0 \leq s\leq t \leq 1}$ is mean-square differentiable, then for any $\varphi \in C_b(\mathbb{R}^d), t \in (0,1]$ and $x, v \in \mathbb{R}^d$,
\begin{align}
\label{eq:bel}
    \big\langle 
    D(P_{s,t} \varphi(x)), v
    \big\rangle
    =
    \frac{1}{\sqrt{\beta}\,(t-s)} 
    \,
    \mathbb{E}
    \left[
    \int_{s}^t
    \big\langle
    D X_r^{s,x} v, 
    \mathrm{\, d} W_r
    \big\rangle 
    \varphi(X_t^{s,x})
    \right],
\end{align}
where $P_{s,t}(\cdot)$ defines a transition evolution operator of the corresponding ${SDE}$ \eqref{eq:sde} as 
$
P_{s,t} \varphi(x)
=
\mathbb{E}
\big[
\varphi(X_t^{s,x})
\big]
$.
Moreover, \eqref{eq:bel} remains true if $\varphi$ is a Borel function with polynomial growth.
\end{lem}

%%%%%%%%%%%%%%%%%%%%%%%%%%%%%%%%%%%%
\iffalse
\begin{lem}
[Bismut-Elworthy-Li formula]
\label{lem:Bismut-Elworthy-Li-formula-ori}
Let $\{X_t^x\}_{t \in [0,1]}$ be the solutions of the SDE \eqref{eq:sde} with the initial state $X_0^x=x$. If $\{X_t^x\}_{t \in (0,1]}$ is mean-square differentiable, then for any $\varphi \in C_b(\mathbb{R}^d), t \in [0,1]$ and $x, v \in \mathbb{R}^d$,
\begin{align}
\label{eq:bel}
    \big\langle 
    D(P_t \varphi(x)), v
    \big\rangle
    =
    \frac{1}{\sqrt{\beta}\,t} 
    \mathbb{E}
    \left[
    \int_0^t
    \big\langle
    D X_s^x v, 
    \mathrm{\, d} W_s
    \big\rangle 
    \varphi(X_t^x)
    \right],
\end{align}
where $P_t(\cdot)$ defines a transition semigroup of the corresponding ${SDE}$ \eqref{eq:sde} as 
$
P_t \varphi(x)
=
\mathbb{E}
\big[
\varphi(X_t^x)
\big]$.
Moreover, \eqref{eq:bel} remains true if $\phi$ is a Borel function with polynomial growth.
\end{lem}
%Recalling the fact that $P_t \phi(x)=u(t, x)$ by \eqref{eq-u-expression}, the Bismut-Elworthy-Li formula paves the way to presenting an expression for the derivatives of $u(t, \cdot)$ where it rely on $\phi$ only, and not of its derivatives.
%
\fi
%%%%%%%%%%%%%%%%%%%%%%%%%%%%%%%%%%%%

Combining with Lemma \ref{lem:X-mean-square-differential} and the above Bismut-Elworthy-Li formula, we are now in the position to obtain some regularization estimates with regard to the derivatives of $u(t, \cdot)=P_{1-t,t}(\cdot), \quad  t \in (0,1]$.

\begin{lem}
\label{lem:regularization-estimates-u}
Let Assumptions \ref{ass:abs-continuous} and \ref{ass:higher-g-lipschitz} hold. Then the function $u(t, \cdot) \in C_b^3(\mathbb{R}^d)$. 
Moreover, for any $ 
\phi \in C_b(\mathbb{R}^d)
$, there exists a constant $C$ independent of $d$ and $t$, such that, for any $x,v_1,v_2,v_3 \in \mathbb{R}^d$ and $t \in(0,1]$,
%\begin{align}
%    \big|
%    D(P_{1-t,1} \varphi(x))v
%    \big| 
%    & \leq 
%    \frac{C |\phi|_0}{\sqrt{t}}
%    \|v_1\|,\\
%    \big|
%    D^2(P_{1-t,1} \varphi(x)) (v_1,v_2)
%    \big|
%    & \leq 
%    \frac{C |\phi|_0 }{t}
%    \|v_1\|\cdot\|v_2\|,\\
%    \big|
%    D^3 (P_{1-t,1} \varphi(x))
%    (v_1,v_2,v_3)
%    \big|
%    & \leq 
%    \frac{C  |\phi|_0}
%    {t^{3/2}}
%    \|v_1\|\cdot\|v_2\|
%    \cdot\|v_3\|.  
%\end{align}
\begin{align}
\label{eq:D-u-est}
    \big|
     D u(t, x)v_1
    \big| 
    & \leq 
    \frac{C |\phi|_0}{\sqrt{t}}
    \|v_1\|,
\end{align}
\begin{align}
\label{eq:D2-u-est}
    \big|
    D^2 u(t, x) (v_1,v_2)
    \big|
    & \leq 
    \frac{C |\phi|_0 }{t}
    \|v_1\|\cdot\|v_2\|,
\end{align}
and
\begin{equation}
\label{eq:D-D2-D3-u-est}
  \big|
    D^3 u(t, x) 
    (v_1,v_2,v_3)
    \big|
    \leq 
    \frac{C  |\phi|_0}
    {t^{3/2}}
    \|v_1\|\cdot\|v_2\|
    \cdot\|v_3\|.  
\end{equation}
%\item[(2)]
%For every $ 
%\phi\in \mathrm{\rm Lip}(\mathbb{R}^d)
%$, there exists a constant $C$ independent of $d$, such that, for any $x,v_1,v_2,v_3 \in \mathbb{R}^d$,
%\begin{align}
%     \big|
%     D u(t, x) v_1
%     \big| 
%    & \leq 
%     C \sqrt{d}
%     \,
%     \|v_1\|,\nonumber\\
%     \big|
%     D^2 u(t, x) (v_1,v_2)
%     \big|
%     & \leq
%     C
%     \frac{\sqrt{d}}{\sqrt{t}}
%     \|v_1\|
%     \cdot\|v_2\|,,
%\end{align}    
%and
%\begin{equation}
%    \big|
%    D^3 u(t, x) (v_1,v_2,v_3)
%    \big| 
%    \leq 
%    C
%    \frac{\sqrt{d}}
%    {t}
%    \|v_1\|
%    \cdot\|v_2\|
%    \cdot\|v_3\|.
%\end{equation}
%\end{enumerate}
\end{lem}

The proof of Lemma \ref{lem:regularization-estimates-u} is given in Appendix \ref{appendix:proof-lemms-estimates-u}.

\subsection{Introduction of Malliavin calculus}
\label{subsec:intr-mal}
For the error analysis, we need to use Malliavin calculus.
We briefly introduce the necessary preliminaries in our setting.
More details can be found in \cite{nualart2006malliavin}.
Let $H=L^2\left([0, 1], \mathbb{R}^d\right)$ denote the space of all measurable functions
$\varphi: [0, 1] \to \mathbb{R}^d$ such that $\int_0^1 \|\varphi(t)\|^2 \dd t < \infty$. 
And $\mathbf{W}=(W^1_t,\cdots,W^d_t)_{t \in [0,1]}$ be a $d$-dimensional Wiener process on a probability space $(\Omega,\mathscr{F},\mathbb{P}).$ 
%where $\mathscr{F}$ is the $\sigma$-field generated by $\mathbf{W}$.
%
For any $h=(h_1,\cdots,h_d)\in H$, define the Wiener integral as 
\begin{equation}
    W(h):=
    \sum_{i=1}^d
    \int_0^1 h_i(t) 
    \dd W_t^i.
\end{equation}
To rigorously define the derivative operator $\mathscr{D}$, one may start with a class of special random variables
\begin{equation}
    F \in \mathcal{S}
    :=
    \left\{f\left(W\left(h_1\right), \ldots, W\left(h_n\right)\right):
    f \in C_p^{\infty}(\mathbb{R}^n), 
    n \geq 1,
    h_i \in H, 
    1 \leq i \leq n
    \right\},
\end{equation}
where $C_p^{\infty}\left(\mathbb{R}^n\right)$ is the set of all infinitely differentiable functions $f: \mathbb{R}^n \rightarrow \mathbb{R}$ such that all the partial derivatives have at most polynomial growth. 
We are in a position to give the definition of the derivative operator $\mathscr{D}$.
\begin{defn}
[Derivative operator $\mathscr{D}$]
\label{def:Derivative-operator} 
    Consider a variable $F=f\left(W\left(h_1\right), \ldots, W\left(h_n\right)\right) \in \mathcal{S}$ for some $n \geq 1$. For any $t \in [0,1]$, the Malliavin derivative at $t$ is defined as
    \begin{equation}
        \mathscr{D}_t F=\sum_{j=1}^n \frac{\partial f}{\partial x_j}\left(W\left(h_1\right), \ldots, W\left(h_n\right)\right) h_j(t)
    \end{equation}
    Moreover, $\mathscr{D}F:=\left(\mathscr{D}_tF\right)_{t\in[0,1]}$ is the process in $L^2(\Omega ,H)$.
\end{defn}
For any $p \geq  1$, we denote by $\mathbb{D}^{1,p}$ as the closure of the class of smooth random variables in $\mathcal{S}$ with respect to the norm
\begin{equation}
    \|F\|_{\mathbb{D}^{1,p}}
    :=
    \left(\mathbb{E}
    \left[
    |F|^p\right]
    +
    \mathbb{E}
    \left[
    \int_0^1
    \left\|
    \mathscr{D}_s F\right\|^p \mathrm{\, d} s
    \right]
    \right)^{1 / p}.
\end{equation}
Using the duality, one can define the following divergence operator $\delta$.
\begin{defn}
\label{def:divergence-operator}
    The domain of the divergence operator $\operatorname{Dom} \delta$ in $L^2(\Omega)$ is the set of processes $u \in L^2(\Omega , H)$ such that there exists $\delta(u) \in L^2(\Omega)$ satisfying the duality relationship
    \begin{equation}
        \mathbb{E}
        \big[
        \langle \mathscr{D} F, u\rangle_H
        \big]
        =
        \mathbb{E}\big[F \delta(u)\big],
    \end{equation}
    for any $F \in \mathbb{D}^{1,2}$.
    Here, $\langle f, g\rangle_H=\int_0^1 f \cdot g d t$ is the inner product in $L^2\left([0, 1] ,\mathbb{R}^d\right)$.
\end{defn}
Definition \ref{def:divergence-operator} is in fact an integration by parts formula. 
Given $\mathbf{h} \in H,$ we can define the Malliavin derivative along the direction $\mathbf{h}$, denoted by $\mathscr{D}^{\mathbf{h}}F,$ as the following 
\begin{equation}
    \mathscr{D}^{\mathbf{h}}F
    =
    \langle
    \mathscr{D}F, \mathbf{h}
    \rangle_H
    =
    \int_0^1
    \langle
    \mathbf{h}_s,\mathscr{D}_s F
    \rangle
    \mathrm{\, d} s.
\end{equation}
The following are some basic properties in Malliavin calculus which we will use frequently
later.
\begin{prop}
    1. (Chain rule for Malliavin derivative) Let 
    $\varphi: \mathbb{R} \rightarrow \mathbb{R}$ 
    be a continuous differentiable function such that $\left|\varphi^{\prime}(x)\right| \leq C\left(1+|x|^\alpha\right)$ for some $\alpha \geq 0$. 
    Let $F \in \mathbb{D}^{1, p}$ for some $p \geq \alpha+1$. Then, $\varphi(F)$ belongs to $\mathbb{D}^{1, q}$, where $q=p /(\alpha+1)$, and
    \begin{equation}
    \label{eq:Chain-rule-for-Malliavin-derivative}
       \mathscr{D}(\varphi(F))
       =
       \varphi^{\prime}(F)
       \mathscr{D} F .
    \end{equation}

    2. (Product rule for divergence) 
    Let $G \in \mathbb{D}^{1,2}$ be a real-valued random variable such that $G u \in L^2\left(\Omega, L^2[0, 1]\right)$, then
    \begin{equation}
    \label{eq:Product-rule-for-divergence}
        \delta(G u)=G \delta(u)-\langle 
        \mathscr{D} G, u
        \rangle_H,
    \end{equation}
    where we suppose that the right-hand side is integrable.
\end{prop}
To facilitate the subsequent error analysis, we define $N_1 \in \mathbb{N}$ such that 
$0 < N_1h < \frac{1}{2} \leq (N_1+1)h<1.$
For the grid points $t_n, n=N_1+1,\cdots,N-1$, one is able to find a grid point $t_m$ such that 
\begin{equation}
\label{eq:distance-t_n-t_m}
    0 
    \leq t_m 
    < t_n < t_{n+1} 
    \leq 1, 
    \quad 
    \text{with} 
    \quad
    \frac{1}{2} \leq t_{n} -t_m \leq 1.
\end{equation}
We introduce the following process 
\begin{equation}
\label{eq:defination-direction-h}
    h_{\ell}(t)
    :=
    \tfrac{1}{\sqrt{\beta}}
    \sum_{i=m}^n \mathbf{1}_{[t_i, t_{i+1})}(t)
    \left(
    \frac{e_{\ell}}
    {t_{n}-t_m}
    -
    \frac{t_i-t_m}
    {t_{n}-t_m} 
    D \tilde{f}_{\beta}
    \big(
    \tilde{Y}_{t_i},t_i
    \big) 
    e_{\ell}\right), 
    \quad t \in[t_m, t_{n+1}),
\end{equation}
where 
$\{e_{\ell}\}_{\ell \in [d]}$ denotes the orthonormal basis of $\mathbb{R}^d$. 
We directly observe that the process \eqref{eq:defination-direction-h} is adapted with respect to $\mathcal{F}_t, t \in[t_m, t_{n+1})$, and $h_{\ell}(\cdot) \in L^2(\Omega, H)$ due to Proposition \ref{prop:tilde-Df-bound}.
Thanks to \cite [Proposition 1.3.11]{nualart2006malliavin}, we deduce that $h_{\ell} \in \operatorname{Dom} \delta$ and
\begin{equation}
    \delta(h_{\ell})
    =
    \int_{t_m}^{t_{n+1}}
    \big\langle 
    h_{\ell}, \mathrm{\, d} W_s
    \big\rangle .
\end{equation}

Moreover, we are in the position to show the Malliavin regularity estimates of the process \eqref{eq:Euler-Maruyama-contimuous-version} with respect to the direction $h_{\ell}$.
\begin{lem}
\label{lem:Malliavin-regularity-estimates-numerical}
   Let Assumptions 
   \ref{ass:abs-continuous} and \ref{ass:higher-g-lipschitz} hold. 
   Given the direction $h_{\ell}$ as defined by \eqref{eq:defination-direction-h}, for any $p \in[1, \infty)$, there exists a positive constant $C$ depending on $p$ such that,
   \begin{equation}
   \label{eq:delta-h-estimation}
       \big\|
       \delta
       (h_{\ell})
       \big\|_{L^{2 p}(\Omega, \mathbb{R})} 
       \leq C,
   \end{equation}
   and, for $t_m \leq t_i \leq s \leq t_{i+1} \leq t_{n+1}$,
   \begin{equation}
   \label{eq:numerical-Malliavin-estimation}
   \mathscr{D}^{h_{\ell}} 
   \tilde{Y}_s
   =
   \frac{s-t_m}
   {t_{n}-t_m} 
   e_{\ell}.
   \end{equation}
   Moreover, there exists a positive constant $C$, depending on $L_f$, such that
   \begin{equation}
   \label{eq:f-Malliavin-estimation}
    \begin{aligned}
    \big\|
   \mathscr{D}^{h_{\ell}} \tilde{f}_{\beta}
   (\tilde{Y}_{t_n},{t_n})
    \big\|
    & \leq 
    C .
    \end{aligned}
   \end{equation}
\end{lem}
\begin{proof}
\textbf{Proof.}
The proofs of the first two formulas in Lemma \ref{lem:Malliavin-regularity-estimates-numerical} can be found in \cite[Lemma 6.4]{PWW25}.
In the following, we only need to prove \eqref{eq:f-Malliavin-estimation}.
Using Proposition \ref{prop:tilde-Df-bound} and the chain rule for Malliavin derivative \eqref{eq:Chain-rule-for-Malliavin-derivative}, we get
\begin{align*}
   \big\|
   \mathscr{D}^{h_{\ell}} \tilde{f}_{\beta}
   (\tilde{Y}_{t_n},{t_n})
   \big\|
   & =
   \big\|
    D \tilde{f}_{\beta}
   (\tilde{Y}_{t_n},{t_n})
    \mathscr{D}^{h_{\ell}}
    \tilde{Y}_{t_n}
    \big\|
    \leq
    L_f
    \big\|
    \mathscr{D}^{h_{\ell}}
    \tilde{Y}_{t_n}
    \big\|
    =
    \frac{L_f(t_n-t_m)}
    {t_{n}-t_m} 
    \|e_{\ell}\|
    \leq 
    C.
\end{align*}
\end{proof}

\section{Proof of Theorem \ref{thm:exact-f} and Theorem \ref{thm-main-rerult-w0-w1}}
\label{section:proof-theorem}
\noindent
In the subsequent analysis, for simplicity, we denote, for any $n \in [N-1]_0$, 
\begin{align*}
        & 
        p(x,t_n)
        :=
        \beta\,
        \mathbb{E}_{\xi} 
        \big[
        \nabla 
        g_\beta(x+\sqrt{(1-t_n)\beta} \,\xi)
        \big], 
        \quad 
        p_{M}(x,t_n,\xi)
        :=
        \tfrac{\beta}
        {M}
        \sum_{j=1}^M 
        \nabla g_\beta(x+\sqrt{(1-t_n)\beta} \,\xi_{t_n,j}), 
        \nonumber\\
        & 
        q(x,t_n)
        := 
        \mathbb{E}_{\xi}
        \big[
        g_\beta(x+\sqrt{(1-t_n)\beta} \,\xi)
        \big], 
        \quad \quad
        q_{M}(x,t_n,\xi)
        :=
        \tfrac{1}{M}
        \sum_{j=1}^{M} 
        g_\beta(x+\sqrt{(1-t_n)\beta} \,\xi_{t_n,j}).
\end{align*}

Before proving Theorem \ref{thm:exact-f} and Theorem \ref{thm-main-rerult-w0-w1}, we recall the following lemma (see \cite[Lemma A.6]{MRjiao} for a similar proof).
\begin{lem}
\label{lem:jiao-M-error}
    Let Assumptions \ref{ass:abs-continuous} and \ref{ass:g-nabla-g-lip} hold. Then for any $x\in \mathbb{R}^d$ and $n \in [N-1]_0$,
    \begin{align*}
        \mathbb{E}_{\xi}
        \Big[
        \big\|
        p(x,t_n)-p_M(x,t_n,\xi)
        \big\|^2
        \Big]
        \leq 
        \frac{L' d}
        {M},
        \quad 
        \mathbb{E}_{\xi}
        \Big[
        \big\|q(x,t_n)-q_M(x,t_n,\xi)
        \big\|^2
        \Big]
        \leq 
        \frac{L' d}
        {M},
    \end{align*}
    with  
    $L'
    :=
    2 L_g^2\, \beta^2$, 
    where $L_g$ comes from Assumption \ref{ass:g-nabla-g-lip},
    and $M$ is the number of samples used in the Monte Carlo estimator \eqref{eq:drift-est}.
\end{lem}

Now we first prove Theorem \ref{thm-main-rerult-w0-w1}.

\begin{proof}
\textbf{Proof of Theorem \ref{thm-main-rerult-w0-w1}}.
%%%%%%%%%%%%%%%%%%%
\iffalse
Thanks to Kantorovich's dual formula, the total variation distance and $\mathcal{W}_1$ distance between the law of the SFS and the target distribution $\mu$ can be rewritten as following,
\begin{align}
         \big\|\mathcal{L}aw
         \big(
         \tilde{Y}_1^{0}
         \big), \mu
         \big\|_{TV}
         & = 
         \sup_{\phi \in c_b(\mathbb{R}^d),[\phi]_0 \leq 1}
         \Big|
         \mathbb{E}
         \big[
         \phi\big(
         \tilde{Y}_1^0
         \big)
         \big]
         -
         \mathbb{E}
         \big[
         \phi\big(X_1^0\big)
         \big]
         \Big|,\nonumber\\
         \mathcal{W}_1
         \big(
         \mathcal{L}aw
         (\tilde{Y}_1^{0}), \mu
         \big)
         & =
         \sup_{\phi \in C(\mathbb{R}^d),[\phi]_{Lip} \leq 1}
         \Big|
         \mathbb{E}
         \big[
         \phi\big(
         \tilde{Y}_1^0
         \big)
         \big]
         -
         \mathbb{E}
         \big[
         \phi(X_1^0)
         \big]
         \Big|.
    \end{align}
%It is crucial to highlight that when $p = 0$, i.e., the total variation distance, the $\phi$ reads on bounded Borel-measurable functions, and when $p = 1$, the $\phi$ reads on Lipschitz continuous functions.
\fi
%%%%%%%%%%%%%%%%%%%%%%%%%%
By \eqref{eq:definition-tv-rv}, the total variation distance between the law of the SFS algorithm and the target distribution induced by \eqref{eq:target-distribution} boils down to the weak error analysis of the SFS scheme \eqref{eq:e-m-mc} for the SDEs \eqref{eq:sde}:
\begin{equation}
\label{eq:tv-weak-error}
\big\|
    \mathcal{L}aw
    \big(
    \tilde{Y}_1
    \big)
    - 
    \mu
    \big\|_{\mathrm{TV}}
    :=
    \sup_{\phi \in c_b(\mathbb{R}^d),|\phi|_0 \leq 1}
         \Big|
         \mathbb{E}
         \big[
         \phi\big(
         \tilde{Y}_1
         \big)
         \big]
         -
         \mathbb{E}
         \big[
         \phi\big(X_1
         \big)
         \big]
         \Big|.
\end{equation}
By the telescoping argument, \eqref{eq:u-definition} and \eqref{eq:Euler-Maruyama-contimuous-version}, we decompose the weak error between them as follows
    \begin{align}
    \label{eq:weak-error-decomposition}
        &
        \Big|
        \mathbb{E}
        \big[
        \phi
        \big(
        \tilde{Y}_1
        \big)
        \big]
        -
        \mathbb{E}
        \big[
        \phi
        \big(X_1\big)
        \big]
        \Big|\notag \\
        & = 
        \Big|
        u\big(0, \tilde{Y}_1
        \big)
        -
        u\big(1, 0\big)
        \Big|\notag\\
        & =
        \Bigg|
        \sum_{n=0}^{N-1} 
        \mathbb{\,E}
        \Big[
        u\big(
        1-t_{n+1},
        \tilde{Y}_{t_{n+1}}
        \big)
        \Big]
            -
            \mathbb{E}
            \Big[
            u\big(1-t_n, \tilde{Y}_{t_n}\big)
            \Big]
        \Bigg|\notag \\
        & =
        \Bigg|
        \sum_{n=0}^{N-1} 
        \mathbb{\,E}
        \bigg[
          \int_{t_n}^{t_{n+1}} 
           D u\big(1-s,\tilde{Y}_s
           \big)
           \Big(
           \tilde{f}
           _{\beta}
           \big(
           \tilde{Y}_{t_n}, t_n
           \big)
           -
           f_{\beta}
           \big(
           \tilde{Y}_s,s
           \big)
           \Big) 
           \mathrm{\, d} s
           \bigg]
         \Bigg|\notag \\
%        & =
%        \Bigg|
%          \sum_{n=0}^{N-1} \mathbb{\,E}
%          \bigg[
%          \int_{t_n}^{t_{n+1}} 
%            D u\big
%            (1-s, \tilde{Y}_s
%            \big)
%            \Big(
%            \tilde{f}
%            \big(
%            \tilde{Y}_{t_n},
%             t_n
%            \big)
%             -f_{\beta}\big(\tilde{Y}_{t_n},t_n\big)
%             +
%             f
%             \big(
%             \tilde{Y}_{t_n},
%             t_n
%             \big)
%             -
%             f_{\beta}\big(
%             \tilde{Y}_s,s
%             \big)
%             \Big) 
%             \mathrm{\, d} s
%             \bigg]\Bigg|\notag \\
        & \leq 
        \underbrace{\Bigg|
        \sum_{n=0}^{N-1} \mathbb{E}
        \bigg[
        \int_{t_n}^{t_{n+1}} 
        D u\big(1-s, \tilde{Y}_s\big)
        \Big(
        \tilde{f}_{\beta}
        \big(
        \tilde{Y}_{t_n},
        t_n
        \big)
        -f_{\beta}\big(\tilde{Y}_{t_n},t_n\big)
        \Big) 
        \mathrm{\, d} s
        \bigg]\Bigg|}_{=:J_1}
        \notag \\
        & \quad +
        \underbrace{\Bigg| 
        \sum_{n=0}^{N-1} \mathbb{E}
        \bigg[
        \int_{t_n}^{t_{n+1}} 
        D u\big(1-s, \tilde{Y}_s\big)
        \Big(
        f_{\beta}\big(
        \tilde{Y}_{t_n},t_n
        \big)
        -
        f_{\beta}\big(
        \tilde{Y}_s,s
        \big)
        \Big) 
        \mathrm{\, d} s
        \bigg]
        \Bigg|}_{=:J_2},
    \end{align}
the third equation is due to the It\^o formula and the corresponding Kolmogorov equation \eqref{eq:Kolmogorov-equation}, as follows
\begin{align}
        &
        \mathbb{E}
        \Big[ 
        u
        \big(
        1-t_{n+1},
        \tilde{Y}_{t_{n+1}}
        \big)
        -
        u
        \big(
        1-t_n, \tilde{Y}_{t_n}
        \big)
        \Big]\nonumber\\
        & = 
        \mathbb{E}
        \bigg[
        \int^{t_{n+1}}_{t_n}
        \Big(
        -\partial_{s} 
        u
        \big(
        1-s,\tilde{Y}_s
        \big)
        + 
        D u
        \big(
        1-s, \tilde{Y}_s
        \big)
        \tilde{f}_{\beta}
        \big(
        \tilde{Y}_{t_n},t_n
        \big)
        +
        \tfrac{\beta}{2}
        \sum_{j=1}^d 
        D^2 u
        \big(
        1-s, \tilde{Y}_s
        \big)
        \big(e_j,e_j\big) 
        \Big)
        \mathrm{\, d} s
        \bigg]\nonumber \\
        & = 
        \mathbb{E}
        \bigg[
        \int^{t_{n+1}}_{t_n} 
        D u
        \big(1-s, \tilde{Y}_s
        \big)
        \Big(
        \tilde{f}_{\beta}
        \big(
        \tilde{Y}_{t_n},t_n
        \big)
        -
        f_{\beta}\big(
        \tilde{Y}_s,s
        \big)
        \Big) 
        \mathrm{\, d} s
        \bigg].
\end{align}
Regarding $J_1$, one begins with the following decomposition:
\begin{align}
    J_1
    & \leq  
    \Bigg|
    \sum_{n=0}^{N-1} 
    \mathbb{E}
    \bigg[
    \int_{t_n}^{t_{n+1}}
    \Big( 
    D u
    \big(1-s, \tilde{Y}_s
    \big)
    -
    D u
    \big(
    1-s, \tilde{Y}_{t_n}
    \big)
    \Big)
    \Big(
    \tilde{f}_{\beta}
    \big(
    \tilde{Y}_{t_n},t_n
    \big)
    -
    f_{\beta}\big(
    \tilde{Y}_{t_n},t_n
    \big)
    \Big)         
    \mathrm{\, d} s
    \bigg]
    \Bigg| 
    \nonumber \\
    & \quad + 
    \Bigg|
    \sum_{n=0}^{N-1} 
    \mathbb{E}
    \bigg[
    \int_{t_n}^{t_{n+1}} 
    D u
    \big(1-s, \tilde{Y}_{t_n}
    \big)
    \Big(
    \tilde{f}_{\beta}
    \big(
    \tilde{Y}_{t_n},t_n
    \big)
    -
    f_{\beta}\big(
    \tilde{Y}_{t_n},t_n
    \big)
    \Big) 
    \mathrm{\, d} s
    \bigg]
    \Bigg| 
    \nonumber\\
    & \leq 
    \underbrace{
    \Bigg|
    \sum_{n=0}^{N-1} 
    \mathbb{E}
    \bigg[
    \int_{t_n}^{t_{n+1}}
    \Big( 
    D u
    \big(
    1-s, \tilde{Y}_s
    \big)
    -
    D u
    \big(
    1-s, \tilde{Y}_{t_n}
    \big)
    \Big)
    \tilde{f}_{\beta}
    \big(
    \tilde{Y}_{t_n},t_n
    \big) 
    \mathrm{\, d} s
    \bigg]
    \Bigg| }_{=:J_{1,1}}
    \notag \\
    & \quad +
        \underbrace{
        \Bigg|
        \sum_{n=0}^{N-1} 
        \mathbb{E}
        \bigg[
    \int_{t_n}^{t_{n+1}}
        \Big( 
        D u
        \big(
        1-s, \tilde{Y}_s
        \big)
        -
        D u
        \big(
        1-s, \tilde{Y}_{t_n}
        \big)
        \Big)
        f_{\beta}
        \big(
    \tilde{Y}_{t_n},t_n
    \big) 
    \mathrm{\, d} s
    \bigg]
    \Bigg|}_{=:J_{1,2}}
    \notag \\
    & \quad + 
    \underbrace{
    \Bigg|
    \sum_{n=0}^{N-1}
    \mathbb{E}
    \left[
    \int_{t_n}^{t_{n+1}} 
    D u
    \big(1-s, \tilde{Y}_{t_n}
    \big)
    \left(
    \tilde{f}_{\beta}
    \big(
    \tilde{Y}_{t_n},t_n
    \big)
    -
    f_{\beta}\big(
    \tilde{Y}_{t_n},t_n
    \big)
    \right) 
    \mathrm{\, d} s
    \right]
    \Bigg|}_{=:J_{1,3}}.
\end{align}
%To carry out the error analysis, we first decompose the relevant terms and establish several preliminary estimates.
%
We firstly present an estimate of $J_{1,3}$.
For notational simplicity, we write
$p=p(x,t_n),\,q=q(x,t_n),\,
p_M=p_M(x,t_n,\xi)$ and $q_M=q_M(x,t_n,\xi)$, for fixed $x\in\mathbb{R}^d$ and $n\in[N-1]_0$.
Applying the multivariate Taylor expansion to the 
$
S: \mathbb{R}^d
\times
\mathbb{R} \setminus \{ 0 \}
\rightarrow \mathbb{R}^d$, 
we note that
\begin{align}
\label{eq:f-taylor-exp}
    S(p_M,q_M)
    =
    \frac{p_M}{q_M}
    & =
    \frac{p}{q}
     +
    \frac{1}{q}
    \big(
    p_M-p
    \big)
    -
    \frac{p}{q^2}
    \big(
    q_M-q
    \big)
    +
    \mathcal{R}_S
    (x,t_n,\xi),
\end{align}
where, $r\in[0,1]$,
\begin{align}
   \mathcal{R}_S(x,t_n,\xi)
   & :=
   2
   \int_0^1
   (1-r)
   \left[
   -
   \frac{(p_M-p)(q_M-q)}
   {
   \big(
   q+r(q_M-q)
   \big)^2}
   +
   \frac{(p+r(p_M-p))(q_M-q)^2}
   {\big(
   q+r(q_M-q)
   \big)^3}
   \right]
   \mathrm{\, d} r.
%   & =
%   2
%   \int_0^1
%   (1-r)
%   \bigg[
%   -
%   \frac{(p_M-p)(q_M-q)}
%   {\big(
%   (1-r)q+rq_M
%   \big)^2}
%   +
%   \frac{(1-r)p(q_M-q)^2}
%   {\big((1-r)q+rq_M
%  \big)^3}
%    +
%   \frac{rp_M(q_M-q)^2}
%   {\big((1-r)q+rq_M
%   \big)^3}
%   \bigg]
%   \mathrm{\, d} r. 
\end{align}
Moreover, in light of 
%the inequality $(a+b)^2 \geq a^2+b^2 $  for any $a,b \geq 0 $, 
Assumption \ref{ass:g-nabla-g-lip}, Lemma \ref{lem:jiao-M-error} and the H\"older inequality, we can get
\begin{align}
    \left\|
    \frac{(p_M-p)(q_M-q)}
    {\big(
    (1-r)q+rq_M
    \big)^2}
    \right\|_{L^1(\Omega,
    \mathbb{R}^d)}
    & \leq 
    \frac{1}
    {\rho^2}
    \big\|
    p_M-p
    \big\|_{L^2(\Omega,
    \mathbb{R}^d)}
    \big\|
    q_M-q
    \big\|_{L^2(\Omega,
    \mathbb{R})}
    \nonumber\\
    & \leq
    \frac{L' d}
    {\rho^2 M}.
\end{align}
%where the last inequality follows from  $(1-r)^2+r^2 \geq \frac{1}{2}$, for any $r\in [0,1]$.
%
Due to Assumption \ref{ass:g-nabla-g-lip}, we have, for any $x\in \mathbb{R}^d, \, n\in [N]_0$,
\[
q(x,t_n)
+
r
\big(
q_M(x,t_n,\xi)-q(x,t_n)
\big)
=
(1-r)q(x,t_n)+r q_M(x,t_n,\xi)
\geq 
\rho,
\qquad r\in[0,1].
\]
Using the fact that $g_{\beta}, \nabla g_{\beta}$ are Lipschitz continuous and the definition of the operator norm \eqref{SFS-W2-eq:operator-norm-def}, one can get, for any $x\in \mathbb{R}^d, \, n\in [N]_0$,
$$
\|p(x,t_n)\|\leq \beta L_g, 
\quad 
\|p_M(x,t_n,\xi)\|\leq \beta L_g ,
$$
whose proof can be found in \cite[Appendix A]{Wang2026Multimodal}. Therefore,
\begin{align}
\left\|
\frac{(p+r(p_M-p))(q_M-q)^2}
{\big(q+r(q_M-q)\big)^3}
\right\|_{L^1(\Omega,\mathbb R^d)}
\nonumber
& \leq
\frac{1}{\rho^3}
\,
\mathbb E_\xi
\Big[
\big\|(1-r)p+rp_M\big\|
\|q_M-q\|^2
\Big]
\nonumber\\
& \leq
\frac{\beta L_g}{\rho^3}
\mathbb E_\xi
\Big[
\big\|q_M-q\big\|^2
\Big]
\nonumber\\
& \leq
\frac{\beta L_g L'd}
{\rho^3 M}.
\end{align}
Noting that 
$$
\mathbb{E}_\xi
\big[p_M-p\big]
=
0,
\quad
\text{and}
\quad
\mathbb{E}_\xi
\big[q_M-q\big]
=
0,
$$
taking a conditional expectation, putting \eqref{eq:f-taylor-exp} into the estimate of $J_{1,3}$ and using Lemma \ref{lem:regularization-estimates-u}, we can directly get
\begin{align}
\label{eq:J_13-est}
    J_{1,3}
    & = 
    \Bigg|
    \sum_{n=0}^{N-1} \mathbb{E}
    \bigg[
    \int_{t_n}^{t_{n+1}} 
    D u
    \big(
    1-s, \tilde{Y}_{t_n}
    \big)
    \Big(
    \tilde{f}_{\beta}
    \big(
    \tilde{Y}_{t_n},t_n
    \big)
    -
    f_{\beta}
    \big(
    \tilde{Y}_{t_n},t_n
    \big)
    \Big) 
    \mathrm{\, d} s
    \bigg]
    \Bigg|
    \nonumber\\
    & =
    \Bigg|
    \sum_{n=0}^{N-1}
    \mathbb{E}
    \bigg[
    \mathbb{E}
    \bigg[
    \int_{t_n}^{t_{n+1}} 
    D u
    \big(
    1-s, \tilde{Y}_{t_n}
    \big)
    \Big(
    \tilde{f}_{\beta}
    \big(
    \tilde{Y}_{t_n},t_n
    \big)
    -
    f_{\beta}
    \big(
    \tilde{Y}_{t_n},t_n
    \big)
    \Big) 
    \mathrm{\, d} s
    \bigg]
    \bigg|
    \sigma(
    \tilde{Y}_{t_n})
    \bigg]
    \Bigg|
    \nonumber\\
    & =
    \Bigg|
    \sum_{n=0}^{N-1} 
    \mathbb{E}
    \bigg[
    \mathbb{E}_\xi
    \bigg[
    \int_{t_n}^{t_{n+1}} 
    D u
    (1-s, x)
    \big(
    \tilde{f}_{\beta}(x,t_n)
    -
    f_{\beta}
    (x,t_n)
    \big)
    \mathrm{\, d} s 
    \bigg]
    _{x=\tilde{Y}_{t_n}}  
    \bigg]
    \Bigg|
    \nonumber\\
    & =
    \Bigg|
    \sum_{n=0}^{N-1} 
    \mathbb{E}
    \bigg[
    \mathbb{E}_\xi
    \bigg[
    \int_{t_n}^{t_{n+1}} 
    D u(1-s, x)
    \mathcal{R}_S
    (x,t_n,\xi)
    \mathrm{\, d} s 
    \bigg]
    _{x=\tilde{Y}_{t_n}}  
    \bigg]
    \Bigg|
    \nonumber\\
    & \leq 
    C
    \sum_{n=0}^{N-1} 
    \int_{t_n}^{t_{n+1}} 
    \frac{1}{\sqrt{1-s}}
    \big\|
    \mathcal{R}_S(x,t_n,\xi)
    \big\|_{L^1(\Omega_\xi,
    \mathbb{R}^d)}
    \dd s
    \nonumber\\
    & \leq 
    \frac{C d}{M},
\end{align}
where 
\[
 \sum_{n=0}^{N-1} 
    \int_{t_n}^{t_{n+1}} 
    \frac{1}{\sqrt{1-s}}
    \,
    \dd s
    =
    \int_0^1
    \frac{1}{\sqrt{t}}
    \,
    \dd t
    = 
    2.
\]
%the second equation is due to the nature of conditional expectation, and the third equation is due to $\tilde{Y}_{t_n}$ and $\xi$ in the term of $\frac{p_m(x)}{q_{\hat{m}}(x)}$ are independent combined with Lemma \ref{lem-independent-expection}.
%
For the estimate of $J_{1,1}$, to overcome the possible singularities in Lemma \ref{lem:regularization-estimates-u}, we split the time interval $[0,1]$ into $[0,1-h^2]$ and $[1-h^2,1]$, as follows 
\begin{align}
     \label{eq:J_11-TV-experssion}
         J_{1,1}
         \leq &
         \Bigg|
         \sum_{n=0}^{N-2} 
         \mathbb{E}
         \bigg[
         \int_{t_n}^{t_{n+1}}
         \Big( 
         D u
         \big(
         1-s, \tilde{Y}_s
         \big)
         -
         D u
         \big(
         1-s, \tilde{Y}_{t_n}
         \big)
         \Big)
         \tilde{f}_{\beta}
         \big(
         \tilde{Y}_{t_n},t_n
         \big) 
         \mathrm{\, d}s
         \bigg]
         \nonumber\\
         & 
         \underbrace{
         \quad \quad +
         \mathbb{E}
         \bigg[
         \int_{1-h}^{1-h^2}
         \Big( 
         D u
         \big(
         1-s, \tilde{Y}_s
         \big)
         -
         D u
         \big(
         1-s, \tilde{Y}_{t_{N-1}}
         \big)
         \Big)
         \tilde{f}_{\beta}
         \big(
         \tilde{Y}_{t_{N-1}},t_{N-1}
         \big) 
         \mathrm{\, d} s
         \bigg]
         \Bigg|}_{=:J_{1,1,1}}
         \nonumber\\
         & 
         \underbrace{+
         \Bigg|
         \mathbb{E}
         \bigg[
         \int_{1-h^2}^{1}
         \Big( 
         D u
         \big(
         1-s, \tilde{Y}_s
         \big)
         -
         D u
         \big(
         1-s, \tilde{Y}_{t_{N-1}}
         \big)
         \Big)
         \tilde{f}_{\beta}
         \big(
         \tilde{Y}_{t_{N-1}},t_{N-1}
         \big) 
         \mathrm{\, d} s
         \bigg]
         \Bigg|}_{=:J_{1,1,2}}.
\end{align}
In order to treat $J_{1,1,1}$, applying the generalized It\^o formula \cite[Theorem 4.2.1]{oksendal2013stochastic} to $Du(1-s,\cdot)$, one note that
\iffalse
\begin{equation}
    \begin{aligned}
        D u\left(1-s, \tilde{Y}_s\right)
        & =
        D u\left(1-s, \tilde{Y}_{t_n}\right)
        +
        \int_{t_n}^s
        D^2 u\left(1-s,\tilde{Y}_r\right)
        \tilde{f}
        \left(
        \tilde{Y}_{t_n},t_n
        \right) 
        \mathrm{\, d} r\\
        & \quad +
        \int_{t_n}^s
        D^2 u\left(1-s,\tilde{Y}_r\right)
        \mathrm{\, d} W_r
        +
        \frac{1}{2}
        \sum_{j=1}^d
        \int_{t_n}^s
        D^3 u\left(1-s,\tilde{Y}_r\right)
        \left(
        e_j,e_j
        \right)
        \mathrm{\, d} r,
    \end{aligned}
\end{equation}
thus the term $J_{1,1,1}$ can be rewritten as 
\fi
    \begin{align}
        J_{1,1,1}
         & \leq
         \Bigg|
         \sum_{n=0}^{N-2} \mathbb{E}
         \bigg[
         \int_{t_n}^{t_{n+1}}
         \int_{t_n}^s
         D^2 u
         \big(1-s,\tilde{Y}_r
         \big)
         \Big(
         \tilde{f}_{\beta}
         \big(
         \tilde{Y}_{t_n},t_n
         \big), 
         \tilde{f}_{\beta}
         \big(
         \tilde{Y}_{t_n},t_n
         \big)
         \Big)
         \mathrm{\, d} r
         \mathrm{\, d} s 
         \bigg]
         \notag\\
         & \quad
         \underbrace{
         \quad +
         \mathbb{E}
         \bigg[
         \int_{1-h}^{1-h^2}
         \int_{1-h}^s 
         D^2 u
         \big(
         1-s,\tilde{Y}_r
         \big)
         \Big(
         \tilde{f}_{\beta}
         \big(
         \tilde{Y}_{t_{N-1}},t_{N-1}
         \big), 
         \tilde{f}_{\beta}
         \big(
         \tilde{Y}_{t_{N-1}},t_{N-1}
         \big)
         \Big) 
         \mathrm{\, d} r
         \mathrm{\, d} s 
         \bigg]
         \Bigg|}
        _{=:J_{1,1,1}^{(1)}}
        \notag\\
        & \quad +
        \frac{\beta}{2}
        \Bigg|
        \sum_{n=0}^{N-2}
        \sum_{j=1}^{d}
         \mathbb{E}
         \bigg[
         \int_{t_n}^{t_{n+1}}
         \int_{t_n}^s
         D^3 u
         \big(
         1-s,\tilde{Y}_r
         \big)
         \Big(e_j,e_j,
         \tilde{f}_{\beta}
         \big(
         \tilde{Y}_{t_n},t_n
         \big)
         \Big)
         \mathrm{\, d} r
         \mathrm{\, d} s 
         \bigg]\notag\\
         & \quad
         \underbrace{ 
         \quad +
         \sum_{j=1}^{d}
         \mathbb{E}
         \bigg[
         \int_{1-h}^{1-h^2}
         \int_{1-h}^s
         D^3 u
         \big(
         1-s,\tilde{Y}_r
         \big)
         \Big(e_j,e_j,
         \tilde{f}_{\beta}
         \big(
         \tilde{Y}_{t_{N-1}},t_{N-1}
         \big)
         \Big)
         \mathrm{\, d} r
         \mathrm{\, d} s 
         \bigg]
         \Bigg|}
         _{=:J_{1,1,1}^{(2)}}
         \nonumber\\
         & \quad +
         \sqrt{\beta}
         \,
         \Bigg|
         \sum_{n=0}^{N-2} 
         \mathbb{E}
         \bigg[
         \int_{t_n}^{t_{n+1}}
         \int_{t_n}^s
         D^2 u
         \big(
         1-s,\tilde{Y}_r
         \big)
         \Big(
         \tilde{f}_{\beta}
         \big(
         \tilde{Y}_{t_n},t_n
         \big), 
         \mathrm{\, d} W_r
         \Big)
         \mathrm{\, d} s 
         \bigg]
         \notag\\
         & \quad
         \underbrace{
         \quad +
         \mathbb{E}
         \bigg[
         \int_{1-h}^{1-h^2}
         \int_{1-h}^s
         D^2 u
         \big(
         1-s,\tilde{Y}_r
         \big)
         \Big(
         \tilde{f}_{\beta}
         \big(
         \tilde{Y}_{t_{N-1}},t_{N-1}
         \big), 
         \mathrm{\, d} W_r
         \Big)
         \mathrm{\, d} s 
         \bigg]
         \Bigg|}
         _{=:J_{1,1,1}^{(3)}}.
    \end{align}
It immediately follows from the property of the It\^o integral that
$J_{1,1,1}^{(3)}=0$. 
Using the H\"older inequality, Lemma \ref{lem:regularization-estimates-u} and Proposition \ref{prop:tilde-Df-bound}, one can easily get 
\begin{align}
\label{eq:J111-1-est}
    J_{1,1,1}^{(1)}
    & \leq
    C  h
    \bigg(
    \sum_{n=0}^{N-2} 
    \int_{t_n}^{t_{n+1}}
    \frac{1}
    {1-s}
    \mathrm{\, d} s
    +
    \int_{1-h}^{1-h^2}
    \frac{1}{1-s}
    \mathrm{\, d} s
    \bigg)
    \sup_{n \in [N-2]_0}
    \big\|
    \tilde{f}_{\beta}
    (\tilde{Y}_{t_n},t_n)
    \big\|^2_{L^2(\Omega,\mathbb{R}^d)}
    \nonumber\\
    & \leq 
    C h| \ln h|,
\end{align}
where we noted
\begin{align}
    &\sum_{n=0}^{N-2} 
    \int_{t_n}^{t_{n+1}}
    \frac{1}{1-s}
    \mathrm{\, d} s
    +
    \int_{1-h}^{1-h^2}
    \frac{1}{1-s}
    \mathrm{\, d} s
    =
    \int_0^{1-h^2}
    \frac{1}{1-s}
    \mathrm{\, d} s
    %\nonumber\\
%    =
%    \int_{h^2}^{1}
%%    \frac{1}{t}
%    \mathrm{\, d} t
    %\nonumber\\
    \leq
    2 |\ln h|.
\end{align}
To estimate $J_{1,1,1}^{(2)}$, we split
$[0,1-h^2]$ into $[0,1/2] $ and $[1/2,1-h^2]$.
On the second subinterval, the Malliavin derivative along a direction the $\big\{\tilde{Y}_r\big\}_{r \in[\frac{1}{2}, 1-h^2]}$ has a positive lower bound.
This lower bound is the key to the subsequent argument.
As in Subsection \ref{subsec:intr-mal}, we choose $N_1\in\mathbb{N}$ satisfying $0<N_1 h<\frac{1}{2}\leq (N_1+1)h<1$.
And we can rewrite $J_{1,1,1}^{(2)}$ as follows
\begin{align}
        J_{1,1,1}^{(2)}
        & \leq
        \underbrace{
        \frac{\beta}{2}
        \Bigg|
        \sum_{n=0}^{N_1}
        \sum_{j=1}^{d}
        \mathbb{E}
        \bigg[
        \int_{t_n}^{t_{n+1}}
        \int_{t_n}^s
        D^3 u
        \big(
        1-s,\tilde{Y}_r
        \big)
        \Big(
        e_j,e_j,
        \tilde{f}_{\beta}
        (\tilde{Y}_{t_n},
        t_n)
        \Big)
        \mathrm{\, d} r
        \mathrm{\, d} s 
        \bigg]
        \Bigg|}
        _{=:T_1}
        \nonumber\\
        & \quad +
        \frac{\beta}{2}
        \Bigg|
        \sum_{n=N_1+1}^{N-2}
        \sum_{j=1}^{d}
        \mathbb{E}
        \bigg[
        \int_{t_n}^{t_{n+1}}
        \int_{t_n}^s
        D^3 u
        \big(1-s,\tilde{Y}_r
        \big)
        \Big(
        e_j,e_j,
        \tilde{f}_{\beta}
        (\tilde{Y}_{t_n},
        t_n)
        \Big)
        \mathrm{\, d} r
        \mathrm{\, d} s 
        \bigg]
        \nonumber\\
        & \quad
        \underbrace{
        \quad +
        \sum_{j=1}^{d}
         \mathbb{E}
         \bigg[
         \int_{1-h}^{1-h^2}
         \int_{1-h}^s
         D^3 u
         \big(1-s,\tilde{Y}_r
         \big)
         \Big(
         e_j,e_j,
         \tilde{f}_{\beta}
        (\tilde{Y}_{t_{N-1}},t_{N-1})
         \Big)
         \mathrm{\, d} r
         \mathrm{\, d} s 
         \bigg]
         \Bigg|}
         _{=:T_2}.
\end{align}
Applying the H\"older inequality, Lemma \ref{lem:regularization-estimates-u} and Proposition \ref{prop:tilde-Df-bound}, one can easily get the estimate of $T_1$ as follows, for $h \in (0,1/4]$,
\begin{align}
\label{eq:T1-est}
        T_1
        & \leq
        C h 
        \sum_{n=0}^{N_1}
        \sum_{j=1}^d
        \int_{t_n}^{t_{n+1}}
        \frac{1}{(1-s)^{3/2}}
        \big\|
        \|e_j\|^2
        \big\|
        \tilde{f}_{\beta}
        \big(
        \tilde{Y}_{t_n},t_n
        \big)
        \big\|
        \big\|
        _{L^1(\Omega,
        \mathbb{R})}
        \mathrm{\, d} s
        \nonumber\\
        & \leq
        C d h 
        \int_0^{\frac{3}{4}}
        \frac{1}
        {(1-s)^{3/2}}
        \mathrm{\, d} s
        \sup_{n \in [N_1]_0}
        \big\|
        \tilde{f}_{\beta}
        (\tilde{Y}_{t_n},t_n)
        \big\|_{L^2(\Omega,\mathbb{R}^d)}
         \nonumber\\
         & \leq
        C d h,
\end{align}
where
$\int_0^{\frac{3}{4}}
\frac{1}
{(1-s)^{3/2}}
\mathrm{\, d} s
=2$.
Regarding $T_2$, we need to use Malliavin calculus.
Now we first use the chain rule for Malliavin derivative \eqref{eq:Chain-rule-for-Malliavin-derivative} and \eqref{eq:numerical-Malliavin-estimation} to derive, for any $j, \ell \in$ $[d]$,
\begin{equation}
\label{eq:Malliavin-D^2U}
    \begin{aligned}
        \mathscr{D}^{h_{\ell}}
        \big(
        D^2 u
        (1-s,\widetilde{Y_r}) (e_j,e_j) 
        \big)
        & =
        D^3 u
        (1-s, \widetilde{Y_r})
        \big(
        e_j, e_j,
        \mathscr{D}^{h_{\ell}} \widetilde{Y_r}
        \big) \\
        & =
        \frac{r-t_m}
        {t_{n}-t_m} 
        D^3 u
        (1-s, \widetilde{Y_r})
        (e_j, e_j,e_{\ell}),
    \end{aligned}
\end{equation}
where $h_{\ell}$ is defined by \eqref{eq:defination-direction-h} and $\left\{e_{\ell}\right\}_{{\ell} \in[d]}$ is the orthonormal basis of $\mathbb{R}^d$. 
Further, it follows from \eqref{eq:distance-t_n-t_m} that,
\begin{equation}
\label{eq:Malliavin-lower-bound}
    \frac{r-t_m}{t_{n}-t_m}
    \geq
    r-t_m
    \geq
    t_n-t_m
    \geq
    \frac{1}{2}.
\end{equation}
By Definition \ref{def:Derivative-operator}, the integration by parts formula for Malliavin derivative, i.e. Definition \ref{def:divergence-operator}, \eqref{eq:Malliavin-D^2U} and \eqref{eq:Malliavin-lower-bound}, we rewrite $T_2$ as, and for $\tilde{f}_{\beta}(\cdot,\cdot)=([\tilde{f}_{\beta}(\cdot,\cdot)]_1,\cdots,[\tilde{f}_{\beta}(\cdot,\cdot)]_d)^T$,
\begin{align}
\label{eq:T_2-expression}
    T_2
    & =
    \frac{\beta}{2}
    \Bigg|
    \sum_{n=N_1+1}^{N-2}
    \sum_{j=1}^{d}
    \sum_{{\ell}=1}^{d}
    \mathbb{E}
    \bigg[
    \int_{t_n}^{t_{n+1}}
    \int_{t_n}^s
    D^3 u
    \big(
    1-s,\tilde{Y}_r
    \big)
    \big(
    e_j,e_j,e_{\ell}
    \big)
    \cdot
    \big[
    \tilde{f}_{\beta}
    (\tilde{Y}_{t_n},t_n)
    \big]_{\ell}
    \mathrm{\, d} r
    \mathrm{\, d} s 
    \bigg] 
    \notag \\
    & \quad +
    \sum_{j=1}^{d}
    \sum_{{\ell}=1}^{d}
    \mathbb{E}
    \bigg[
    \int_{1-h}^{1-h^2}
    \int_{1-h}^s
    D^3 u
    \big(
    1-s,\tilde{Y}_r
    \big)
    \big(
    e_j,e_j,e_{\ell}
    \big)
    \cdot
    \big[
    \tilde{f}_{\beta}
    (\tilde{Y}_{t_{N-1}},t_{N-1})
    \big]_{\ell}
    \mathrm{\, d} r
    \mathrm{\, d} s 
    \bigg]
    \Bigg| 
    \notag \\
    & \leq
    \beta
    \Bigg|
    \sum_{n=N_1+1}^{N-2}
    \sum_{j=1}^{d}
    \sum_{{\ell}=1}^{d}
    \mathbb{E}
    \bigg[
    \int_{t_n}^{t_{n+1}}
    \int_{t_n}^s
    \mathscr{D}^{h_{\ell}}
    \big(
    D^2 u
    (1-s,\tilde{Y}_r)
    (e_j,e_j)
    \big)
    \cdot
    \big[
    \tilde{f}_{\beta}
    (\tilde{Y}_{t_n},t_n)
    \big]_{\ell}
    \mathrm{\, d} r
    \mathrm{\, d} s 
    \bigg] 
    \notag \\
    & \quad +
    \sum_{j=1}^{d}
    \sum_{{\ell}=1}^{d}
    \mathbb{E}
    \bigg[
    \int_{1-h}^{1-h^2}
    \int_{1-h}^s
    \mathscr{D}^{h_{\ell}}
    \big(
    D^2 u
    (1-s,\tilde{Y}_r)
    ( e_j,e_j)
    \big)
    \cdot
    \big[
    \tilde{f}_{\beta}
    (\tilde{Y}_{t_{N-1}},t_{N-1})
    \big]_{\ell}
    \mathrm{\, d} r
    \mathrm{\, d} s 
    \bigg]
    \Bigg| \notag \\
    & =
    \beta
    \Bigg|
    \sum_{n=N_1+1}^{N-2}
    \sum_{j=1}^{d}
    \sum_{{\ell}=1}^{d}
    \mathbb{E}
    \bigg[
    \int_{t_n}^{t_{n+1}}
    \int_{t_n}^s
    \int_{t_m}^{t_{n+1}}
    \big \langle
    h_{\ell}(\tilde{r}),
    \mathscr{D}
    _{\tilde{r}}
    \big(
    D^2 u
    (1-s,\tilde{Y}_r)
    (e_j,e_j)
    \big \rangle
    \cdot
    \big[
    \tilde{f}_{\beta}
    (\tilde{Y}_{t_n},t_n)
    \big]_{\ell}
    \mathrm{\, d} \tilde{r}
    \mathrm{\, d} r
    \mathrm{\, d} s 
    \bigg] 
    \notag \\
    & \quad +
    \sum_{j=1}^{d}
    \sum_{{\ell}=1}^{d}
    \mathbb{E}
    \bigg[
    \int_{1-h}^{1-h^2}
    \int_{1-h}^s
    \int_{t_m}^{t_{n+1}}
    \big \langle
    h_{\ell}(\tilde{r}),
    \mathscr{D}
    _{\tilde{r}}
    \big(
    D^2 u
    (1-s,\tilde{Y}_r)
    (e_j,e_j)
    \big \rangle
    \cdot
    \big[
    \tilde{f}_{\beta}
    (\tilde{Y}_{t_{N-1}},t_{N-1})
    \big]_{\ell}
    \mathrm{\, d} \tilde{r}
    \mathrm{\, d} r
    \mathrm{\, d} s 
    \bigg]
    \Bigg| 
    \notag \\
    & =
    \beta
    \Bigg|
    \sum_{n=N_1+1}^{N-2}
    \sum_{j=1}^{d}
    \sum_{{\ell}=1}^{d}
    \mathbb{E}
    \bigg[
    \int_{t_n}^{t_{n+1}}
    \int_{t_n}^s
    D^2 u
    (1-s,\tilde{Y}_r)
    (e_j,e_j)
    \cdot
    \delta
    \big(
    h_{\ell}
    \big[
    \tilde{f}_{\beta}
    (\tilde{Y}_{t_n},t_n)
    \big]_{\ell}
    \big)
    \mathrm{\, d} r
    \mathrm{\, d} s 
    \bigg] 
    \notag \\
    & \quad +
    \sum_{j=1}^{d}
    \sum_{{\ell}=1}^{d}
    \mathbb{E}
    \bigg[
    \int_{1-h}^{1-h^2}
    \int_{1-h}^s
    D^2 u
    (1-s,\tilde{Y}_r)
    (e_j,e_j)
    \cdot
    \delta
    \big(
    h_{\ell}
    \big[
    \tilde{f}_{\beta}
    (\tilde{Y}_{t_{N-1}},t_{N-1})
    \big]_{\ell}
    \big)
    \mathrm{\, d} r
    \mathrm{\, d} s 
    \bigg]
    \Bigg|,
\end{align}
where it is straightforward that $D^2 u(t, \cdot) (e_j,e_j) \in \mathbb{D}^{1,2}$ and $h_{\ell}[\tilde{f}_{\beta}(\tilde{Y}_{t_n},t_n)]_\ell \in $ Dom $\delta$ due to Proposition \ref{prop:f-in-C^3}, Proposition \ref{prop:tilde-Df-bound}, and Lemma \ref{lem:regularization-estimates-u}.
Using the product rule for divergence \eqref{eq:Product-rule-for-divergence} yields 
\begin{equation}
\label{eq:delta-divergence}
    \begin{aligned}
         \delta
         \Big(
         h_{\ell}
         \big[
         \tilde{f}_{\beta}
         (\tilde{Y}_{t_n},
         t_n)
         \big]_{\ell}
         \Big)
         & = 
         \big[
         \tilde{f}_{\beta}
         (\tilde{Y}_{t_n},
         t_n)
         \big]_{\ell} 
         \delta(h_{\ell})
         -
         \mathscr{D}
         ^{h_{\ell}}
         \big[
         \tilde{f}_{\beta}
         (\tilde{Y}_{t_n},
         t_n)
         \big]_{\ell} .
    \end{aligned}
\end{equation}
Plugging this equation into \eqref{eq:T_2-expression} shows 
\begin{align}
        T_2
        & \leq
        \beta
        \Bigg|
        \sum_{n=N_1+1}^{N-2}
        \sum_{j=1}^{d}
        \sum_{{\ell}=1}^{d}
        \mathbb{E}
        \bigg[
         \int_{t_n}^{t_{n+1}}
         \int_{t_n}^s
         D^2 u
         (1-s,\tilde{Y}_r)
         (e_j,e_j)
         \cdot
         \big[
         \tilde{f}_{\beta}
         (\tilde{Y}_{t_n},
         t_n)
         \big]_{\ell}
         \delta(h_{\ell})
         \mathrm{\, d} r
         \mathrm{\, d} s 
         \bigg]
         \notag\\
         & \quad
         \underbrace{ \quad +
         \sum_{j=1}^{d}
         \sum_{{\ell}=1}^{d}
         \mathbb{E}
         \bigg[
         \int_{1-h}^{1-h^2}
         \int_{1-h}^s
         D^2 u
         (1-s,\tilde{Y}_r)
         (e_j,e_j)
         \cdot
         \big[
         \tilde{f}_{\beta}
         (\tilde{Y}_{t_{N-1}},
         t_{N-1})
         \big]_{\ell}
         \delta(h_{\ell})
         \mathrm{\, d} r
         \mathrm{\, d} s 
         \bigg]
         \Bigg|}
         _{=:T_{2,1}}\notag\\
         & \quad +
         \beta
         \Bigg|
         \sum_{n=N_1+1}^{N-2}
         \sum_{j=1}^{d}
         \sum_{{\ell}=1}^{d}
         \mathbb{E}
         \bigg[
         \int_{t_n}^{t_{n+1}}
         \int_{t_n}^s
         D^2 u
         (1-s,\tilde{Y}_r)
         (e_j,e_j)
         \cdot
         \mathscr{D}^{h_{\ell}}
         \Big(
         \big[
         \tilde{f}_{\beta}
         (\tilde{Y}_{t_n},
         t_n)
         \big]_{\ell}
         \Big)
         \mathrm{\, d} r
         \mathrm{\, d} s 
         \bigg]\notag\\
         & \quad 
         \underbrace{
         \quad +
         \sum_{j=1}^{d}
         \sum_{{\ell}=1}^{d}
         \mathbb{E}
         \bigg[
         \int_{1-h}^{1-h^2}
         \int_{1-h}^s
         D^2 u
         (1-s,\tilde{Y}_r)
         (e_j,e_j)
         \cdot
         \mathscr{D}
         ^{h_{\ell}}
         \Big(
         \big[
         \tilde{f}_{\beta}
         (\tilde{Y}_{t_{N-1}},
         t_{N-1})
         \big]_{\ell}
         \Big)
         \mathrm{\, d} r
         \mathrm{\, d} s 
         \bigg]
         \Bigg|}_{=:T_{2,2}}.
\end{align}
Taking the H\"older inequality, Lemma \ref{lem:regularization-estimates-u}, Proposition \ref{prop:tilde-Df-bound} and Lemma \ref{lem:Malliavin-regularity-estimates-numerical} into account, one can easily get
\begin{align}
\label{eq:T21-est}
    T_{2,1}
%    & =
%    \beta
%        \Bigg|
%        \sum_{n=N_1+1}^{N-2}
%        \sum_{j=1}^{d}
%        \mathbb{E}
%        \bigg[
%         \int_{t_n}^{t_{n+1}}
%         \int_{t_n}^s
%         D^2 u
%         (1-s,\tilde{Y}_r)
%         (e_j,\tilde{f}_{\beta}
%         (\tilde{Y}_{t_n},
%         t_n))
%         \cdot
%         \delta(h_{\ell})
%         \mathrm{\, d} r
%         \mathrm{\, d} s 
%         \bigg]
%         \notag \\
 %        & \quad +
%         \sum_{j=1}^{d}
%         \mathbb{E}
%         \bigg[
%         \int_{1-h}^{1-h^2}
%         \int_{1-h}^s
%         D^2 u
%         (1-s,\tilde{Y}_r)
%         (e_j,\tilde{f}_{\beta}
%         (\tilde{Y}_{t_{N-1}},
%         t_{N-1}))
%%         \cdot
 %        \delta(h_{\ell})
%         \mathrm{\, d} r
%         \mathrm{\, d} s 
%         \bigg]
%         \Bigg|
%         \notag \\
    & \leq
    C d h
    \sum_{n=N_1+1}^{N-2}
    \int_{t_n}^{t_{n+1}}
    \frac{1}{1-s}
    \big\|
    \tilde{f}_{\beta}
    (\tilde{Y}_{t_n},t_n)
    \big\|
    _{L^2(\Omega,
    \mathbb{R}^d)}
    \big\|
    \delta(h_{\ell})
    \big\|_{L^2(\Omega,
    \mathbb{R})}
    \mathrm{\, d} s
    \nonumber\\
    & \quad +
    C d h
    \int_{1-h}^{1-h^2}
    \frac{1}
    {1-s}
    \big\|
    \tilde{f}_{\beta}
    (\tilde{Y}_{t_{N-1}},t_{N-1})
    \big\|_{L^2(\Omega,
    \mathbb{R}^d)}
    \big\|
    \delta(h_{\ell})
    \big\|_{L^2(\Omega,
    \mathbb{R})}
    \mathrm{\, d} s
    \nonumber\\
%    & \leq
%    C d h
%    \int_{\frac{1}{2}}^{1-h^2}
%    \frac{1}{1-s}
%    \mathrm{\, d} s
%    \sup_{n \in [N-1]_0}
%        \big\|
%        \tilde{f}_{\beta}
%        (\tilde{Y}_{t_n},t_n)
%        \big\|^2_{L^2(\Omega,\mathbb{R}^d)}
%    \nonumber\\
    & \leq
    C d  h|\mathrm{ln}h|,
\end{align}
and 
\begin{align}
\label{eq:T22-est}
    T_{2,2}
    & \leq
    C d h
    \sum_{n=N_1+1}^{N-2}
    \int_{t_n}^{t_{n+1}}
    \frac{1}
    {1-s}
    \big\|
    \mathscr{D}^{h_{\ell}}
    \big(
    \tilde{f}_{\beta}
    \big(
    \tilde{Y}_{t_n},t_n
    \big)
    \big)
    \big\|_{L^2(\Omega,
    \mathbb{R}^d)}
    \mathrm{\, d} s
    \nonumber\\
    & \quad +
    C d h
    \int_{1-h}^{1-h^2}
    \frac{1}
    {1-s}
    \big\|
    \mathscr{D}^{h_{\ell}}
    \big(
    \tilde{f}_{\beta}
    \big(
    \tilde{Y}_{t_{N-1}},t_{N-1}
    \big)
    \big)
    \big\|_{L^2(\Omega,
    \mathbb{R}^d)}
    \mathrm{\, d} s
    \nonumber\\
%    & \leq
%%    C d h 
%    \int_{\frac{1}{2}}^{1-h^2}
%    \frac{1}
%    {1-s}
%    \mathrm{\, d} s
%    \nonumber\\
    & \leq
    C d h|\ln h|.
\end{align}
Thus, combining \eqref{eq:J111-1-est}, \eqref{eq:T1-est}, \eqref{eq:T21-est} with \eqref{eq:T22-est} yields
\begin{equation}
    J_{1,1,1} 
    \leq 
    C d h |\mathrm{ln}h|.
\end{equation}
We now use Lemma \ref{lem:regularization-estimates-u}, Lemma \ref{lem:numerical-solution-monents-bounded} and Proposition \ref{prop:tilde-Df-bound} to estimate $J_{1,1,2}$ directly
\begin{equation}
\label{eq:J112-est}
    \begin{aligned}
        J_{1,1,2} 
        & \leq  
        \Bigg|
        \mathbb{E}
        \bigg[
        \int_{1-h^2}^{1}
        D u
        (1-s, \tilde{Y}_s)
        \tilde{f}_{\beta}
        (\tilde{Y}_{t_{N-1}},t_{N-1}) 
        \mathrm{\, d} s
        \bigg]
        \Bigg|\\
        & \quad +
        \Bigg|
        \mathbb{E}
        \bigg[
        \int_{1-h^2}^{1} 
         D u
         (1-s, \tilde{Y}_{t_{N-1}})
         \tilde{f}_{\beta}
         (\tilde{Y}_{t_{N-1}},t_{N-1}) 
         \mathrm{\, d} s
         \bigg]
         \Bigg|\\
         & \leq 
         C 
         \int_{1-h^2}^{1}
         \frac{1}
         {\sqrt{1-s}}
         \big\|
         \tilde{f}_{\beta}
         (\tilde{Y}_{t_{N-1}},t_{N-1})
         \big\|
         _{L^2(\Omega,
         \mathbb{R}^d)}
         \mathrm{\, d} s\\
%         \\
%         & \leq
%         C 
%         \int_{1-h^2}^{1}
%         \frac{1}{\sqrt{1-s}} 
%         \mathrm{\, d} s
%        \big\|
%        \tilde{f}_{\beta}
%        (\tilde{Y}_{t_{N-1}},t_{N-1}) 
%        \big\|^2_{L^2(\Omega,\mathbb{R}^d)}\\
          & \leq 
         C h,
    \end{aligned}
\end{equation}
where it is straightforward to show 
$
\int_{1-h^2}^1 
\tfrac{1}
{\sqrt{1-s}}
\mathrm{\, d} s
=
2h.$ 
Therefore, 
\begin{equation}
\label{eq:J_11-est}
J_{1,1}
\leq 
C d h | \ln h|.
\end{equation}
The estimate of $J_{1,2}$ can be obtained in the same way as that of $J_{1,1}$, and is therefore omitted.
Putting together the estimates \eqref{eq:J_11-est} and \eqref{eq:J_13-est}, we arrive at
\begin{align}
\label{eq:J1-est}
    J_1
    \leq 
    C d h|\mathrm{ln}h|
    +
    \tfrac{Cd}{M}.
\end{align}
Regarding $J_2$, we begin with the following decomposition:
\begin{align}
J_2
& \leq
\underbrace{
\left|
\sum_{n=0}^{N-1}
\mathbb{E}
\left[
\int_{t_n}^{t_{n+1}}
D u(1-s,\tilde Y_{t_n})
\Big(
f_\beta(\tilde Y_{t_n},t_n)
-
f_\beta(\tilde Y_s,s)
\Big)
\,\mathrm{d}s
\right]
\right|
}_{=:J_{2,1}}
\nonumber\\
&\quad +
\underbrace{
\left|
\sum_{n=0}^{N-1}
\mathbb{E}
\left[
\int_{t_n}^{t_{n+1}}
\left(
D u(1-s,\tilde Y_s)
-
D u(1-s,\tilde Y_{t_n})
\right)
\Big(
f_\beta(\tilde Y_{t_n},t_n)
-
f_\beta(\tilde Y_s,s)
\Big)
\,\mathrm{d}s
\right]
\right|
}_{=:J_{2,2}}.
\end{align}
For the estimate of $J_{2,1}$, Proposition \ref{prop:f-in-C^3} reveals that
$\partial_t f_\beta(x,t)$ has a singularity at the terminal time $t=1$, which poses a difficulty in improving the convergence rate in the subsequent error analysis. To overcome this difficulty, we divide the integration interval $[0,1]$ into $[0,1-h]$ and $[1-h,1]$, as follows:
\begin{align}
J_{2,1}
& \leq
\underbrace{
\left|
\sum_{n=0}^{N-2}
\mathbb{E}
\left[
\int_{t_n}^{t_{n+1}}
D u(1-s,\tilde Y_{t_n})
\Big(
f_\beta(\tilde Y_{t_n},t_n)
-
f_\beta(\tilde Y_s,s)
\Big)
\,\mathrm{d}s
\right]
\right|
}_{=:J_{2,1,1}}
\nonumber\\
&\quad +
\underbrace{
\left|
\mathbb{E}
\left[
\int_{1-h}^{1}
D u(1-s,\tilde Y_{t_{N-1}})
\Big(
f_\beta(\tilde Y_{t_{N-1}},t_{N-1})
-
f_\beta(\tilde Y_s,s)
\Big)
\,\mathrm{d}s
\right]
\right|
}_{=:J_{2,1,2}} .
\label{eq:J21-decomposition}
\end{align}
Applying the It\^o formula to $f_\beta(x,t)$ in $J_{2,1,1}$, we obtain
\begin{align}
\label{eq:J211-ito}
J_{2,1,1}
& \leq
\underbrace{
\left|
\sum_{n=0}^{N-2}
\mathbb{E}
\left[
\int_{t_n}^{t_{n+1}}
\int_{t_n}^{s}
D u(1-s,\tilde Y_{t_n})
\partial_r f_\beta(\tilde Y_r,r)
\,\mathrm{d}r
\,\mathrm{d}s
\right]
\right|
}_{=:J^{(1)}_{2,1,1}}
\nonumber\\
& \quad +
\underbrace{
\left|
\sum_{n=0}^{N-2}
\mathbb{E}
\left[
\int_{t_n}^{t_{n+1}}
\int_{t_n}^{s}
D u(1-s,\tilde Y_{t_n})
\Big(
D f_\beta(\tilde Y_r,r)
\tilde f_\beta(\tilde Y_{t_n},t_n)
\Big)
,\mathrm{d}r
\,\mathrm{d}s
\right]
\right|
}_{=:J^{(2)}_{2,1,1}}
\nonumber\\
&\quad+
\underbrace{
\frac{\beta}{2}
\left|
\sum_{n=0}^{N-2}
\sum_{j=1}^{d}
\mathbb{E}
\left[
\int_{t_n}^{t_{n+1}}
\int_{t_n}^{s}
D u(1-s,\widetilde Y_{t_n})
\Big(
D^2 f_\beta(\widetilde Y_r,r)(e_j,e_j)
\Big)
\,\mathrm{d}r
\,\mathrm{d}s
\right]
\right|
}_{=:J^{(3)}_{2,1,1}}.
\end{align}
In the following we cope with the above three items separately.
Now we are in position to give an estimate of $J^{(1)}_{2,1,1}$.
Once again, combining with Proposition \ref{prop:f-in-C^3} and Lemma \ref{lem:regularization-estimates-u}, we get
\begin{equation}
    \begin{aligned}
        J^{(1)}_{2,1,1}
        & \leq
        C 
        \sum_{n=0}^{N-2} 
        \int_{t_n}^{t_{n+1}} 
        \int_{t_n}^{s}
        \frac{1}
        {\sqrt{1-s}}
        \big\|
        \partial_t 
        f_{\beta}(\tilde{Y}_r,r)
        \big\|_{L^2(\Omega,\mathbb{R}^d)}
        \mathrm{\, d} r
        \mathrm{\, d} s \\
        & \leq 
        C d^{\frac{1}{2}} h
        \sum_{n=0}^{N-2} 
        \int_{t_n}^{t_{n+1}} 
        \frac{1}
        {1-s}
        \mathrm{\, d} s \\
        & \leq 
        C d^{\frac{1}{2}} h |\mathrm{ln}h|,
    \end{aligned}
\end{equation}
where we used the fact that $r \in [t_n,s]$ and 
\begin{equation}
\sum_{n=0}^{N-2} 
\int_{t_n}^{t_{n+1}}
\frac{1}{\sqrt{1-s}}
\frac{1}{\sqrt{1-r}}
\mathrm{d}s
=
\int_0^{1-h}
\frac{1}{\sqrt{1-s}}
\frac{1}{\sqrt{1-r}}
\mathrm{d}s
\leq 
\int_0^{1-h}
\frac{1}{1-s}
\mathrm{d}s
\leq
|\mathrm{ln}h|.
\end{equation}
Next, we handle the last two items. Using Proposition \ref{prop:f-in-C^3} and Proposition \ref{prop:tilde-Df-bound} yields
\begin{equation}
\begin{aligned}
J^{(2)}_{2,1,1}
& \leq 
        C 
        \sum_{n=0}^{N-2} 
        \int_{t_n}^{t_{n+1}} 
        \int_{t_n}^{s}
        \frac{1}
        {\sqrt{1-s}}
        \big\|
        Df_{\beta}
    (\tilde{Y}_r,r)
    \tilde{f}_{\beta}
    (\tilde{Y}_{{t_n}},
    t_n) \big\|_{L^2(\Omega,\mathbb{R}^d)} 
    \mathrm{\, d} r
        \mathrm{\, d} s   \\
        & \leq 
        C h
        \sum_{n=0}^{N-2} 
        \int_{t_n}^{t_{n+1}} 
        \frac{1}
        {\sqrt{1-s}}
        \mathrm{\, d} s 
        \sup_{n\in[N-1]_0}
         \big\|
    \tilde{f}_{\beta}
    (\tilde{Y}_{{t_n}},
    t_n) 
    \big\|_{L^2(\Omega,\mathbb{R}^d)} \\
    & \leq 
    C h,
\end{aligned}
\end{equation}
and 
\begin{equation}
\begin{aligned}
J^{(3)}_{2,1,1}
& \leq 
C 
        \sum_{n=0}^{N-2}
        \sum_{j=1}^d
        \int_{t_n}^{t_{n+1}} 
        \int_{t_n}^s
        \frac{1}
        {\sqrt{1-s}}
        \big\|
        D^2 f_{\beta}
        (\tilde{Y}_r,r)
        (e_j,e_j) \big\|_{L^2(\Omega,\mathbb{R}^d)}
        \mathrm{\, d} r
        \mathrm{\, d} s\\
        & \leq 
        C d h 
        \sum_{n=0}^{N-2} 
        \int_{t_n}^{t_{n+1}} 
        \frac{1}
        {\sqrt{1-s}}
        \mathrm{\, d} s \\
        & \leq 
        C d h ,
\end{aligned}
\end{equation}
where 
$\sum_{n=0}^{N-2} 
\int_{t_n}^{t_{n+1}} 
\frac{1}
{\sqrt{1-s}}
\mathrm{\, d} s
=
\int_0^{1-h}
\frac{1}
{\sqrt{1-s}}
\mathrm{\, d} s
=
\int_h^1 
\frac{1}
{\sqrt{t}}
\mathrm{\, d} t
\leq 2$.
Regarding $J_{2,1,2}$, thanks to Lemma \ref{lem:regularization-estimates-u}, Lemma \ref{lem:numerical-solution-monents-bounded}  and \eqref{eq:b-lipschitz-in-x-houder-in-t}, one can get
\begin{align}
        J_{2,1,2}
        \leq & 
        C
        \int_{1-h}^{1}
        \frac{1}
        {\sqrt{1-s}}
        \big\|
        f_{\beta}(
        \tilde{Y}_{t_{N-1}},t_{N-1})
        -
        f_{\beta}(
        \tilde{Y}_s,s)
        \big\|
        _{L^2(\Omega,
        \mathbb{R}^d)} 
        \mathrm{\, d} s 
        \nonumber\\
       \leq &
       C
       \int_{1-h}^{1}
       \frac{1}
       {\sqrt{1-s}}
       \left(
       \big\|
       \tilde{Y}_s-\tilde{Y}_{t_{N-1}} \big\|_{L^2(\Omega,
        \mathbb{R}^d)}
       +
       d^\frac{1}{2}
       (s-1+h)^\frac{1}{2} 
       \right)
       \mathrm{\, d} s
       \nonumber\\
%       \leq & 
%       C h^{\frac{1}{2}}
%       \int_{1-h}^{1}
%       \frac{1}
%       {\sqrt{1-s}} 
%       \sup_{0 \leq r \leq {t_N}}
%       \left(1+ \| \tilde{Y}_r\|_{L^2\left(\Omega,\mathbb{R}^d\right)}
%       \right) 
%       \mathrm{\, d} s\\
%       & +
%       C d^\frac{1}{2}
%       h^{\frac{1}{2}}
%       \int_{1-h}^{1}
%       \frac{1}{\sqrt{1-s}}
%       \mathrm{\, d} s\\
       \leq & 
       C d^{\frac{1}{2}}h, 
\end{align}
where 
\[
    \int_{1-h}^1 
    \frac{1}{\sqrt{1-s}} 
    \mathrm{\, d} s
    =
    \int_0^h
    \frac{1}{\sqrt{t}} 
    \mathrm{\, d} t
    \leq
    2 \sqrt{h}.
\]
Based on the above estimates, we get 
\begin{equation}
    J_{2,1} 
    \leq C  d h.
\end{equation}
For the estimate of $J_{2,2}$, to overcome the possible singularities, we also need to split the time interval $[0,1]$ into $[0,1-h]$ and $[1-h,1]$. Applying the Taylor expansion to $D u(1-s,\cdot)$ leads to, for
$v(r):=\tilde Y_{t_n}+r\big(\tilde Y_s-\tilde Y_{t_n}\big),\, r\in [0,1],$
\begin{align}
\label{eq:J22-decomposition}
J_{2,2}
&\leq
\underbrace{
\left|
\sum_{n=0}^{N-2}
\mathbb{E}
\left[
\int_{t_n}^{t_{n+1}}
\int_0^1
D^2u(1-s,v(r))
\Big(
\tilde Y_s-\tilde Y_{t_n},
f_\beta(\tilde Y_s,s)
-
f_\beta(\tilde Y_{t_n},t_n)
\Big)
\,\mathrm{d}r
\,\mathrm{d}s
\right]
\right|
}_{=:J_{2,2,1}}
\nonumber\\
&\quad+
\underbrace{
\left|
\mathbb{E}
\left[
\int_{1-h}^{1}
\left(
D u(1-s,\tilde Y_s)
-
D u(1-s,\tilde Y_{t_{N-1}})
\right)
\Big(
f_\beta(\tilde Y_s,s)
-
f_\beta(\tilde Y_{t_{N-1}},t_{N-1})
\Big)
\,\mathrm{d}s
\right]
\right|
}_{=:J_{2,2,2}}.
\end{align}
With regard to $J_{2,2,1}$ and $J_{2,2,2}$, using the same technique as before and applying the H\"older inequality, Lemma \ref{lem:regularization-estimates-u} and Lemma \ref{lem:numerical-solution-monents-bounded}, we deduce that
\begin{equation}
     \begin{aligned}
      J_{2,2,1}
      \leq & 
      C
      \sum_{n=0}^{N-2}
     \int_{t_n}^{t_{n+1}} 
      \frac{1}{1-s}
      \big\|
      \tilde{Y}_s-\tilde{Y}_{t_{n}}
      \big\|
      _{L^2(\Omega,\mathbb{R}^d)}
      \big\|
      f_{\beta}(\tilde{Y}_s,s)
      -f_{\beta}(
      \tilde{Y}_{t_n},t_n)
      \big\|
      _{L^2(\Omega,\mathbb{R}^d)}
      \mathrm{\, d} s\\
      \leq & 
      C h 
      \int_0^{1-h}
      \frac{1}
      {1-s} 
       \mathrm{\, d} s
       \sup_{0 \leq r \leq 1}
       \left(1+ \| \tilde{Y}_r\|^2_{L^2\left(\Omega,\mathbb{R}^d\right)}
        \right)\\
      \leq &
      C d h |\mathrm{ln}h|.
      \end{aligned}
\end{equation}
and
\begin{equation}
       \begin{aligned}
        J_{2,2,2}
        & \leq
        \Bigg|
        \mathbb{E}
        \bigg[
        \int_{1-h}^1 
        D u
        (1-s, \tilde{Y}_s)
        \big(
        f_{\beta}(\tilde{Y}_s,s)
        -
        f_{\beta}(\tilde{Y}_{t_{N-1}},t_{N-1})
        \big) 
        \mathrm{\, d} s
        \bigg]
        \Bigg|\\
        & \quad +
        \Bigg|
        \mathbb{E}
         \bigg[
         \int_{1-h}^1 
         D u
         (1-s, \tilde{Y}_{t_{N-1}})
         \big(
         f_{\beta}(\tilde{Y}_s,s)
         -
         f_{\beta}(
         \tilde{Y}_{t_{N-1}},t_{N-1})
         \big) 
         \mathrm{\, d} s
         \bigg]
         \Bigg|\\
       & \leq 
       C \int_{1-h}^1 
       \frac{1}
       {\sqrt{1-s}}
       \big\|
       f_{\beta}(\tilde{Y}_s,s)
       -
       f_{\beta}(
       \tilde{Y}_{t_{N-1}},t_{N-1})
       \big\|
       _{L^2(\Omega,\mathbb{R}^d)}
       \mathrm{\, d} s \\
       & \leq 
       C d^{\frac{1}{2}} h.
      \end{aligned}
\end{equation}
To conclude
\begin{equation}
\label{eq:J2-est}
    J_{2}
    \leq 
    C d 
    h|\mathrm{ln}h|.
\end{equation}
Plugging the estimates \eqref{eq:J1-est} and \eqref{eq:J2-est} into \eqref{eq:tv-weak-error} yields
\begin{equation}
 \Big|
        \mathbb{E}
        \big[
        \phi
        \big(
        \tilde{Y}_1
        \big)
        \big]
        -
        \mathbb{E}
        \big[
        \phi
        \big(X_1\big)
        \big]
        \Big| 
        \leq 
        C d 
    h|\mathrm{ln}h|
    +
    \tfrac{Cd}{M}.
\end{equation}
The proof of Theorem \ref{thm-main-rerult-w0-w1} is completed.

\end{proof}

At this stage, we are well-prepared to derive the convergence rate of the SFS algorithm with exact drift \eqref{eq:Euler-scheme-no-mc}.

\begin{proof}
\textbf{Proof of Theorem \ref{thm:exact-f}}.
When the drift is available in closed form, we set
\(\widetilde f_\beta=f_\beta\), and the scheme
\eqref{eq:e-m-mc} reduces to the SFS algorithm
\eqref{eq:Euler-scheme-no-mc}.
According to the proof of Theorem \ref{thm-main-rerult-w0-w1}, the weak error can be decomposed by using the Kolmogorov equation as in \eqref{eq:weak-error-decomposition}. 
In this case, the Monte Carlo error term vanishes, namely,
\[
J_1=0.
\]
Consequently, the total weak error is bounded only by the time discretization term $J_2$. By the estimate of $J_2$ obtained in \eqref{eq:J2-est}, we have
\[
\Big|
        \mathbb{E}
        \big[
        \phi
        ({Y}_1)
        \big]
        -
        \mathbb{E}
        \big[
        \phi
        (X_1)
        \big]
        \Big| 
        \leq  
        Cdh |\ln h|.
\]
The proof of Theorem \ref{thm:exact-f} is completed.
\end{proof}

\section{Numerical experiments}
\label{section:Numerical-experiment}
\noindent
In this section, we conduct some numerical experiments to validate our theory and investigate the performance of SFS.
We consider sampling from some two-dimensional and high-dimensional Gaussian mixture distributions and other classic examples.

\subsection{Two-dimensional distribution}
\begin{figure}[htbp]
    \centering
    \begin{subfigure}{0.4\linewidth}  
        \includegraphics[width=\linewidth]{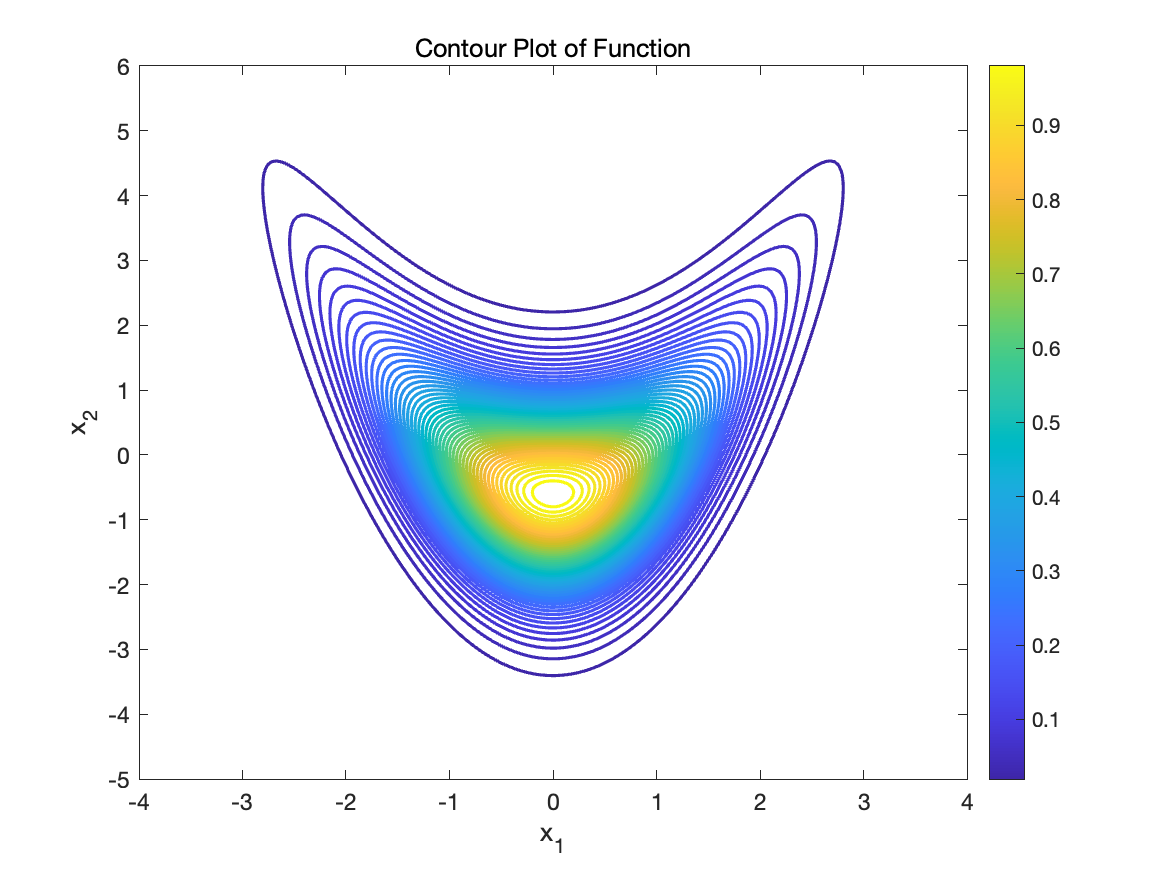} 
        \caption{True contour plot}	
    \end{subfigure}
    \hspace{-0.5cm} 
    \begin{subfigure}{0.4\linewidth}
        \includegraphics[width=\linewidth]{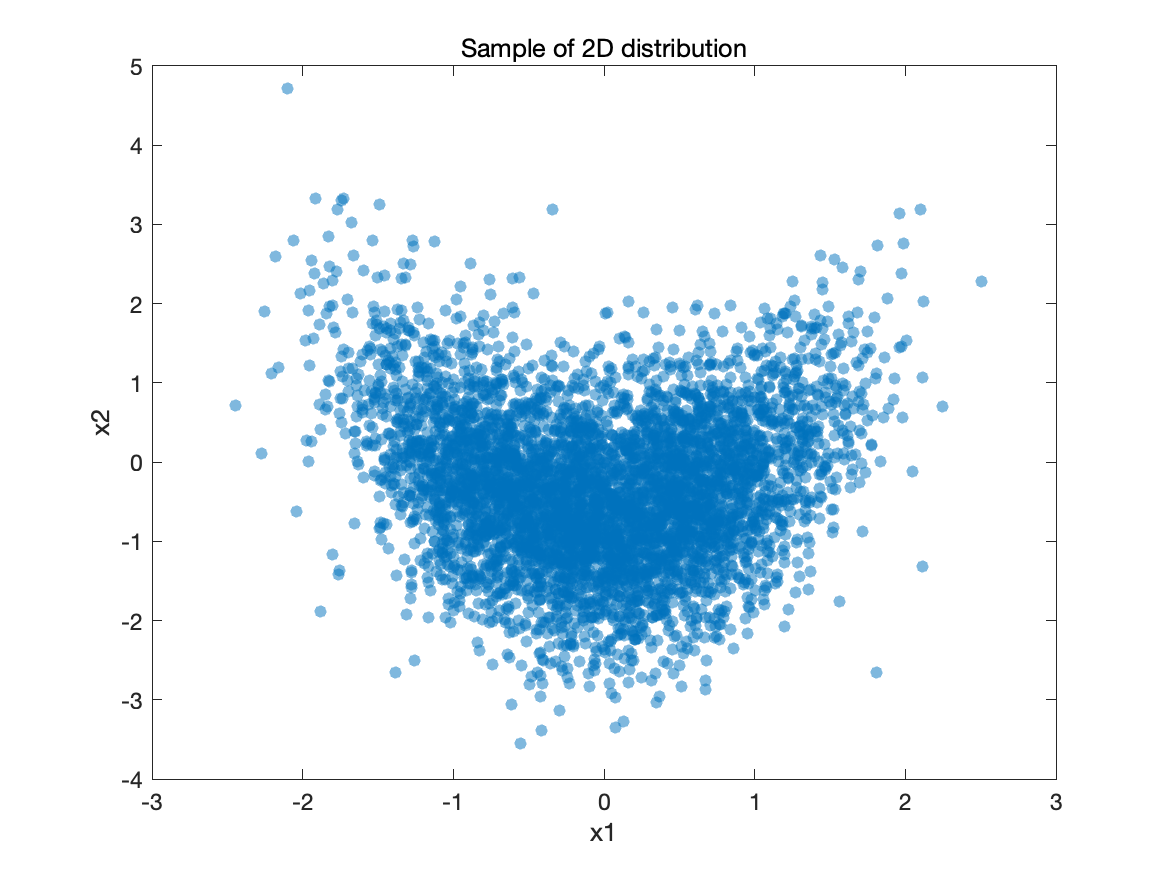}
        \caption{Sample point via SFS }	
    \end{subfigure}
    \caption{Sampling from $p_1(x)$ with SFS.}
    \label{figure:two_dim}
\end{figure}

In this subsection, we consider sampling from the following two-dimensional distributions with probability density functions (PDF):
\begin{equation}
    p(x)
    :=
    c\,
    \mathrm{e}^{
    -(
    x_1^2+
    (
    x_2-3/5(x_1^2-1)
    )^2)/2},
    \quad 
    (x_1,x_2)\in \mathbb{R}^2,
\end{equation}
where $c$ is normalizing constant.
We employ the SFS algorithm \eqref{eq:e-m-mc} with $\beta = 2$ and step size $h = 2^{-7}$ to generate $5000$ sample points, which are visualized as two-dimensional scatter plots in Figure \ref{figure:two_dim}.
The left picture presents the true contour plots, while the right one depicts the scatter plots of sampled points.
Evidently, the SFS gives a good performance in sampling from complex distributions.

\subsection{Double-well potential}
\begin{figure}[htbp]
    \centering
    \begin{subfigure}{0.4\linewidth}  
        \includegraphics[width=\linewidth]{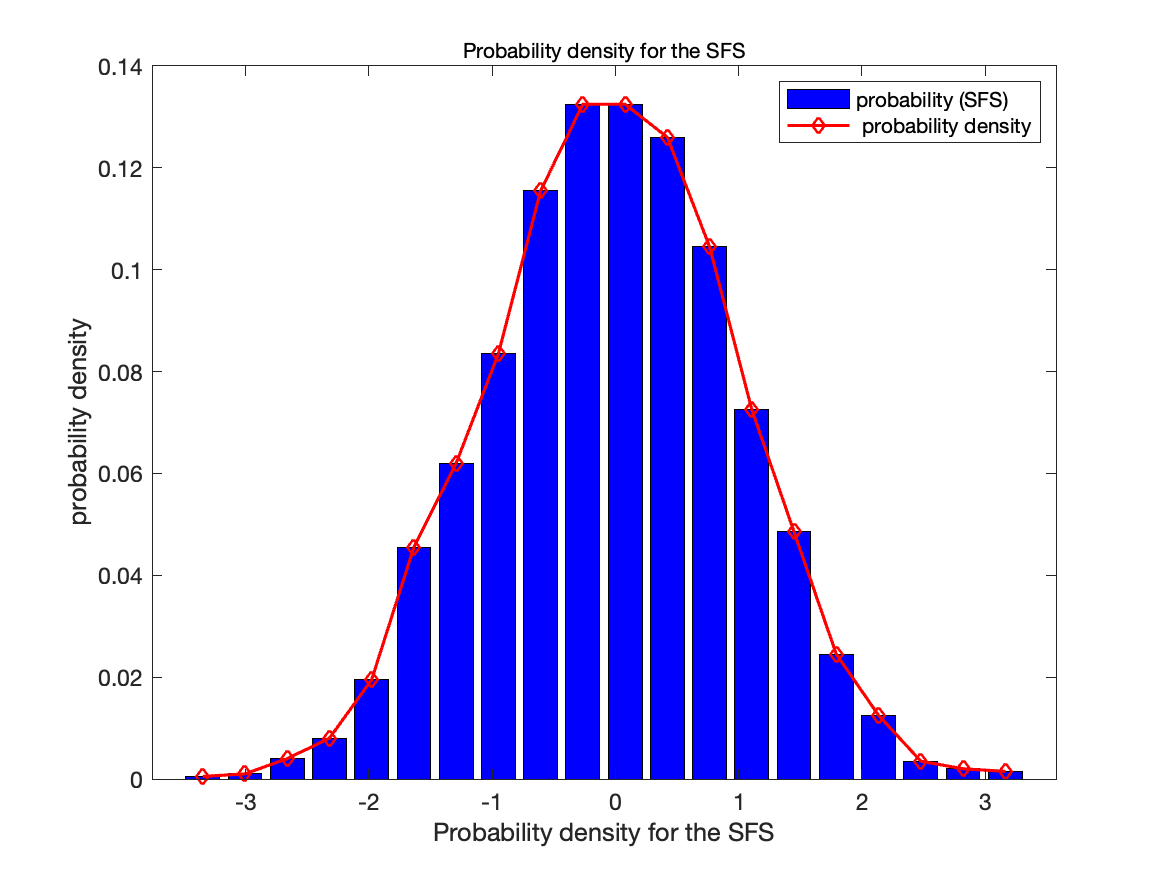} 
        \caption{SFS}	
    \end{subfigure}
    \hspace{-0.5cm} 
    \begin{subfigure}{0.4\linewidth}
        \includegraphics[width=\linewidth]{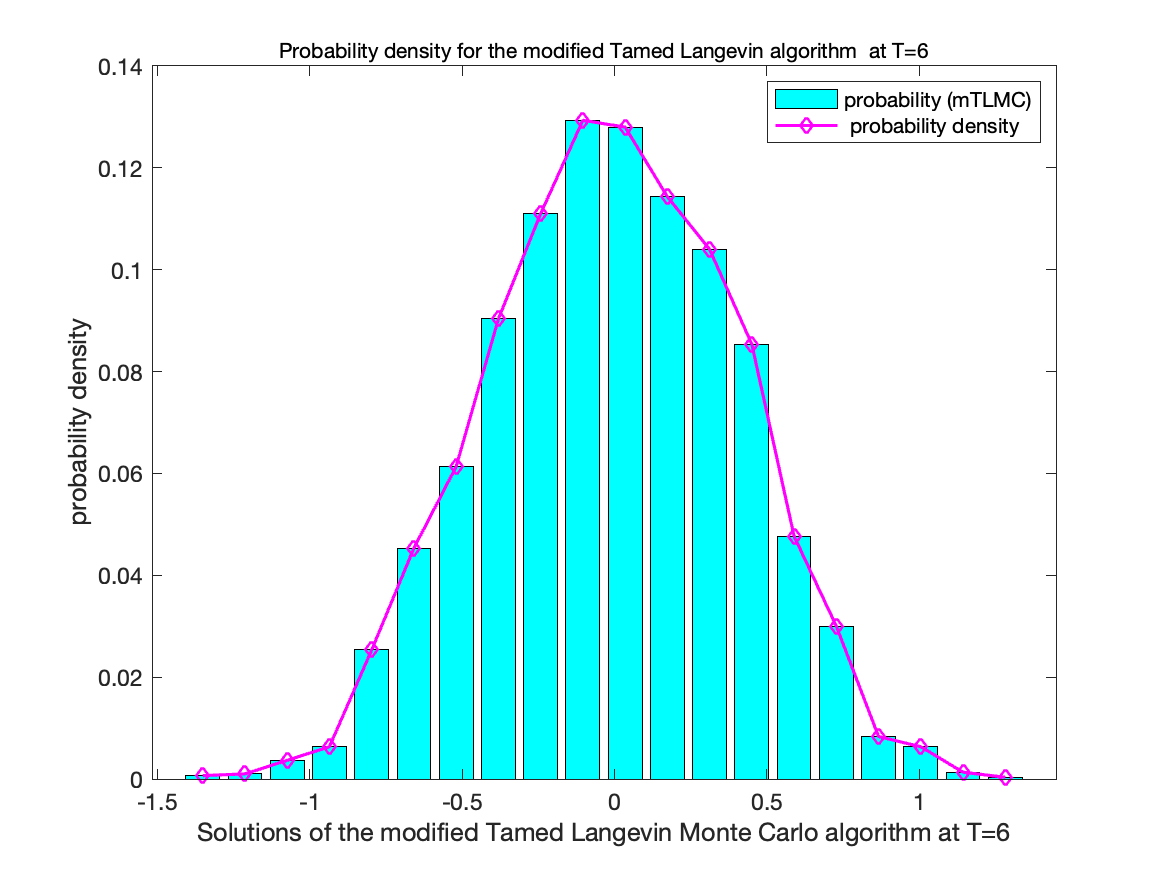}
        \caption{MTLMC}	
    \end{subfigure}
    \caption{Probability density of the first component for the Double-well potential.}
    \label{figure:dw_dens}
\end{figure}

In this subsection, we consider sampling from the following high-dimensional distribution with its potential given by:
\begin{equation}
    V(x)
    :=
    \tfrac{1}{4}
    \|x\|^4
    -
    \tfrac{1}{2}
    \|x\|^2,
    \quad x\in \mathbb{R}^d.
\end{equation}
For $d=100$ and a fixed step-size $h=2^{-7}$, we simulate $2000$ independent samples using the SFS \eqref{eq:e-m-mc} with $\beta=1$.
In Figure \ref{figure:dw_dens}, we plot the density curves of the first components of SFS and modified tamed Langevin Monte Carlo (MTLMC) algorithm proposed by \cite{NNZ2025} as a reference. 
As illustrated by Figure \ref{figure:dw_dens}, the normalized histogram plots of samples as well as the marginal probability density
curves generated by the SFS algorithm \eqref{eq:e-m-mc} and the MTLMC algorithm are very close.

\subsection{ Gaussian mixture distributions}
\label{subsec:Mog}
In this subsection, we focus on sampling from the following high-dimensional  Gaussian mixture distributions,
\begin{equation}
    p_2(x) 
    =
    \sum_{i=1}^\kappa
    \theta_i
    \,
    \mathcal{N}(x;\alpha_i, \Sigma_i),
    \quad \sum_{i=1}^\kappa \theta_i=1 
    \text { and } 
    0 \leq \theta_i \leq 1,
    \quad \kappa \in \mathbb{N},
   % \quad  \alpha_i=4(\sin (2(i-1) \pi / 4), \cos (2(i-1) \pi / 4)), \, \Sigma_i=\tfrac{1}{4}\cdot I_2.
\end{equation}
where $\alpha_i\in \mathbb{R}^d$ and $\Sigma_i\in\mathbb{R}^{d\times d}$ are the mean and covariance matrix of the Gaussian component, respectively.

\subsubsection{Sampling performance for Gaussian mixtures distributions}

\begin{figure}[htbp]
    \centering
    \begin{subfigure}{0.4\linewidth}
        \includegraphics[width=\linewidth]{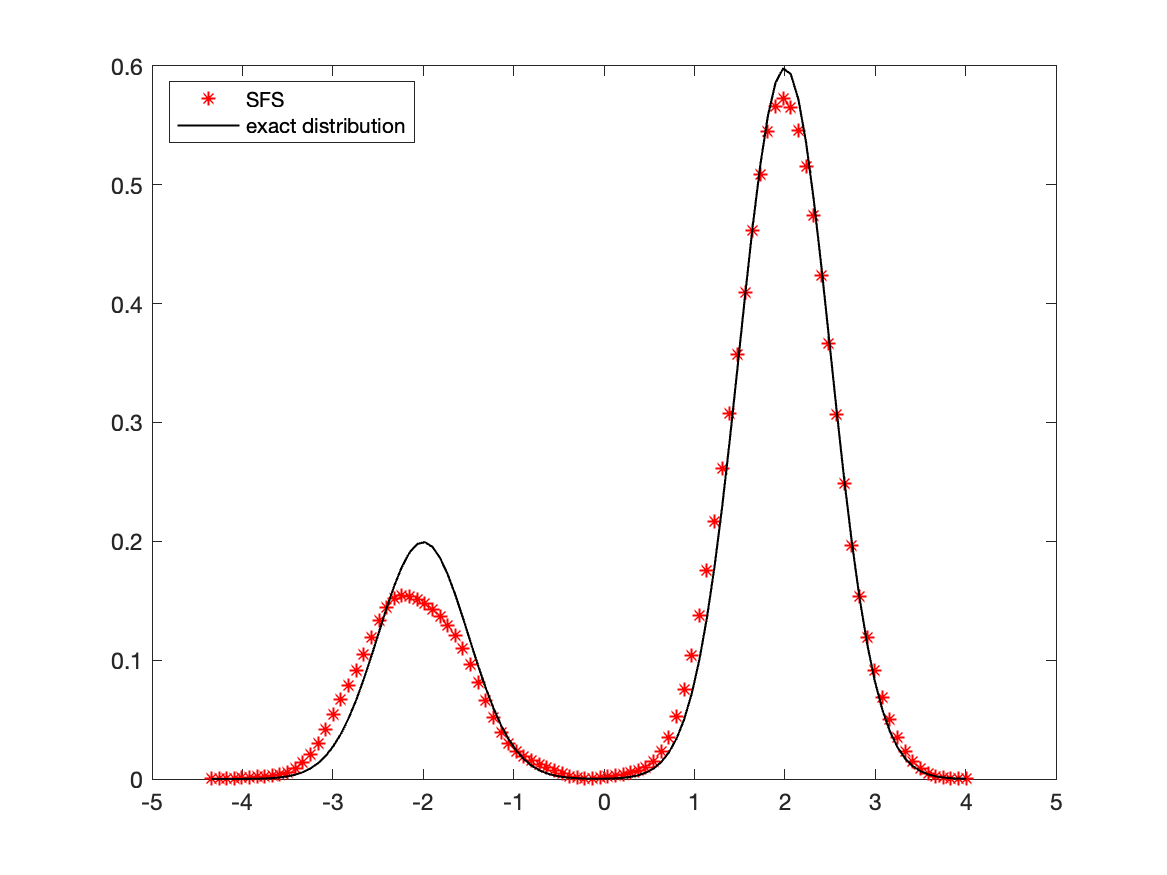}
        \caption{$\alpha_1=-2\cdot\mathbf{1}_d,\,\alpha_2=2\cdot\mathbf{1}_d,\,d=100$}	
    \end{subfigure}
    \hspace{-0.5cm}
    \begin{subfigure}{0.4\linewidth}  
        \includegraphics[width=\linewidth]{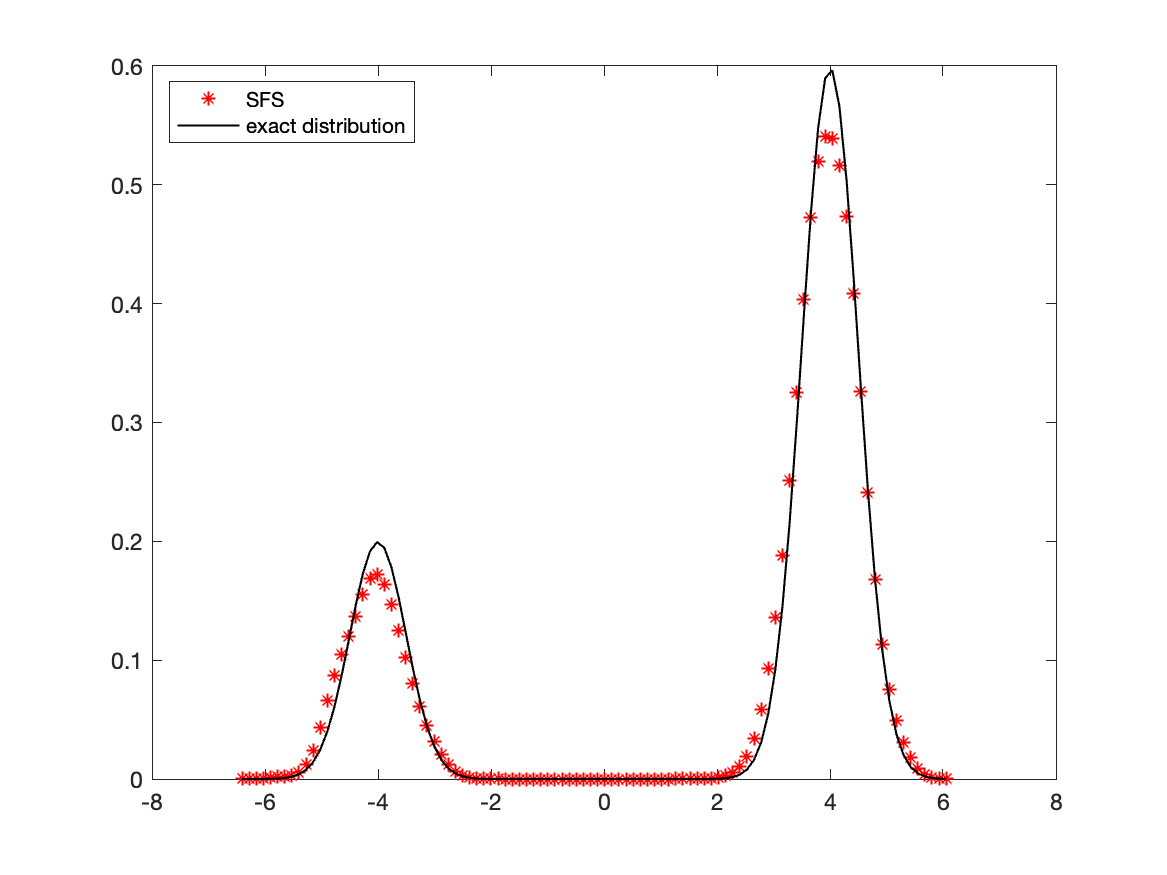} 
        \caption{$\alpha_1=-4\cdot\mathbf{1}_d,\,\alpha_2=4\cdot\mathbf{1}_d,\,d=40$}	
    \end{subfigure}
    \hspace{-0.5cm} 
    \caption{Probability density of the first component for the  Gaussian mixture distributions.}
    \label{figure:mog_dens_2}
\end{figure}

\textbf{Sampling via SFS with different dimensions $d$:}
Given $\kappa=2, \,\theta_1=\tfrac{1}{4},\,\theta_2=\tfrac{3}{4},\, \Sigma_1=\Sigma_2=\tfrac{1}{4}\cdot\mathbf{I}_d$ and a uniform time step-size $h=10^{-3}$, we simulate $1000$ independent Markov chains using the SFS algorithm \eqref{eq:Euler-scheme-no-mc} with $\beta=1$ in different dimensions $d$.
In Figure \ref{figure:mog_dens_2}, the probability density of the first component obtained via SFS \eqref{eq:Euler-scheme-no-mc} is plotted using red stars, while the exact density curves are also plotted using black solid lines.

It is evident that, in the case $\alpha_1=-2\cdot\mathbf{1}_d,\,\alpha_2=2\cdot\mathbf{1}_d$, the SFS exhibit bimodal distributions with dimension $d=100$.
However, as the distance between the means of the Gaussian components increases (e.g., $\alpha_1=-4\cdot\mathbf{1}_d, \,\alpha_2=4\cdot\mathbf{1}_d$), only SFS samples with dimension $d = 40$ are able to accurately recover the underlying bimodal target distribution. 
At higher dimensions, the SFS with $\beta=1$ fails to capture the two-mode structure effectively.
 
\begin{figure}[htbp]
    \centering
    \begin{subfigure}{0.4\linewidth}  
        \includegraphics[width=\linewidth]{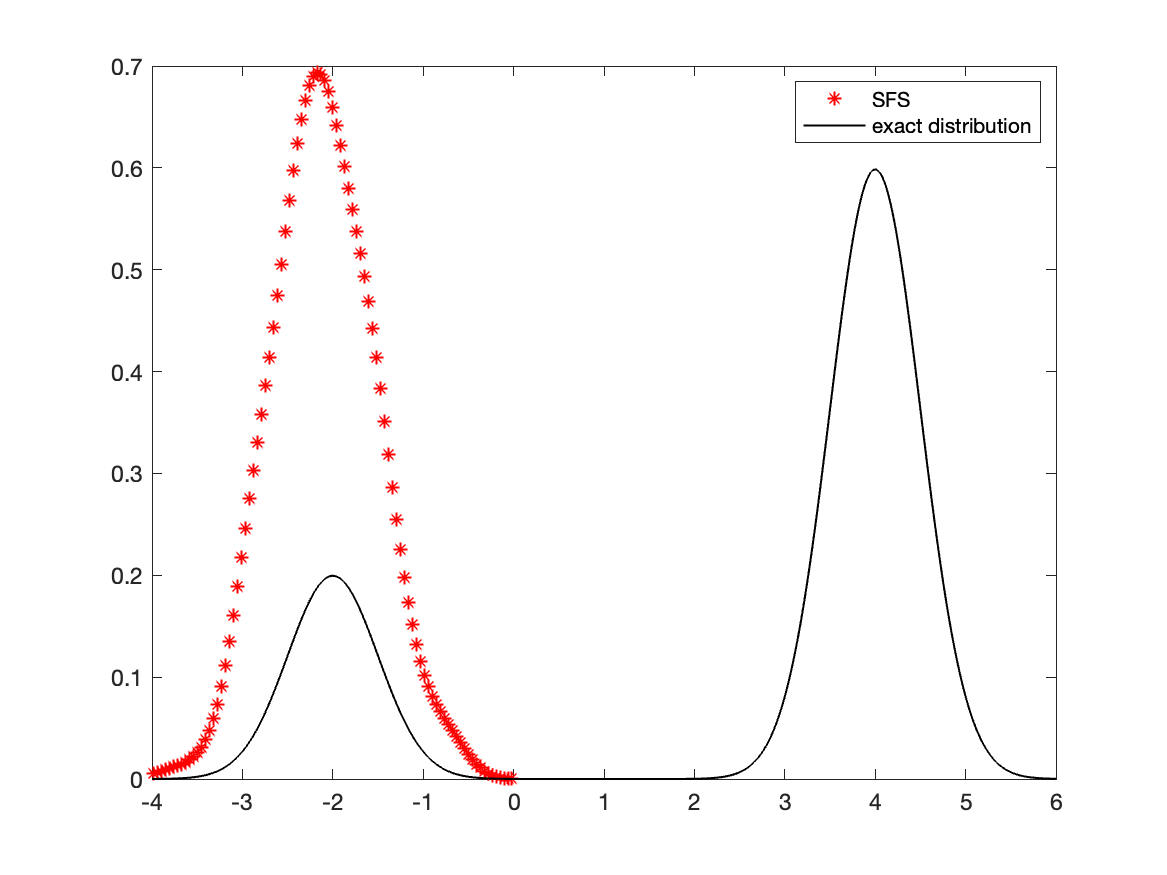} 
        \caption{$\beta=1$}	
    \end{subfigure}
    \hspace{-0.5cm} 
    \begin{subfigure}{0.4\linewidth}
        \includegraphics[width=\linewidth]{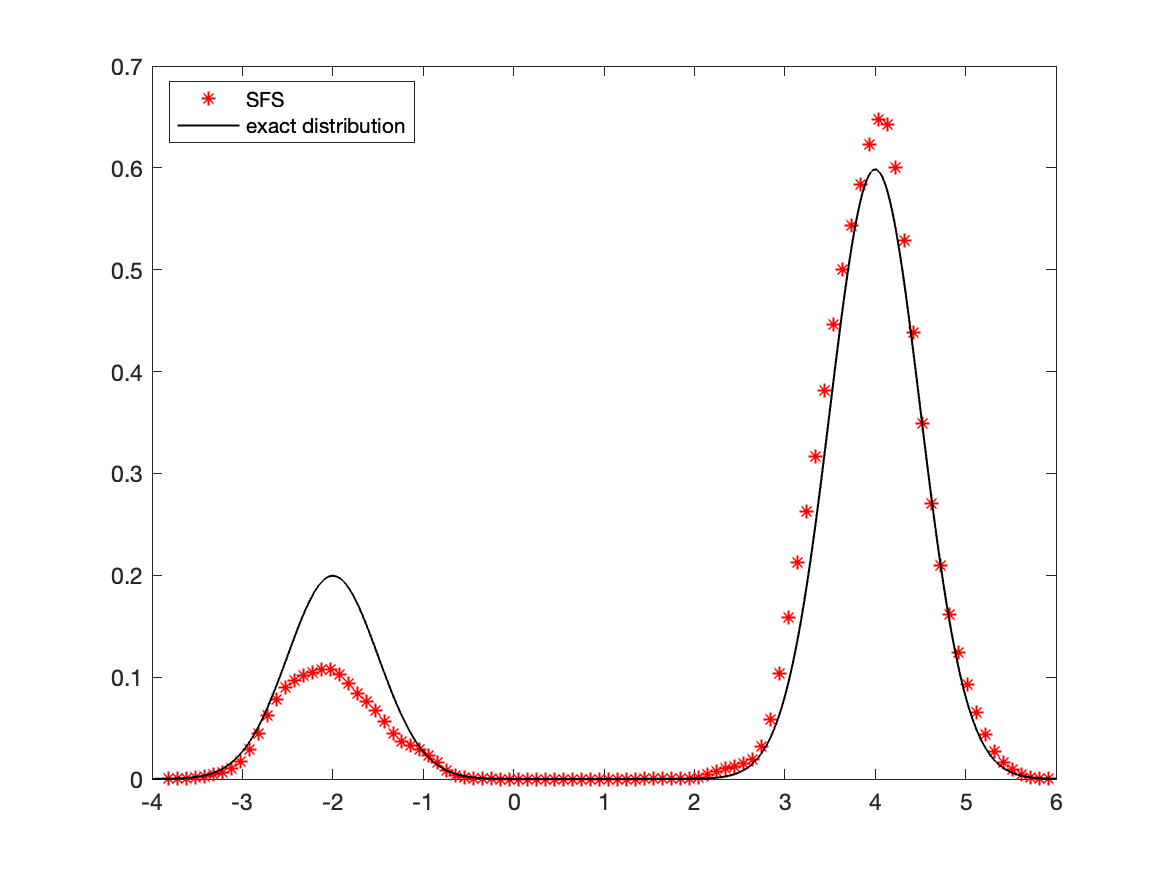}
        \caption{$\beta=5$}	
    \end{subfigure}
    \caption{Probability density of the first component for the  Gaussian mixture distributions.}
    \label{figure:mog_dens_kappa_2_different_d}
\end{figure}

\begin{figure}[htbp]
    \centering
    \begin{subfigure}{0.32\linewidth}  
        \includegraphics[width=\linewidth]{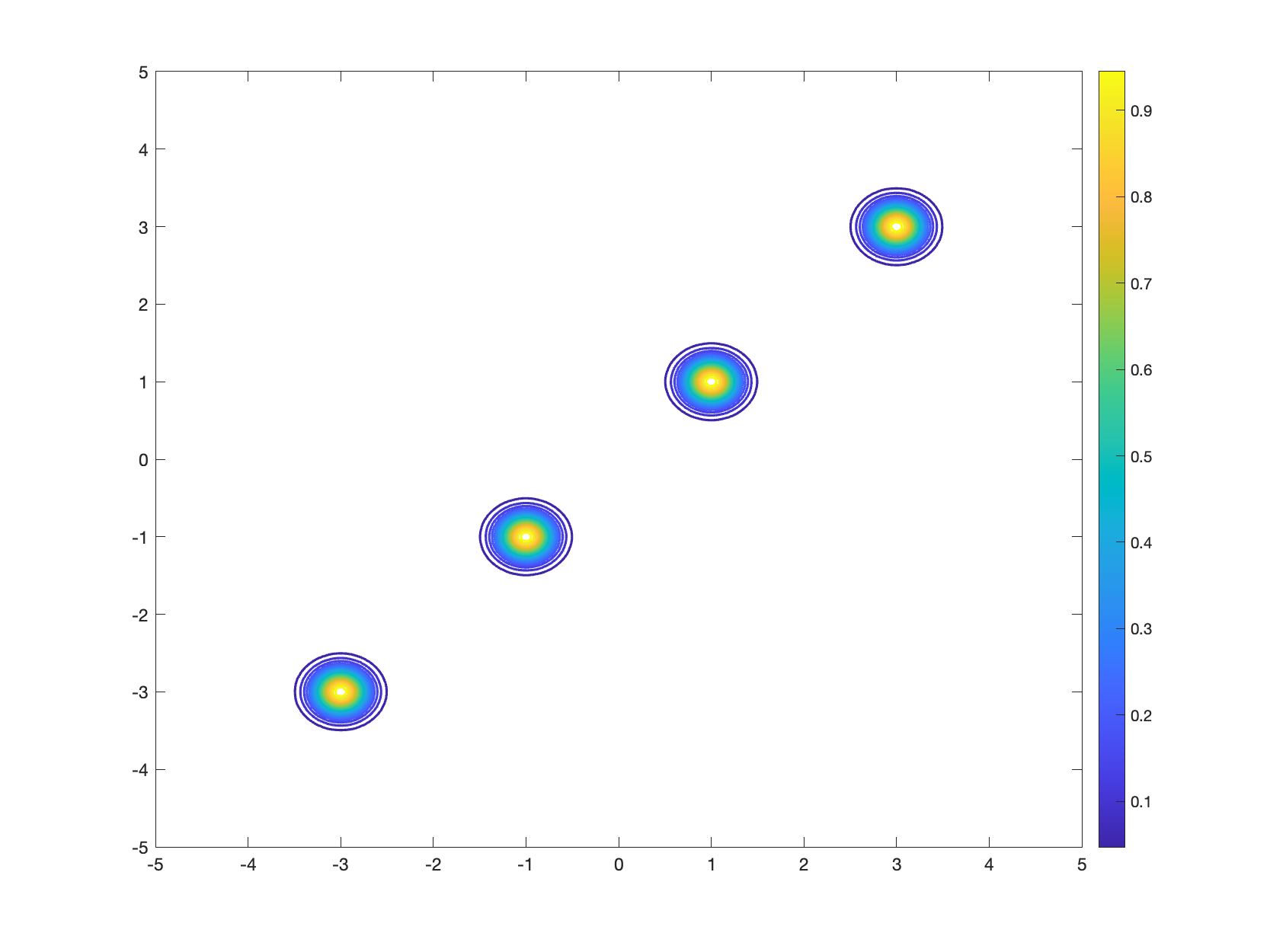} 
        \caption{True contour plot}	
    \end{subfigure}
    \hspace{-0.5cm} 
    \begin{subfigure}{0.32\linewidth}
        \includegraphics[width=\linewidth]{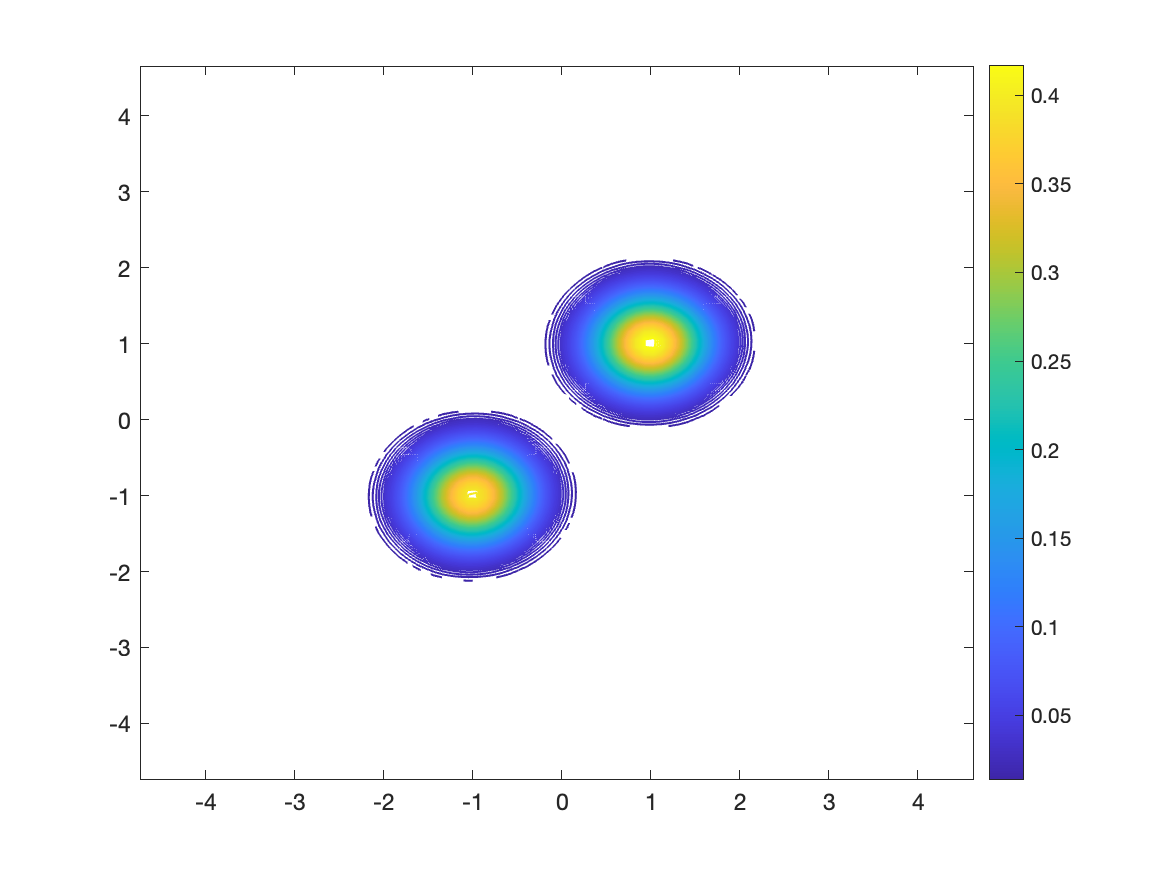}
        \caption{$\beta=1$}	
    \end{subfigure}
    \hspace{-0.5cm}
    \begin{subfigure}{0.32\linewidth}
        \includegraphics[width=\linewidth]{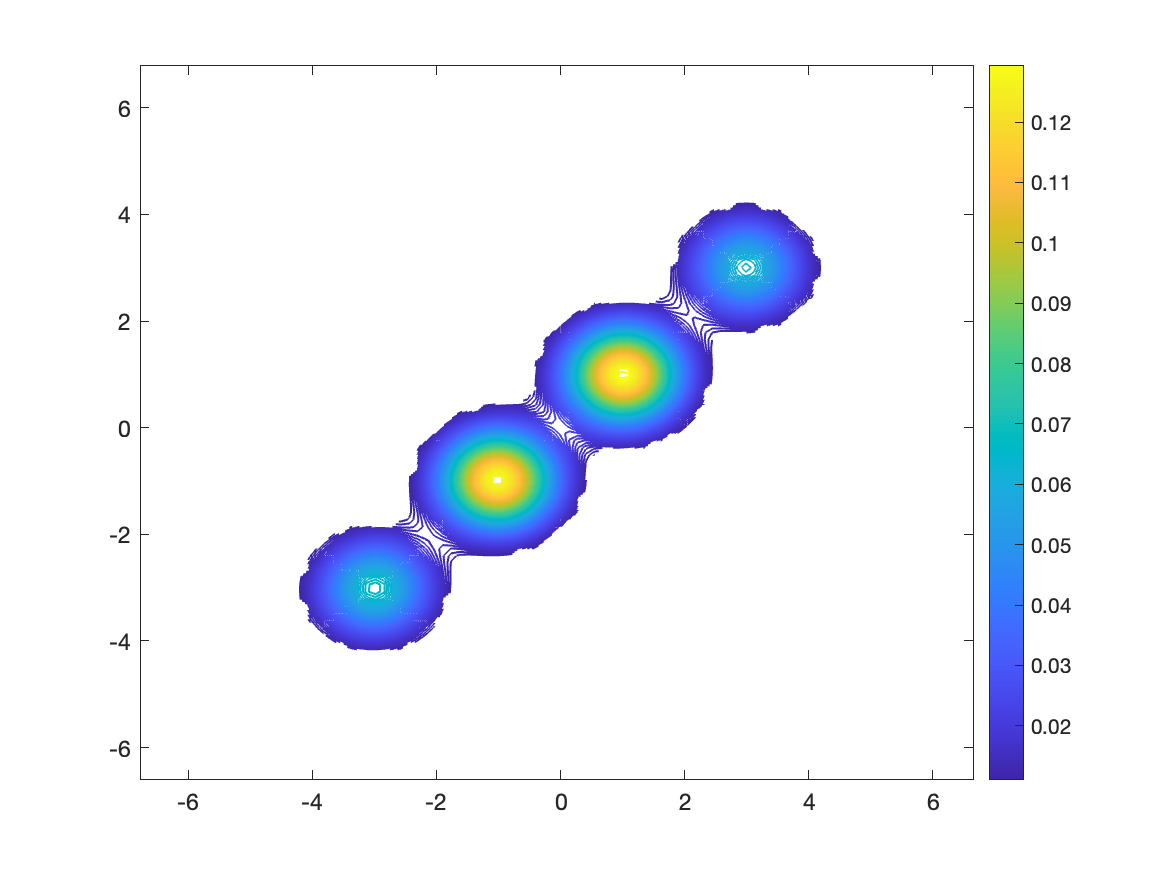}
        \caption{$\beta=5$}	
    \end{subfigure}
    \caption{Sampling via SFS with different temperatures $\beta = 1,  5$.}
    \label{figure:mog_dens_kappa_4_different_d}
\end{figure}

\textbf{Sampling via SFS with different temperatures $\beta$:}
By fixing the dimension $d=100$ and a step-size $h=10^{-3}$, we simulate $1000$ independent Markov chains.
In Figure \ref{figure:mog_dens_kappa_2_different_d}, we set $\kappa=2,\,\theta_1=\tfrac{1}{4},\, \theta_2=\tfrac{3}{4},\,\alpha_1=-2\cdot \mathbf{1}_d,\alpha_2=4\cdot \mathbf{1}_d$ and  $\Sigma_1=\Sigma_2=\tfrac{1}{4}\cdot\mathbf{I}_d$.
In Figure \ref{figure:mog_dens_kappa_4_different_d}, we set $\kappa=4,\,\alpha_1=-3\cdot \mathbf{1}_d,\alpha_2=-\mathbf{1}_d,\,\alpha_3=\mathbf{1}_d,\alpha_4=3\cdot \mathbf{1}_d$ and  $\theta_i=\tfrac14,\,\Sigma_i=\tfrac{1}{5}\cdot\mathbf{I}_d, \, i\in [4]$.
Figure \ref{figure:mog_dens_kappa_4_different_d} (a) displays the true contour plots.
For the case $\beta =1$, the SFS algorithm \eqref{eq:Euler-scheme-no-mc} reduces into the SFS algorithm with exact drift in \cite{MRjiao}. As shown in Figures \ref{figure:mog_dens_kappa_2_different_d} (a) and \ref{figure:mog_dens_kappa_4_different_d} (b), the generated samples may collapse onto only one or two modes.
As the temperature $\beta$ increases to  $\beta =5$, the SFS algorithms show significant improvements in sampling performance and generate high-quality samples with the desired numbers of modes, as observed in Figure \ref{figure:mog_dens_kappa_2_different_d} (b) and Figure \ref{figure:mog_dens_kappa_4_different_d} (c).

\subsubsection{Convergence rate of $h$ for Gaussian mixtures distributions}
%\textbf{Convergence rate:}
To test the convergence rates of the SFS algorithm \eqref{eq:Euler-scheme-no-mc}, we fix $\kappa=2,\,\theta_1=\theta_2=\tfrac{1}{2},\,\alpha_1=-2\cdot\mathbf{1}_d, \, \alpha_2=2\cdot\mathbf{1}_d,\, \Sigma_1=\Sigma_2=\tfrac{1}{4}\cdot\mathbf{I}_d$ and run the scheme using different step-sizes $h\in\{2^{-5},2^{-6},2^{-7},2^{-8},2^{-9}\}$.
The exact solution is identified as the numerical one using a fine step-size $h_{ref}=2^{-14}$ and the expectations are approximated by computing averages over $5000$ samples. \\

\begin{figure}[htbp]
    \centering
    \begin{subfigure}{0.4\linewidth}  
        \includegraphics[width=\linewidth]{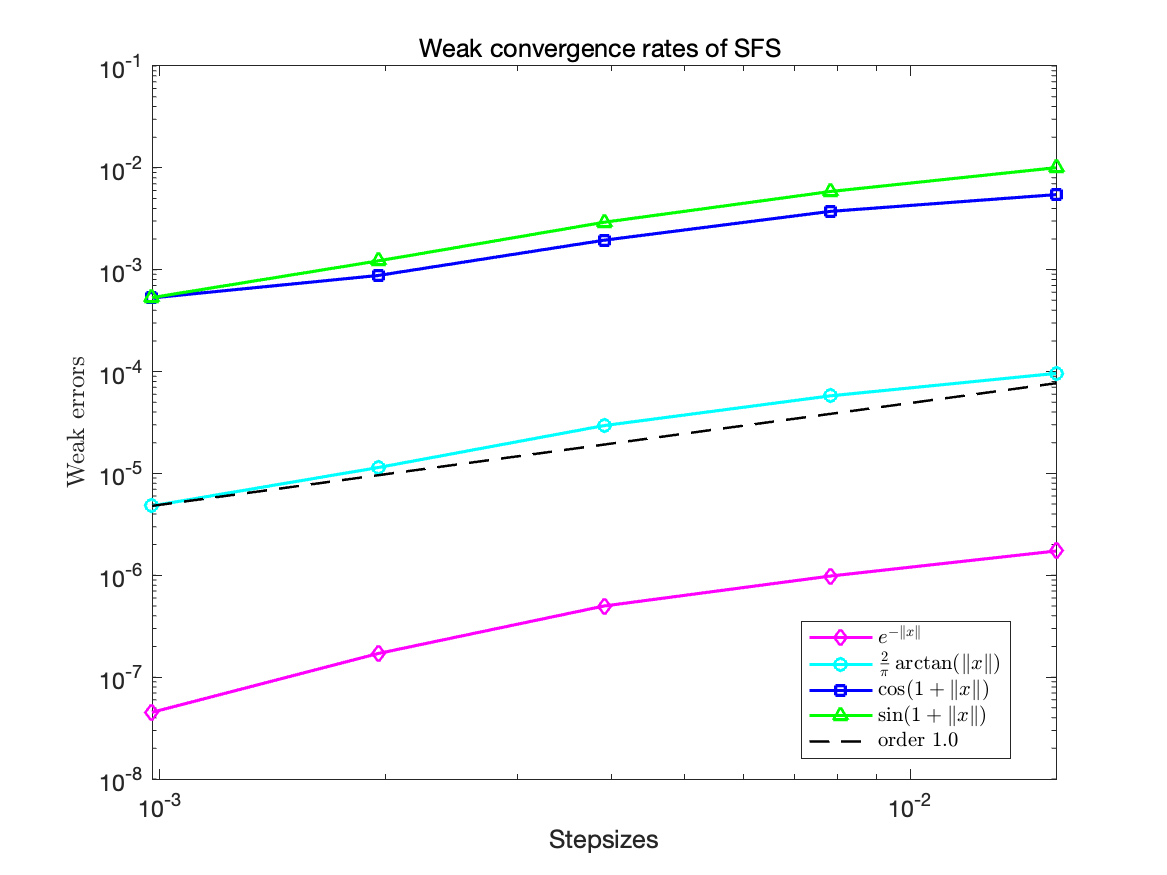} 
        \caption{$\beta=1,\,d=20$}	
    \end{subfigure}
    \hspace{-0.5cm} 
    \begin{subfigure}{0.4\linewidth}
        \includegraphics[width=\linewidth]{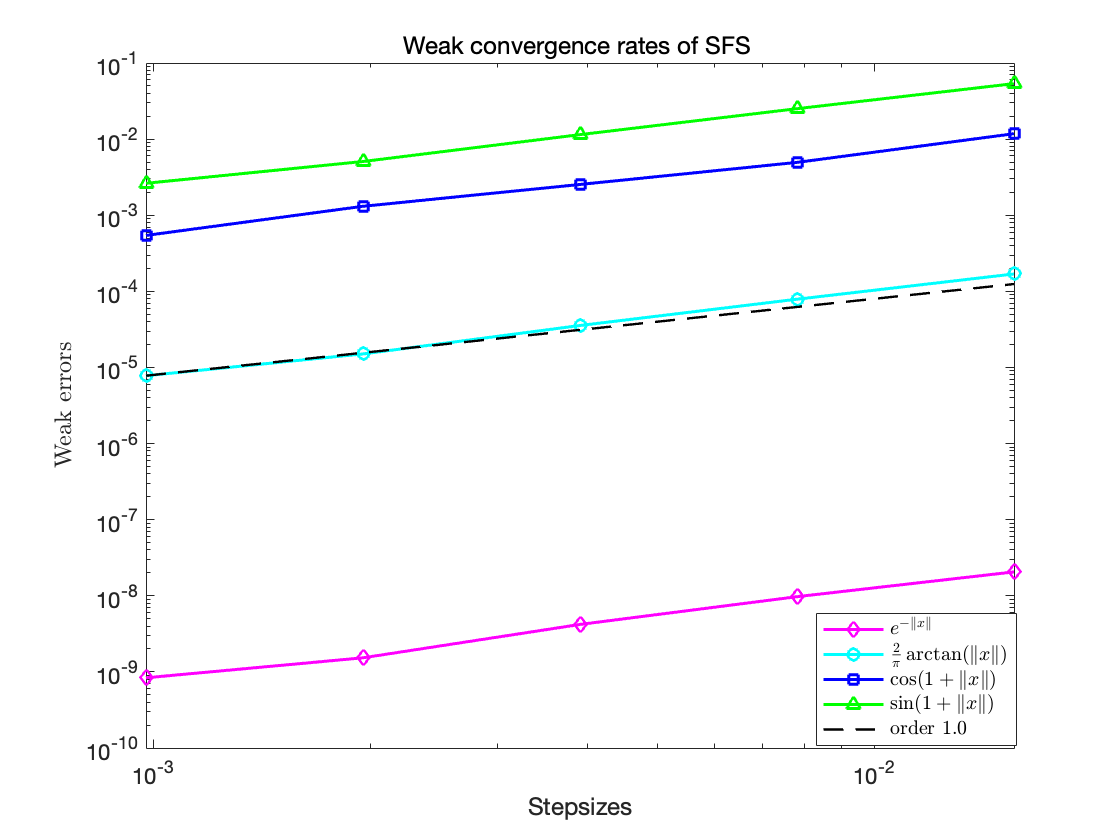}
        \caption{$\beta=5,\,d=50$}	
    \end{subfigure}
    \caption{Weak convergence rates of SFS of the  Gaussian mixture distributions.
    %for $d = 50$ (left) and $d = 80$ (right)
    }
    \label{figure:mog_conv}
\end{figure}

\noindent
\textbf{Test functions:} 
Throughout this subsection, we consider the following four test functions:
\begin{equation}
\label{eq:test-function}
    \phi(x) 
    \in
    \big\{ e^{-\|x\|},
    \tfrac{2}{\pi}\arctan (\|x\|),
    \cos(1+\|x\|),
    \sin(1+\|x\|)
    \big\}
    \quad \forall x \in \mathbb{R}^d.
\end{equation}
Note that all test functions $\phi \in B_b(\mathbb{R}^d)$ with $|\phi|_0 \leq  1$.
%(ii) For $\mathcal{W}_1$-distance, we consider the following two functions
%\begin{equation}
%    \phi_2(x) 
%    \in
%    \left\{\cos(\|x\|),
%    \sin(1+\|x\|)
%    \right\}, \quad \forall x \in \mathbb{R}^d.
%\end{equation}
%Note that all test functions $\phi_2 \in \mathrm{\rm Lip}(\mathbb{R}^d)$ with $|\phi|_{\mathrm{Lip}}=1$.
%
The reference lines of slope $1$ are also presented, depicted as black solid lines in Figure \ref{figure:mog_conv}. 
Approximate convergence rates for each test function are also provided in Figure \ref{figure:mog_conv}.

\subsubsection{Dependence of the errors on the Monte Carlo sample size $M$}
To study the  Monte Carlo sample size $M$ dependence of sampling errors, we let $\kappa=3,\,\beta=1, \,\theta_1=\tfrac{1}{5},\,\theta_2=\theta_3=\tfrac{2}{5}, \alpha_1=\mathbf{0}_d,\,\alpha_2=\mathbf{1}_d, \, \alpha_3=-\mathbf{1}_d,\, \Sigma_1=\mathbf{I}_d, \, \Sigma_2=\Sigma_3=\tfrac{4}{5}\cdot\mathbf{I}_d$ 
and run the scheme \eqref{eq:e-m-mc} using different Monte Carlo sample size $M \in \{4,8,16,32,64,128\}$.
By fixing the dimension $d=10$ and a step-size $h=2^{-6}$, we simulate $5000$ independent Markov chains.

\begin{figure}[htbp]
    \centering
    \begin{subfigure}{0.4\linewidth}  
        \includegraphics[width=\linewidth]{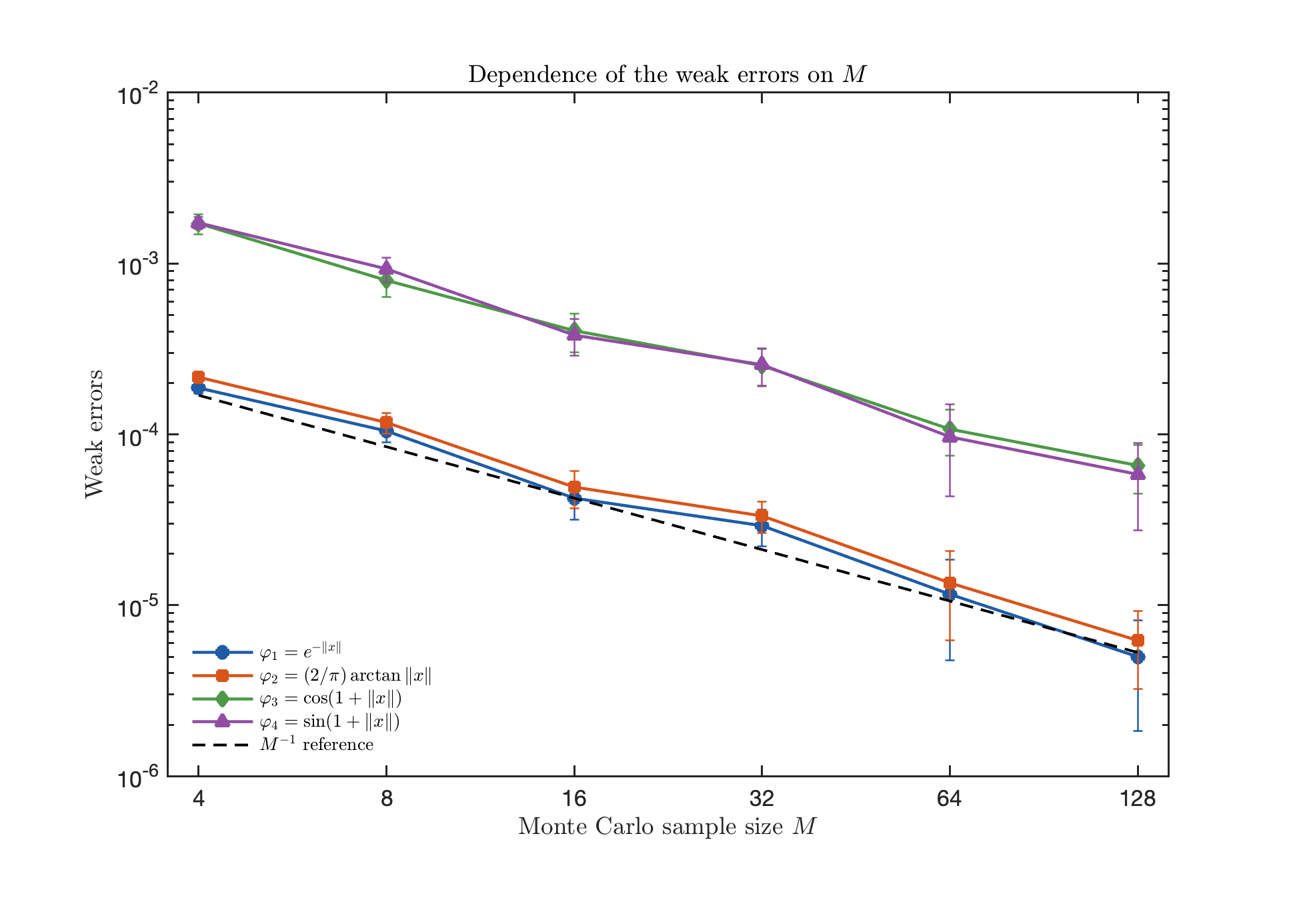} 
        \caption{Weak error convergence}	
    \end{subfigure}
    \hspace{-0.5cm} 
    \begin{subfigure}{0.4\linewidth}
        \includegraphics[width=\linewidth]{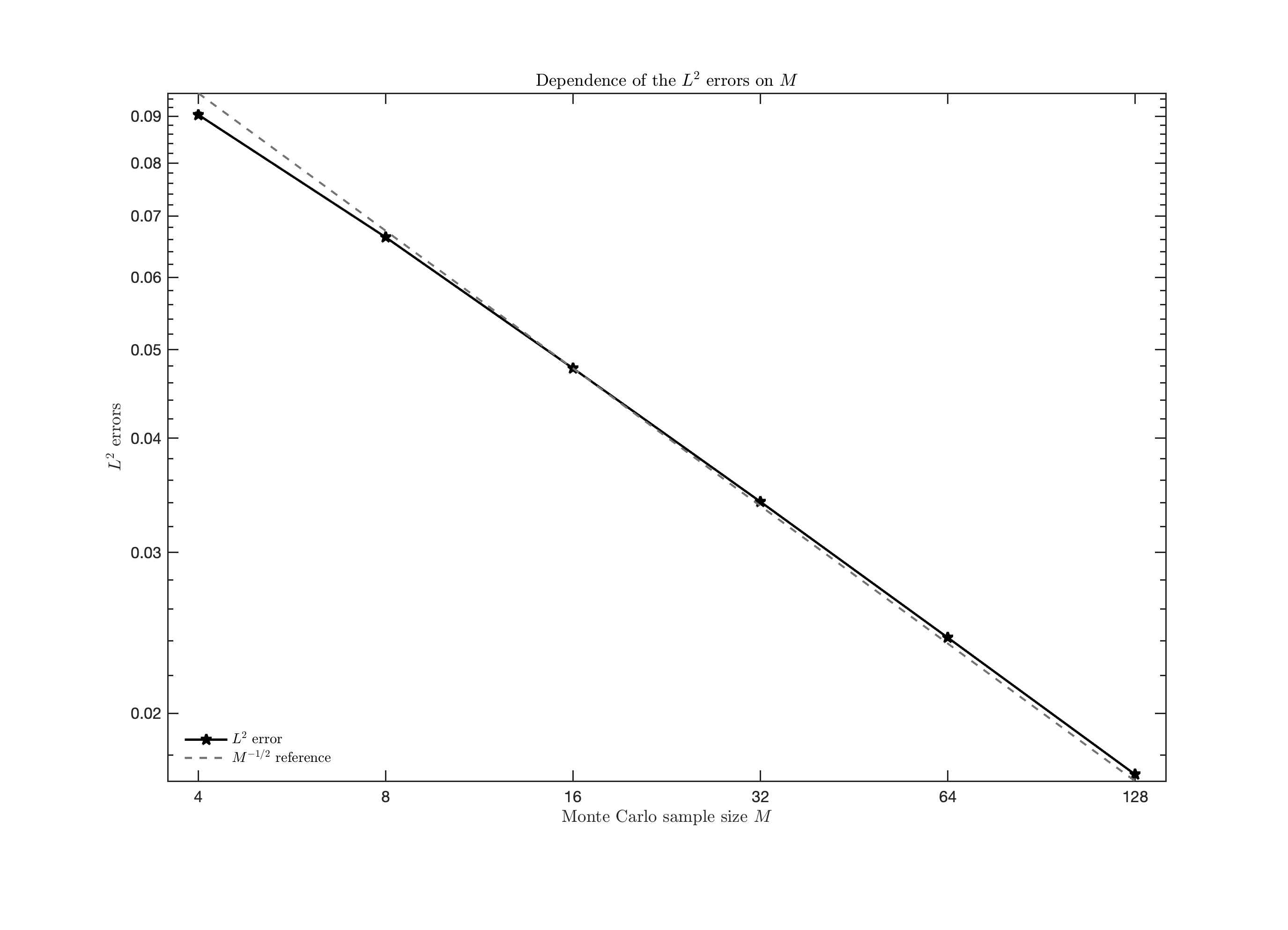}
        \caption{$L^2$-error convergence}	
    \end{subfigure}
    \caption{Monte Carlo convergence rates of SFS for a Gaussian mixture distributions.
    }
    \label{figure:dependence_M}
\end{figure}

We assess the total variation distance through weak errors against the bounded test functions in \eqref{eq:test-function}, while using the strong
$L^2$ error as a computable surrogate for
the $L^2$-Wasserstein distance.
As shown in Figure \ref{figure:dependence_M}(a), the weak errors decay approximately at the rate $\mathcal{O}(M^{-1})$, confirming the previous findings concerning the error estimates.
From Figure \ref{figure:dependence_M}(b), one can observe a convergence rate with respect to $M$ close to $\mathcal{O}(M^{-1/2})$ in $L^2$-error.
For comparison, black dashed reference lines corresponding to the $\mathcal{O} (M^{-1})$ and $\mathcal{O} (M^{-1/2})$ in
panels (a) and (b), respectively.
Together, the two panels illustrate that the Monte Carlo error scales as $\mathcal{O} (M^{-1})$ in the total variation distance, significantly improving the conventional rate of order $\mathcal{O} (M^{-1/2})$ in the $L^2$-Wasserstein distance.

\section{Conclusion and future work}
\noindent
In this work, we establish a convergence rate of order $\mathcal{O}(d h |\ln h|)$ for the time discretization of the Schr\"odinger-F\"ollmer sampler (SFS) with exact drift in the total variation distance under suitable conditions. 
When the drift is approximated by Monte Carlo methods, we obtain a convergence rate of order $\mathcal{O}(dh|\ln h|+d/M)$ in the total variation distance.
The Monte Carlo contribution thus has a sharper dependence on $M$ than the existing $\mathcal{O}(1/\sqrt{M})$ bounds in the $L^2$-Wasserstein distance established in \cite{MRjiao,Jiao2026Convergence,Wang2026Multimodal}.
%To the best of our knowledge, this is the first result to consider the convergence rate of SFS algorithm in the total variation distance, achieving an almost first-order rate with respect to the time step size.
%
Numerical experiments demonstrate the effectiveness of the proposed algorithm and the convergence rate with respect to the time step-size $h$ and the Monte Carlo sample size $M.$

Several interesting lines of exploration remain for future work
\cite{Wang2025randomization}.
%
%First, we propose a higher-order (Runge-Kutta) sampling algorithms based on such Schr\"odinger-F\"ollmer diffusion process with temperatures. 
%
%Second, all analyses for Schr\"odinger-F\"ollmer sampler work in $L^2$-Wasserstein distance. 
For example, it would be interesting to propose sampling algorithms based on the Schr\"odinger-F\"ollmer diffusion process with temperatures that achieve an optimal convergence rate of $\mathcal{O}(\sqrt {d}h)$, or even higher-order sampling algorithms, under only a Lipschitz condition on the drift.

%This may be challenging, as it seems different from our proof technique.

\bibliographystyle{abbrv}
\bibliography{reference}

\appendix

\section{Proof of Lemma \ref{lem:X-mean-square-differential}}
\label{appendix:D^3-x-proof}
\textbf{Proof.} 
Without loss of generality, we denote that, for any $ x, v_1, v_2, v_3\in \mathbb{R}^d$,
$$
    \eta^{v_1}(s,x;t)
    :=
    \mathcal{D} X_t^{s,x} v_1, 
    \quad 
    \xi^{v_1, v_2}(s,x;t)
    :=
    \mathcal{D}^2 X_t^{s,x}(v_1, v_2), 
$$
and 
$$
    \zeta^{v_1, v_2,v_3}(s,x;t)
    :=
    \mathcal{D}^3 X_t^{s,x}(v_1, v_2,v_3).
$$
Following a similar argument used in \cite[Lemma 4.3]{PWW25}, one can derive the estimates for $\mathcal{D}X^{s,x}_t$ and $\mathcal{D}^2X^{s,x}_t$. Hence, we omit their proof here.
As a result, in the same manner, we have the third variation process with respect to
SDE \eqref{eq:sde} as follows:
\begin{equation}
    \begin{aligned}
         \mathrm{\, d} \zeta^{v_1,v_2,v_3}
         (s,x; t) 
         & =
         \Big(
         D f_{\beta}(X^{s,x}_t,t) 
         \zeta^{v_1, v_2, v_3}
         (s,x; t)
         +
         D^2 f_{\beta}(X^{s,x}_t,t)
        \big(
        \xi^{v_1, v_2}(s,x; t), 
        \eta^{v_3}(s,x; t)
        \big)\\
        & \quad +
        D^2 f_{\beta}(X^{s,x}_t,t)
        \big(
        \xi^{v_1, v_3}
        (s,x; t),
        \eta^{v_2}(s,x; t)
        \big)\\
        & \quad +
        D^2 f_{\beta}(X^{s,x}_t,t)
        \big(
        \eta^{v_1}(s,x; t),
        \xi^{v_2 ,v_3}(s,x; t)
        \big)
        \Big)
        \mathrm{\, d} t\\
        & \quad +
        D^3 f_{\beta}(X^{s,x}_t,t)
        \big(
        \eta^{v_1}(s,x; t), 
        \eta^{v_2}(s,x; t),
        \eta^{v_3}(s,x; t)
        \big)\\
         & :=
         \big(
         D f_{\beta}(X^{s,x}_t,t) 
         \zeta^{v_1, v_2, v_3}
         (s,x; t)
         +
         H^{v_1, v_2, v_3}(X^{s,x}_t)
         \big)
         \mathrm{\, d} t, 
    \end{aligned}
\end{equation}
where $\zeta^{v_1, v_2, v_3}(s,x;s)=0.$
It follows from \eqref{eq:one-sided-Lipschitz} that, 
\begin{equation}
    \big \langle 
    y,
    Df_{\beta}(x,t) y
    \big \rangle 
    \leq 
    L_f\|y\|^2,
    \quad \forall x,y \in \mathbb{R}^d, 
    \quad t \in [0,1].
\end{equation}
Following the same argument as \cite[Lemma 4.3]{PWW25} and the Young inequality, we deduce that
\begin{equation}
\label{eq:zeta-expression}
    \begin{aligned}
        &
        \mathrm{\, d}
        \|
        \zeta^{v_1, v_2, v_3}
        (s,x; t)
        \|^2 \\
        & =
        \Big(
        2
        \big\langle
        \zeta^{v_1, v_2, v_3}
        (s,x; t), 
        D f_{\beta}(X_t^{s,x},t) 
        \zeta^{v_1, v_2, v_3}
        (s,x; t)
        \big\rangle
        +2
        \big\langle
        \zeta^{v_1, v_2, v_3}
        (s,x; t),
        H^{v_1, v_2, v_3}(X_t^{s,x})
        \big \rangle
        \Big)
        \mathrm{\, d} t\\
        & \leq
        (2 L_f+1)
        \big\|
        \zeta^{v_1, v_2, v_3}
        (s,x; t)
        \big\|^2 
        \mathrm{\, d} t
        +
        \big\|
        H^{v_1, v_2, v_3}(X_t^{s,x})
        \big\|^2 
        \mathrm{\, d} t .
    \end{aligned}
\end{equation}
We begin with the estimation of $\|H^{v_1, v_2, v_3}(X_t^{s,x})\|_{L^2(\Omega,\mathbb{R}^d)}$ as follows:
\begin{align*}
  \big\|
  H^{v_1, v_2, v_3}(X_t^{s,x})
  \big\|
   _{L^2(\Omega,\mathbb{R}^d)} 
    & \leq
    \big \|D^2 f_{\beta}(X_t^{s,x},t)
    \big(
    \xi^{v_1, v_2}(s,x; t), 
    \eta^{v_3}(s,x; t)
    \big\|_{L^2(\Omega,\mathbb{R}^d)}
    \nonumber\\
    & \quad +
    \big\|
    D^2 f_{\beta}(X_t^{s,x},t)
    \big(
    \xi^{v_1, v_3}(s,x; t),
    \eta^{v_2}(s,x; t)
    \big)
    \big\|_{L^2(\Omega,\mathbb{R}^d)}
    \nonumber\\
    & \quad +
    \big\|
    D^2 f_{\beta}(X_t^{s,x},t)
    \big(
    \eta^{v_1}(s,x; t),
    \xi^{v_2 ,v_3}(s,x; t)
    \big)
    \big\|_{L^2(\Omega,\mathbb{R}^d)}\\
    & \quad +
    \big\|
    D^3 f_{\beta}(X_t^{s,x},t)
    \big(
    \eta^{v_1}(s,x; t), 
    \eta^{v_2}(s,x; t),
    \eta^{v_3}(s,x; t)
    \big)
    \big\|_{L^2(\Omega,\mathbb{R}^d)}.
\end{align*}
The key issue is to estimate 
$
\| D^2 
f_{\beta}(X_t^{s,x},t)
(\cdot,\cdot)
\|_{L^2(\Omega,\mathbb{R}^d)}
$
and 
$
\| D^3 
f_{\beta}(X_t^{s,x},t)
(\cdot,\cdot,\cdot)
\|_{L^2(\Omega,\mathbb{R}^d)}
$.
Using Proposition~\ref{prop:f-in-C^3} together with
\eqref{eq:DX-estimation}, we obtain the following estimates:
\begin{equation}
    \begin{aligned}
        \big\| 
        D^2 f_{\beta}(X_t^{s,x},t)
        \big(
        \xi^{v_1, v_2}(s,x; t),
         \eta^{v_3}(s,x; t)
         \big)
        \big\|
        _{L^2(\Omega,\mathbb{R}^d)}
        & \leq
        C
        \big\|
        \|\xi^{v_1, v_2}(s,x; t)\|
        \|\eta^{v_3}(s,x; t)\|
        \big\|_{L^2(\Omega,\mathbb{R}^d)}
        \\
%        & \leq
%         C 
%         \|\xi^{v_1, v_2}(t,x)\|
%         _{L^4(\Omega,\mathbb{R}^d)}
%         \cdot
%         \|\eta^{v_3}(t,x)\| \\
        & \leq
        C e^{(4L_f+1/2)(t-s)}
        ~\|v_1\|
        \cdot
        \|v_2\|
        \cdot
        \|v_3\|,
    \end{aligned}
\end{equation}
and 
\begin{equation}
    \begin{aligned}
        &
        \big\| D^3 f_{\beta}(X_t^{s,x},t)
        (\eta^{v_1}(s,x; t),
        \eta^{v_2}(s,x; t),
        \eta^{v_3}(s,x; t))
        \big\| 
        _{L^2(\Omega,\mathbb{R}^d)}\\
        & \leq
        C
        \|\eta^{v_1}(s,x; t)\|
        \cdot
        \|\eta^{v_2}(s,x; t)\|
        \cdot
        \|\eta^{v_3}(s,x; t)\|\\
        & \leq 
        C e^{3L_f (t-s)}
        ~\|v_1\|
        \cdot
        \|v_2\|
        \cdot
        \|v_3\|.
    \end{aligned}
\end{equation}
Taking expectations on both sides of \eqref{eq:zeta-expression}, we can easily get, for any $x, v_1, v_2, v_3 \in \mathbb{R}^d$,
\begin{equation}
    \begin{aligned}  
        \mathbb{E}
        \Big[
        \big\|
        \zeta^{v_1, v_2, v_3}
        (s,x; t)
        \big\|^2
        \Big] 
        \leq & 
        (2 L_f+1) 
        \int_s^t 
        \mathbb{E}
        \Big[
        \big\|
        \zeta^{v_1, v_2, v_3}
        (s,x; r)
        \big\|^2
        \Big] 
        \mathrm{\, d} r \\
        & +
        C e^{(8 L_f +1)(t-s)}
        \left\|v_1\right\|^2 
        \cdot
        \left\|v_2\right\|^2
        \cdot
        \left\|v_3\right\|^2, 
        %\quad \forall x, v_1, v_2, v_3 \in \mathbb{R}^d,
    \end{aligned}
\end{equation}
leading to 
\begin{equation}
    \left\|
    \zeta^{v_1, v_2, v_3}(s,x; t)
    \right\|_{L^2(\Omega, \mathbb{R}^d)} 
    \leq 
    C e^{(5 L_f+1) (t-s)}
    \|v_1\| 
    \cdot
    \|v_2\|
    \cdot
    \|v_3\|,
    %\quad \forall x, v_1, v_2, v_3 \in \mathbb{R}^d,
\end{equation}
with the help of the Gronwall inequality. 
The proof is thus completed.

\section{Proof of Lemma \ref{lem:regularization-estimates-u}}
\label{appendix:proof-lemms-estimates-u}
\begin{proof}
\textbf{Proof.}
The Bismut-Elworthy-Li formula (see Lemma \ref{lem:Bismut-Elworthy-Li-formula}) states that, for some function $\phi: \mathbb{R}^d \to \mathbb{R}$ belongs to $C_b(\mathbb{R}^d)$, 
%there exists $\mathcal{K}(\phi, x)>0$, which depends on $\phi$ and $x$, such that,
%\begin{align}
%\label{eq:phi-x}
%        |\phi(x)| 
%        \leq 
%        \mathcal{K}(\phi, x), \quad \forall x \in \mathbb{R}^d,
%\end{align}
one can calculate the first, the second and the third derivative of
$$
u(t, x)
:=
\mathbb{E}
\big[
\phi(X_1^{1-t,x})
\big]
$$
with respect to $x$. According to \eqref{eq:bel}, we have 
\begin{equation}
\label{eq:Du-exp}
D u(t, x)v_1
=
\frac{1}{\sqrt{\beta}\,t}
\mathbb E
\bigg[
\int_{1-t}^{1}
\left\langle
\eta^{v_1}(1-t,x;r),
\mathrm d W_r
\right\rangle
\phi(X^{1-t,x}_1)
\bigg].
\end{equation}
The estimate for the above equation follows from an argument analogous to that in \cite[Theorem 4.5]{PWW25}, by applying the H\"older's inequality and the It\^o isometry together with Lemma~\ref{lem:X-mean-square-differential}. As a result, we omit the proof here. Now, let us concentrate on the estimation of the second derivative \eqref{eq:D2-u-est}. The Markov property of SDE \eqref{eq:sde} \cite[Corollary 9.2]{da2014introduction} yields
\begin{equation}
u(t,x)
=
P_{1-t,1}
\phi(x)
=
P_{1-t,1-t/2}
P_{1-t/2,1}
\phi(x).
\end{equation}
Then by \eqref{eq:Du-exp}, one gets
\begin{equation}
\begin{aligned}
D u(t, x)v_1
& =
\frac{2}{\sqrt{\beta}\,t}
\mathbb E
\bigg[
\int_{1-t}^{1-t/2}
\left\langle
\eta^{v_1}(1-t,x;r),
\mathrm d W_r
\right\rangle
u(
t/2,
X_{1-t/2}^{1-t,x}
)
\bigg].
\end{aligned}
\end{equation}
Accordingly, one can compute the second derivative of $u(t,\cdot)$ using the chain rule as follows:
\begin{equation}
\begin{aligned}
D^2u(t,x)(v_1,v_2)
=
&\,
\frac{2}{\sqrt{\beta}\,t}
\mathbb E
\bigg[
\int_{1-t}^{1-t/2}
\left\langle
\xi^{v_1,v_2}(1-t,x;r),
\mathrm dW_r
\right\rangle
u(t/2,
X_{1-t/2}^{1-t,x})
\bigg]
\\
&+
\frac{2}{\sqrt{\beta}\,t}
\mathbb E
\bigg[
\int_{1-t}^{1-t/2}
\left\langle
\eta^{v_1}(1-t,x;r),
\mathrm dW_r
\right\rangle
Du(t/2,
X_{1-t/2}^{1-t,x})
\eta^{v_2}(1-t,x;1-t/2)
\bigg].
\end{aligned}
\end{equation}
Due to \eqref{eq:D-u-est}, Lemma \ref{lem:X-mean-square-differential}, the H\"older inequality and the It\^o isometry \cite[Lemma 5.4]{mao2007stochastic}, one has
\begin{align}
&
\big|
D^2u(t,x)(v_1,v_2)
\big|
\nonumber\\
& \leq 
\frac{2}{\sqrt{\beta}\,t}
\bigg\|
\int_{1-t}^{1-t/2}
        \big\langle
        \xi^{v_1,v_2}
        (1-t,x;r), 
        \mathrm{\, d} W_r
        \big\rangle
        \bigg\|
        _{L^2(\Omega, \mathbb{R})}
        \Big\|
        u\big(t / 2,
        X_{1-t/2}^{1-t,x}
        \big)
        \Big\|
        _{L^2(\Omega, \mathbb{R})}
        \nonumber\\
& \quad +
\frac{2}
{\sqrt{\beta}\,t}
\bigg\|
\int_{1-t}^{1-t/2}
\left\langle
\eta^{v_1}(1-t,x;r),
\mathrm dW_r
\right\rangle
\bigg\|_{L^2(\Omega, \mathbb{R})}
\Big\|
Du(t/2,
X_{1-t/2}^{1-t,x})
\eta^{v_2}(1-t,x;1-t/2)
\Big\|_{L^2(\Omega, \mathbb{R})} \nonumber\\
& \leq 
C 
\frac{
|\phi|_0}
{t}
\bigg(
\int_{1-t}^{1-t/2}
\big\|
\xi^{v_1,v_2}(1-t,x;r)
\big\|^2_{L^2(\Omega, \mathbb{R}^d)} 
\mathrm{\, d} r
\bigg)^{1 / 2}
\nonumber \\
& \quad +
C 
\frac{
|\phi|_0}
{t^{3/2}}
\bigg(
\int_{1-t}^{1-t/2}
\big\|
\eta^{v_1}(1-t,x;r)
\big\|^2_{L^2(\Omega, \mathbb{R}^d)} 
\mathrm{\, d} r
\bigg)^{1 / 2}
\big\|
\eta^{v_2}(1-t,x;1-t/2)
\big\|
\nonumber \\
& \leq
        C 
        \frac{
        |\phi|_0}
        {t}
        \big(
        e^{(3 L_f+1/2)t}-1
        \big)^{1 / 2}
        \|v_1\| 
        \cdot\|v_2\|
         +
        C
        \frac{
        |\phi|_0}
        { t^{3/2} }
        e^{L_f t/2} 
        \big(
        e^{L_f t} -1
        \big)^{1/2}
        \|v_1\|\cdot
        \|v_2\|
        \notag \\
        & \leq  
        \frac{C |\phi|_0 }
        {t}
        \|v_1\| 
        \cdot\|v_2\|,
\end{align}
where we have also used the fact that $e^{C t}-1 \leq (e^C-1)t$ for any $t\in[0,1], \, C>0$ and $1/\sqrt{t} \leq 1/t, \, t \in (0,1]$.
By proceeding in the same way as \cite{PWW25, cerrai2001second}, it is possible to prove that the following formula holds for the third derivative
{\small
\begin{align}
\label{eq:D^3-u-expression}
        &
        D^3 u(t, x)
        (v_1, v_2,v_3)
        \nonumber\\
        & =  
        \frac{2}{\sqrt{\beta}\,t} 
        \mathbb{E}
        \bigg[
        \int_{1-t}^{1-t/2}
        \big\langle
        \zeta^{v_1, v_2,v_3}
        (1-t,x ;r), 
        \mathrm{\, d} W_r
        \big\rangle 
        u
        (t / 2, X_{1-t / 2}^{1-t,x})
        \bigg] \nonumber\\
        & \quad +
        \frac{2}{\sqrt{\beta}\,t} 
        \mathbb{E}
        \bigg[
        \int_{1-t}^{1-t/2}
        \big\langle
        \xi^{v_1, v_2}(1-t,x ;r), 
        \mathrm{\, d} W_r
        \big\rangle 
        Du
        (t / 2, X_{1-t / 2}^{1-t,x})
        \eta^{v_3}(1-t,x;1-t/2)
        \bigg]\nonumber \\
        & \quad +
        \frac{2}{ \sqrt{\beta}\,t} 
        \mathbb{E}
        \bigg[
        \int_{1-t}^{1-t/2}
        \big\langle
        \xi^{v_1, v_3}(1-t,x ;r), 
        \mathrm{\, d} W_r
        \big\rangle 
        D u
        (t / 2, X_{1-t / 2}^{1-t,x})
        \eta^{v_2}(1-t,x;1-t/2)
        \bigg] \nonumber\\
        & \quad +
        \frac{2}{\sqrt{\beta} \, t} 
        \mathbb{E}
        \bigg[
        \int_{1-t}^{1-t/2}
        \big\langle
        \eta^{v_1}(1-t,x ;r), 
        \mathrm{\, d} W_r
        \big\rangle 
        D^2 u
        (t / 2, X_{1-t / 2}^{1-t,x})
        \big(
        \eta^{v_3}(1-t,x ;1-t/2),
        \eta^{v_2}(1-t,x ;1-t/2)
        \big)
        \bigg] \nonumber\\
        & \quad +
        \frac{2}{\sqrt{\beta} \, t} 
        \mathbb{E}
        \bigg[
        \int_{1-t}^{1-t/2}
        \big\langle
        \eta^{v_1}(1-t,x ;r), 
        \mathrm{\, d} W_r
        \big\rangle 
        D u
        (t / 2, X_{1-t / 2}^{1-t,x})
        \,
        \xi^{v_2, v_3}(1-t,x ;1-t/2)
        \bigg]
        \nonumber\\
        & =:
        I_1+I_2+I_3+
        I_4+I_5.
\end{align}}
Using Lemma~\ref{lem:X-mean-square-differential}, the H\"older inequality,
and the It\^o isometry \cite[Lemma 5.4]{mao2007stochastic}, together with
the fact that
$e^{Ct}-1\leq (e^C-1)t,\, t\in[0,1],\ C>0,$
one can easily get
\begin{align}
\label{eq:I1-est}
        |I_1|
        & \leq
        \frac{2}{\sqrt{\beta}\,t}
        \bigg\|
        \int_{1-t}^{1-t/2}
        \big\langle
        \zeta^{v_1, v_2, v_3}
        (1-t,x ;r), 
        \mathrm{\, d} W_r
        \big\rangle
        \bigg\|
        _{L^2(\Omega, \mathbb{R})}
        \Big\|
        u
        \big(t / 2, 
        X_{1-t / 2}^{1-t,x}\big)
        \Big\|
        _{L^2(\Omega, \mathbb{R})} 
        \nonumber\\
        & \leq 
        C 
        \frac{
        |\phi|_0}
        {t}
        \bigg(
        \int_{1-t}^{1-t/2}
        \big\|
        \zeta^{v_1, v_2, v_3}(1-t,x ;r)
        \big\|^2_{L^2(\Omega, \mathbb{R}^d)} 
        \mathrm{\, d} r
        \bigg)^{1 / 2}
        \nonumber \\
        & \leq
        C 
        \frac{
        |\phi|_0}
        {t}
        \big(
        e^{(5 L_f+1)t}-1
        \big)^{1 / 2}
        \|v_1\| 
        \cdot\|v_2\|
        \cdot\|v_3\|
        \notag \\
        & \leq  
        \frac{C |\phi|_0}
        {\sqrt{t}}
        \|v_1\| 
        \cdot\|v_2\|
        \cdot\|v_3\|,
\end{align}
In the following, we notice that the analysis of $I_2, I_3$ and $I_5$ is similar. Taking the $I_2$ as an example, by Lemma \ref{lem:X-mean-square-differential}, the H\"older inequality and the It\^o isometry \cite[Lemma 5.4]{mao2007stochastic}, we can easily get 
\begin{align}
\label{eq:I2-est}
        |I_2|
        & \leq
        \frac{2}{\sqrt{\beta}\,t}
        \bigg\|
        \int_{1-t}^{1-t/2}
        \big\langle
        \xi^{v_1, v_2}(1-t,x ;r), 
        \mathrm{\, d} W_r
        \big\rangle
        \bigg\|_{L^2(\Omega, \mathbb{R})}
        \big\|
        D u
        (t / 2, X_{1-t / 2}^{1-t,x})
        \eta^{v_3}(1-t,x ;1-t/2)
        \big\|
        _{L^2(\Omega, \mathbb{R})}
        \nonumber\\
        & \leq
        C
        \frac{
        |\phi|_0}
        {t^{3/2}}
        \left(
        \int_{1-t}^{1-t/2}
        \left\|
        \xi^{v_1,v_2}(1-t,x ;r)
        \right\|
        ^2_{L^2(\Omega, \mathbb{R}^d)} 
        \mathrm{\, d} r
        \right)^{1 / 2}
        \left\|
        \eta^{v_3}(1-t,x ;1-t/2)
        \right\| \nonumber\\
        & \leq
        C
        \frac{
        |\phi|_0}
        {t^{3/2}}
        (e^{(3L_f+1/2)t}-1)^{1 / 2}
        e^{L_f t/2}
        \|v_1\|\cdot \|v_2\|\cdot
        \|v_3\|
        \notag \\
        & \leq
        \frac{C |\phi|_0 }
        {t}
        \|v_1\|\cdot \|v_2\|\cdot
        \|v_3\|.
\end{align}
Thanks to Lemma \ref{lem:X-mean-square-differential} and \eqref{eq:D2-u-est}, by using the H\"older inequality get that
\begin{align}
\label{eq:I4-est}
        |I_4|
        & \leq
        \frac{2}{\sqrt{\beta}\,t}
        \bigg\|
        \int_{1-t}^{1-t/2}
        \big\langle
        \eta^{v_1}(1-t,x ;r), 
        \mathrm{\, d} W_r
        \big\rangle
        \bigg\|
        _{L^2(\Omega, \mathbb{R})}\\
        & \quad \cdot
        \big\|
        D^2 u
        (t / 2, X_{1-t / 2}^{1-t,x})
        (\eta^{v_3}(1-t,x ;1-t/2),
        \eta^{v_2}(1-t,x ;1-t/2))
        \big\|
        _{L^2(\Omega, \mathbb{R})} 
        \nonumber\\
        & \leq
        C
        \frac{
        |\phi|_0}
        { t^2 }
        \|\eta^{v_2}(1-t,x ;1-t/2)\| 
        \cdot
        \|
        \eta^{v_3}(1-t,x ;1-t/2)
        \|
        \left(
        \int_{1-t}^{1-t/2}
        \left\|
        \eta^{v_1}(1-t,x ;r)
        \right\|^2_{L^2(\Omega, \mathbb{R}^d)} 
        \mathrm{\, d} r
        \right)^{1 / 2} 
        \nonumber\\
        & \leq 
        C
        \frac{
        |\phi|_0}
        { t^2}
        e^{L_f t} 
        \big(
        e^{L_f t} -1
        \big)^{1/2}
        \|v_1\|\cdot
        \|v_2\|
        \cdot\|v_3\|
        \notag \\
        & \leq 
        \frac{C |\phi|_0 }
        { t^{3/2}}
        \|v_1\|\cdot
        \|v_2\|
        \cdot\|v_3\|.
\end{align}
%Here we choose the function $\Phi=\phi$ with $\mathcal{K}(\Phi, x)=|\phi|_0$. 
%
Combining \eqref{eq:I1-est}, \eqref{eq:I2-est} and \eqref{eq:I4-est} yields, for any $x, v_1, v_2, v_3 \in \mathbb{R}^d$,
\begin{equation}
\label{eq:D^3-u(t, x)-regularity-estimation}
    |D^3 u(t, x)(v_1, v_2,v_3)| 
    \leq 
    \frac{C |\phi|_0 }{t^{3/2}}
    \|v_1\| \cdot\|v_2\|\cdot\|v_3\|,
    %\quad \forall x, v_1, v_2, v_3 \in \mathbb{R}^d.
\end{equation}
where the constant $C$ depends only on $L_f$ and the fact that $1/\sqrt{t} \leq 1/t \leq 1/{t^{3/2}}, \, t \in (0,1]$.
Thus, the proof is completed.

\end{proof}

\section{Proof of Proposition \ref{prop-mixing-time}}
\label{Appendix:mixing-time}
\textbf{Proof.}
Given an error tolerance $\epsilon \in(0,1)$, Theorem \ref{thm-main-rerult-w0-w1} tells that, for $M$ being large enough and $h$ being small enough such that
\begin{equation}
\label{eq:mixing-time-formula}
    C d h |\mathrm{ln}h| 
    \leq 
    \frac{\epsilon}{2},
    \quad 
    \frac{Cd}{M}
    \leq 
    \frac{\epsilon}{2},
\end{equation}
one can arrive at
\[
\big\|
    \mathcal{L}aw
    \big(
    \tilde{Y}_1
    \big)
    - 
    \mu
    \big\|_{\mathrm{TV}}
    \leq
    \epsilon.
\]
Rearranging the first inequality of \eqref{eq:mixing-time-formula} gives 
\begin{equation}
\label{eq:rearrange-mixing}
    \frac{\frac{1}{h}}
    {\mathrm{ln(\frac{1}{h})}}
    \geq
    \frac{2 C d}{\epsilon}.
\end{equation}
For \(\widehat C\geq e\), the condition
$x\geq 2\widehat C\ln\widehat C$
is sufficient to ensure that
$\frac{x}{\ln x}\geq \widehat C.$
For \(0<\widehat C<e\), the latter inequality holds for every \(x>1\).
Therefore, \eqref{eq:rearrange-mixing} can be satisfied if 
\begin{equation}
    N
    =
    \tfrac{1}{h}
    \geq 
    \tfrac{4 C d}{\epsilon}
    \cdot
    \ln \Big( 
    \tfrac{2Cd}{\epsilon}\Big)
    =
    \tilde{\mathcal{O}}
    \big(\tfrac{d}{\epsilon}\big).
\end{equation}
The second part of inequality \eqref{eq:mixing-time-formula} requires
\begin{align*}
    M 
    \geq 
    \frac{2Cd}{\epsilon}.
\end{align*}
Therefore, the number of samples of $\xi$ used at each time step in the Monte Carlo estimator $\tilde{f}_\beta$ is of order ${\mathcal{O}}\big(\tfrac{d}{\epsilon}\big)$, and the total Monte Carlo sample complexity is of order $\tilde{\mathcal{O}}\big(\tfrac{d^{2}}{\epsilon^2}\big)$. This completes the proof of the proposition.
\end{document}